\documentclass[10pt, a4paper]{article}
\usepackage{graphicx}
\usepackage{color}
\usepackage{amsmath}
\usepackage{amssymb}
\usepackage{amscd}
\IfFileExists{bbm.sty}{\usepackage{bbm}}{}
\usepackage{tikz}
\usepackage{pgfplots}
\usepackage{framed}
\IfFileExists{ulem.sty}{\usepackage[normalem]{ulem}}{}
\usetikzlibrary{positioning,arrows.meta,bending,calc}
\usepackage{hyperref}
\hypersetup{
    bookmarks=true,         
    colorlinks=true,       
    linkcolor=red,          
    citecolor=green,        
    filecolor=magenta,      
    urlcolor=cyan           
}
\IfFileExists{tcolorbox.sty}{\usepackage[most]{tcolorbox}}{}
\usepackage[utf8]{inputenc}
\usepackage{amsthm}
\IfFileExists{glossaries.sty}{\usepackage{glossaries}}{}
\IfFileExists{algorithm2e.sty}{\usepackage[linesnumbered,lined,boxed,commentsnumbered]{algorithm2e}}{}

\IfFileExists{marginnote.sty}{\usepackage{marginnote}}{}

\newtheorem{Theorem}{Theorem}
\newtheorem{Assumption}{Assumption}

\newtheorem{Lemma}{Lemma}
\newtheorem{Proposition}{Proposition}
\newtheorem{Corollary}{Corollary}
\newtheorem{Remark}{Remark}
\newtheorem{Example}{Example}

\newenvironment{Lemmabis}[1]
    {\renewcommand{\theLemma}{\ref{#1} (rewrite)}%
    \addtocounter{Lemma}{-1}%
    \begin{Lemma}}
    {\end{Lemma}}

\newenvironment{Examplebis}[1]
    {\renewcommand{\theExample}{\ref{#1} (revisited)}%
    \addtocounter{Example}{-1}%
    \begin{Example}}
    {\end{Example}}

\newcounter{BigConst}                     

\newcounter{SmallConst}                     

\newcounter{gamma}                     

\newcounter{kappa}                     

\newcounter{delta}                     

\newcounter{theta}                     

\newcounter{L}                     

\newcounter{eps}                     

\newcommand{\vertiii}[1]{{\left\vert\kern-0.25ex\left\vert\kern-0.25ex\left\vert #1 
    \right\vert\kern-0.25ex\right\vert\kern-0.25ex\right\vert}}

\newcommand{\vertiiii}[1]{{\left\vert\kern-0.25ex\left\vert\kern-0.25ex\left\vert\kern-0.25ex\left\vert #1 
    \right\vert\kern-0.25ex\right\vert\kern-0.25ex\right\vert\kern-0.25ex\right\vert}}

\DeclareMathOperator*{\argmin}{argmin}

\DeclareMathOperator{\diag}{diag}
\DeclareMathOperator{\Span}{span}

\newcommand{\dist}{\mathrm{dist}}

\DeclareMathOperator*{\sign}{sign}

\newcommand{\Proj}{{\mathrm{Proj}}}
\newcommand{\Prox}{{\mathrm{Prox}}}

\newcommand{\self}{\mathrm{self}}
\newcommand{\FSD}{\mathrm{FSD}}
\def\ds1{\textrm{1\kern-0.25emI}} 

\DeclareMathOperator{\Var}{\mathrm{Var}}
\DeclareMathOperator{\VAR}{\mathrm{VAR}}
\DeclareMathOperator{\Bias}{\mathrm{Bias}}
\DeclareMathOperator{\BIAS}{\mathrm{BIAS}}

\newcommand \cN{{\cal N}}

\newcommand \cT{{\cal T}}

\newcommand \bE{{\mathbb E}}

\newcommand \bN{{\mathbb N}}

\newcommand \bP{{\mathbb P}}

\newcommand \bR{{\mathbb R}}

\newcommand{\bbeta}{{\boldsymbol{\beta}}}

\newcommand{\ve}{{\boldsymbol{e}}}

\newcommand{\vg}{{\boldsymbol{g}}}
\newcommand{\vh}{{\boldsymbol{h}}}

\newcommand{\vu}{{\boldsymbol{u}}}

\newcommand{\vw}{{\boldsymbol{w}}}
\newcommand{\vx}{{\boldsymbol{x}}}
\newcommand{\vy}{{\boldsymbol{y}}}
\newcommand{\vz}{{\boldsymbol{z}}}
\newcommand{\one}{{\boldsymbol{1}}}

\newcommand{\vG}{{\boldsymbol{G}}}

\newcommand{\vgamma}{{\boldsymbol{\gamma}}}

\newcommand{\vzero}{\boldsymbol{0}}

\newcommand{\bxi}{{\boldsymbol{\xi}}}

\DeclareMathOperator*{\Med}{Med}

\newcommand{\mysymbol}[3]{%
\newglossaryentry{#1}{%
      name={\ensuremath{#2}},%
      text={\ensuremath{#2}},%
      description={#3},%
      sort={#1}%
    }%
\expandafter\newcommand\expandafter{\csname smb#1\endcsname}{\gls{#1}}%
\expandafter\newcommand\expandafter{\csname #1\endcsname}{\ensuremath{#2}}%
}

\newcommand{\Mo}{\operatorname{Mo}}
\usepackage[english]{babel}
\usepackage{authblk}
\usepackage{enumitem}
\numberwithin{equation}{section}
\providecommand{\uline}[1]{\underline{#1}}
\allowdisplaybreaks[1]
\begin{document}
\title{On the Instance-Specific Sharp Rates for Constrained Least Squares Estimators in the Sequence Model}
\author{Zong Shang\thanks{School of Mathematics, Georgia Institute of Technology, Atlanta, GA, 30332-0160, USA. E-mail: \texttt{shangzong9@gmail.com}}
        }
\date{\today}
\maketitle
\begin{abstract}
We develop a general methodology based on the Feature Space Decomposition framework for establishing explicit, instance-specific sharp rates for the expected estimation error of convex-constrained least squares estimators in the sequence model, without requiring estimates of localized Gaussian mean widths. We apply this methodology to $\ell_q$ ellipsoids, $\ell_1$ balls, and $kB_1^N\cap B_\infty^N$, obtaining matching upper and lower bounds in the nondegenerate regime. Our methodology also extends beyond Gaussian noise to symmetric noise distributions satisfying a Poincaré inequality. We also provide an elementary proof of the instance-specific sharp rate for isotonic regression for any well-specified signal and any noise level. Combining this result with Han's asymptotic reduction, we further obtain explicit instance-specific sharp convergence rates in probability for isotonic regression in Gaussian linear regression under general signal-dependent asymptotics conditions.
\end{abstract}



\section{Introduction}\label{sec:introduction}

Let $N\in\bN_+$, let $T\subset \bR^N$ be a nonempty closed convex set, and let $\bbeta^\star\in T$ be a deterministic vector, called the signal. Let $\sigma_\xi>0$ and let $\bxi\sim\cN(\vzero,\sigma_\xi^2 I_N)$ be a centered Gaussian random vector. Let $\vy = \bbeta^\star + \bxi$. Equip $\bR^N$ with the Euclidean inner product $\langle\cdot,\cdot\rangle$ and the Euclidean norm $\|\cdot\|_2$. Our objective is to construct $r_\star^2:=r_\star^2(T,\bbeta^\star,\sigma_\xi)$ such that
\begin{align}
\label{eq:def_LSE}
\bE\| \hat\bbeta - \bbeta^\star\|_2^2 \asymp r_\star^2,\quad \mbox{ where } \quad \hat\bbeta\in\argmin\Bigl(\frac{1}{2}\|\vy-\bbeta\|_2^2:\, \bbeta\in T\Bigr) = \Proj_T[\vy]
\end{align}
is the constrained least squares estimator (LSE), where $\Proj_T[\cdot]$ is the Euclidean metric projection onto $T$. Here $\bR^N$ is called the feature space; this model is called the Gaussian sequence model; and $\mathbb{E}\|\hat\bbeta - \bbeta^\star\|_2^2$ is called the (squared risk) expected estimation error of the constrained LSE for the signal $\bbeta^\star$. For any $\boldsymbol z\in\mathbb R^N$, write $\hat{\boldsymbol\beta}[\boldsymbol z]=\Proj_T[\boldsymbol z]$, so that $\hat{\boldsymbol\beta}=\hat{\boldsymbol\beta}[\boldsymbol y]$.

Understanding the estimation error of constrained least squares estimators in the sequence model is a classical problem in mathematical statistics; see, e.g., \cite{zhang_risk_2002,chatterjee_risk_2015,chatterjeeNewPerspectiveLeast2014,bellecSharpOracleInequalities2018,fang_risk_2019,chatterjeeAdaptiveRiskBounds2019,neykov_minimax_2022,hanNoisyLinearInverse2023,prasadan_facts_2025}. A fundamental contribution in this direction is the work of Chatterjee \cite{chatterjeeNewPerspectiveLeast2014}, which establishes matching upper and lower bounds for the estimation error, namely, an \uline{instance-specific sharp rate}\footnote{Since no normalization is performed, this rate usually does not converge as $N\to\infty$, and hence we only call it a sharp rate instead of a sharp convergence rate.}, over a general closed convex set through the maximizer of an offset localized Gaussian width. 
\begin{Theorem}[\cite{chatterjeeNewPerspectiveLeast2014}]\label{theo:chatterjee}
    For Gaussian noise $\bxi\sim\mathcal N(\vzero,\sigma_\xi^2I_N)$, let $t_{\boldsymbol\beta^\star}=\operatorname*{argmax}_{t\geq0}\bigl\{\mathbb E\sup_{\boldsymbol\beta\in T:\,\|\boldsymbol\beta-\boldsymbol\beta^\star\|_2\leq t}\langle\boldsymbol\xi,\boldsymbol\beta-\boldsymbol\beta^\star\rangle-\frac{t^2}{2}\bigr\}.$
    There exist absolute constants $c_0,c_1,C_1>0$ such that, whenever $t_{\boldsymbol\beta^\star}\geq c_0\sigma_\xi$, $c_1t_{\boldsymbol\beta^\star}^2\leq\mathbb E\|\hat{\boldsymbol\beta}-\boldsymbol\beta^\star\|_2^2\leq C_1t_{\boldsymbol\beta^\star}^2.$
\end{Theorem}
This characterization establishes a remarkably general connection between convex-constrained LSE in the Gaussian sequence model and Gaussian processes. However, extracting an explicit rate from this variational quantity requires a precise analysis of the localized geometry of the constraint set around $\bbeta^\star$. We give an example below.
\begin{Example}\label{example:ell_q}
    Let $1\le q\le\infty$, and equip $\mathbb{R}^N$ with the canonical basis $\ve_1,\cdots,\ve_N$. Let $B_q^N$ be the unit ball of the $\|\cdot\|_q$ norm. Given $\Sigma=\diag(\sigma_1,\cdots,\sigma_N)$, a diagonal matrix with $\sigma_1,\cdots,\sigma_N>0$, let $T=\Sigma^{1/2}B_q^N$. Then the localized Gaussian mean width in Theorem~\ref{theo:chatterjee} is $\mathbb{E}\sup(\langle\bxi,\bbeta-\bbeta^\star\rangle-\frac{t^2}{2} : \bbeta\in\Sigma^{1/2}B_q^N\cap B_2(\bbeta^\star;t))$, where $B_2(\bbeta^\star;r)$ is the Euclidean ball centered at $\bbeta^\star$ with radius $r$.
\end{Example}
We emphasize that to apply Theorem~\ref{theo:chatterjee} to obtain the instance-specific sharp rate of the constrained LSE over the index set in Example~\ref{example:ell_q}, one needs to estimate this offset localized Gaussian mean width up to a multiplicative constant. However, even without the offset term, estimating the localized Gaussian mean width of such an anisotropic convex body is a well-known challenge, let alone after adding the offset term $\frac{t^2}{2}$. For example, when $q=2$, such an estimate is already quite complicated; see \cite{wei_from_gaussian_2020}. The contributions of this work are as follows:
\begin{framed}
    To provide a methodology by which one can obtain an explicit instance-specific sharp rate.
\end{framed}
In contrast to Chatterjee \cite{chatterjeeNewPerspectiveLeast2014}, we do not aim to establish a general master theorem applicable to arbitrary convex constraint sets. Although the underlying principles of our analysis are common to different problems, their implementation depends on the geometry of the constraint set. A master theorem encompassing these different geometries would require additional technical assumptions, making its formulation considerably more complicated and potentially obscuring the main mathematical mechanisms of the method. We therefore organize our presentation around concrete examples. We first develop the complete argument for the $\ell_q$ ellipsoid, including the proofs of matching upper and lower bounds, and subsequently apply the method to other convex constraint sets, emphasizing the geometric estimates and modifications specific to each setting. Through this approach, we seek to reveal the common mathematical principles underlying these problems and illustrate how our method can be implemented to establish explicit instance-specific sharp rates in the sequence model.

\subsection{Structure of this work}

The remainder of this work is organized as follows. Section~\ref{sec:main_results} presents the main results for the $\ell_q$ ellipsoid, the $\ell_1$ ball $\rho B_1^N$, the $k$-interpolated body $kB_1^N\cap B_\infty^N$, and isotonic regression, together with a discussion of their statistical implications. The non-asymptotic, explicit, user-friendly instance-specific rate for isotonic regression is placed in Section~\ref{sec:isotonic_regression}, while in Section~\ref{sec:asymptotic_isotonic}, combining Han \cite{hanNoisyLinearInverse2023}, we give the asymptotic, explicit, user-friendly instance-specific convergence rate for Gaussian linear regression. In Section~\ref{sec:exemplar_ell_q}, we take the $\ell_q$ ellipsoid as an exemplar and develop the complete analysis, including the proofs of matching upper and lower bounds. The remaining proofs are placed in the appendix.

\subsection{Notation}

Let $\{\ve_j\}_{j=1}^N$ be the canonical basis of $\bR^N$, and for any $\bbeta\in\bR^N$ and any $j\in[N]$, we write $\beta_j=\langle\bbeta,\ve_j\rangle$. We use ``$:=$'' to mean ``defined as''. For any integer $k\in\bN_+$, we write $[k]:=\{1,\cdots,k\}$. For any real numbers $a,b>0$, we say $a\lesssim b$ if there exists an absolute constant $C>0$ such that $a\leq Cb$; similarly we define $a\gtrsim b$. If $a\lesssim b$ and $a\gtrsim b$, we write $a\asymp b$. For $q\in\bR_+$, we write $c=c(q)$ to represent a constant depending on $q$. If \(a\leq C(q)b\), we write \(a\lesssim_q b\). Similarly, we define \(\gtrsim_q\) and \(\asymp_q\). We use $\partial^-g(\bullet)$ to denote the subdifferential of $g$ at $\bullet$, and $\nabla g(\bullet)$ to denote its gradient at $\bullet$. For a multivariate function $g(\bullet_1,\bullet_2)$, we use $\nabla_{\bullet_1}g(\vx,\vy)$ to denote the gradient of $g(\bullet_1,\vy)$ at $\bullet_1=\vx$ with $\bullet_2=\vy$ fixed.  For any $\vu\in\mathbb{R}^N$ and any set $A\subset\mathbb{R}^N$, denote by $\dist(\vu,A)=\inf\{\|\vu-\vx\|_2:\vx\in A\}=\|\vu-\Proj_T[\vu]\|_2$ the Euclidean distance from $\vu$ to the set $A$. For a proper, lower semicontinuous convex function $g:\bR^N\to\bR\cup\{+\infty\}$, let $\Prox_g[\vx]:=\argmin_{\vu\in\bR^N}\{\frac12\|\vx-\vu\|_2^2+g(\vu)\}$ be its proximal operator, \cite[Definition 12.23]{bauschkeConvexAnalysisMonotone2011}. 

\section{Main Results}\label{sec:main_results}
The contributions of this work are twofold: on the one hand, we propose a method for establishing instance-specific sharp rates for the estimation error of the convex-constrained LSE. We apply this method to several examples to provide exemplars of the method, and the resulting statistical conclusions constitute the second aspect of the contributions of this work. In this section, we first introduce the statistical conclusions obtained by applying this method to the \(\ell_q\) ellipsoid, the \(\ell_1\) ball, and the \(k\)-interpolated body. To the best of our knowledge, these statistical conclusions are novel. The results on isotonic regression in the Gaussian sequence model and Gaussian linear regression are postponed in Section~\ref{sec:isotonic_regression} and Section~\ref{sec:asymptotic_isotonic} respectively. Further discussion is postponed to Section~\ref{sec:comments}. The proof of Theorem~\ref{theo:ell_q} can be found in Section~\ref{sec:exemplar_ell_q}; for Theorem~\ref{theo:ell_1_ball} and Theorem~\ref{theo:k_interpolated_body}, we only briefly describe how their proofs differ from that of Theorem~\ref{theo:ell_q}.

\paragraph{\texorpdfstring{$\ell_q$}{lq} ellipsoid.}
Let $T$ be $\Sigma^{1/2}B_q^N$ considered in Example~\ref{example:ell_q}, where $1<q<\infty$.  For $(\lambda,\vx)\in\mathbb R_+\times\mathbb R^N$, define
\begin{align*}
 \Mo(\lambda,\vx):=\min_{\bbeta\in\mathbb R^N}
 \Bigl\{\frac12\|\vx-\bbeta\|_2^2
       +\lambda\Psi(\bbeta)\Bigr\},\quad\mbox{ where }\quad \Psi(\bbeta) =  \frac{1}{q}\Bigl(\|\Sigma^{-1/2}\bbeta\|_q^q-1\Bigr).
\end{align*}
For any $\lambda\geq0$, $\vx\mapsto \Mo(\lambda,\vx)$ is the Moreau envelope of the function $\lambda\Psi$, \cite[Section 12.4]{bauschkeConvexAnalysisMonotone2011}. For any $\vx\in\bR^N$, let
\begin{equation}\label{eq:def_divergence}
 d_\lambda(\vx):=\frac12\dist^2(\vx,T)-\Mo(\lambda,\vx)+\lambda\Psi\bigl(\hat\bbeta[\vx]\bigr)\in\bR_+,\mbox{ and }
 \lambda_\self:=\argmin_{\lambda\geq0}\bE d_\lambda(\vy).
\end{equation}
We explain the non-negativity of $d_\lambda(\cdot)$ later in Proposition~\ref{prop:implicit_regularization}. We call $\lambda_\self$ the \uline{self-regularization strength}; its existence and uniqueness are proved in Lemma~\ref{lem:calibration_equation_l_q}.  Below we present the main results for the $\ell_q$ ellipsoid case, whose proof may be found in Section~\ref{sec:exemplar_ell_q}.
\begin{Theorem}\label{theo:ell_q}
For any $N\in\bN_+$, any $1<q<\infty$, any $\sigma_\xi>0$ and any $\sigma_1,\cdots,\sigma_N>0$, $\lambda_\self$ defined in \eqref{eq:def_divergence} is the unique solution of the following equation
\begin{align}\label{eq:general_self_calibration}
    \sum_{j=1}^N\sigma_j^{-q/2}\bE|\eta_{\self,j}(y_j)|^q=\bE\min\Bigl\{1,\sum_{j=1}^N\sigma_j^{-q/2}|y_j|^q\Bigr\},\mbox{ where }\eta_{\self,j}(\cdot) = \argmin\Bigl(\frac{1}{2}(\cdot-s)^2+\frac{\lambda_\self}{q}\sigma_j^{-q/2}|s|^q: s\in\bR\Bigr).
\end{align}Let
\begin{align}\label{eq:explicit_ellipsoid_rate}
    r_\star^2 = \sum_{j=1}^N\min\left\{|\beta_j^\star|^2, \frac{\lambda_\self^2}{\sigma_j^q}|\beta_j^\star|^{2q-2}\right\} +\sum_{j=1}^N\min\left\{\sigma_\xi^2, \left(\frac{\sigma_j^{q/2}\sigma_\xi}{\lambda_\self}\right)^{2/(q-1)}\right\}.
\end{align}
    There exist absolute constants $c=c(q),C=C(q)$ and $C'=C'(q)>0$ such that if $r_\star\geq C'\sigma_\xi$, then $cr_\star^2\leq \bE\|\hat\bbeta-\bbeta^\star\|_2^2\leq Cr_\star^2.$
\end{Theorem}

\paragraph{\texorpdfstring{$\ell_1$}{l1} ball.} We consider the case when $T=\rho B_1^N$ for some $\rho>0$.
Replacing $\Psi(\bbeta) = \frac{1}{q}(\|\Sigma^{-1/2}\bbeta\|_q^q-1)$ in the $\ell_q$ ellipsoid case by $\Psi(\bbeta) = \|\bbeta\|_1-\rho$, we obtain
\begin{align*}
    \Mo(\lambda,\vx)&=\min_{\bbeta\in\bR^N}\Bigl(\frac{1}{2}\|\vx-\bbeta\|_2^2+\lambda(\|\bbeta\|_1-\rho)\Bigr),\mbox{ and } d_\lambda(\vx)=\frac{1}{2}\dist^2(\vx,T)-\Mo(\lambda,\vx)+\lambda\bigl(\|\hat\bbeta[\vx]\|_1-\rho\bigr).
\end{align*}
Recall that $\lambda_\self=\argmin(\bE[d_\lambda(\vy)]:\lambda\geq0)$. For every $\lambda\geq0$, $1\leq j\leq N$, and $t\in\bR$, let $\eta_{\lambda,j}(t)=\argmin(\frac{1}{2}(t-s)^2+\lambda|s|:s\in\bR)=\sign(t)(|t|-\lambda)_+$, that is, the soft-thresholding function. Here, since \(\rho B_1^N\) is isotropic, \(\eta_{\lambda,j}\) is the same for all \(j\). Repeating the proof for the $\ell_q$ ellipsoid, we obtain the following result.
\begin{Theorem}\label{theo:ell_1_ball}
For any $\rho>0$, $N\in\bN_+$, and $\sigma_\xi>0$, $\lambda_\self\geq0$ is the unique solution of the following equation
\begin{align}\label{eq:eta_self_ell_1_ball}
    \sum_{j=1}^N\bE\bigl(|y_j|-\lambda_\self\bigr)_+=\bE\min\{\rho,\|\vy\|_1\}.
\end{align}
Let
\begin{align*}
    r_\star^2=\sum_{j=1}^N\min\Bigl\{\min\bigl\{\bigl|\beta_j^\star\bigr|^2,\lambda_\self^2\bigr\}+\sigma_\xi^2,\, \bigl|\beta_j^\star\bigr|^2+\bE\bigl(|\xi_j|-\lambda_\self\bigr)_+^2\Bigr\}.
\end{align*}
Then there exist absolute constants $c,C,C'>0$ such that if $r_\star\geq C'\sigma_\xi$, then $cr_\star^2\leq \bE\|\hat\bbeta-\bbeta^\star\|_2^2\leq Cr_\star^2$.
\end{Theorem}

\paragraph{\texorpdfstring{$k$}{k}-interpolated body.}
We can also change the representation of \(T\) to accommodate more constraints. For example, we consider \(T=kB_{1}^{N}\cap B_{\infty}^{N}\), where \(1\leq k\leq N\). This is a typical class of convex bodies in Asymptotic Geometric Analysis (\cite{gordonGeometrySpacesPolytopes2002}), and we take it as an example here because it provides an example in which \(T\) has two types of constraints. Here, we replace $\Psi(\bbeta) = \frac{1}{q}\bigl(\|\Sigma^{-1/2}\bbeta\|_q^q-1\bigr)$ in the $\ell_q$ ellipsoid case by $\Psi(\bbeta) = \|\bbeta\|_1-k$, and we write $T=\{\bbeta\in B_\infty^N: \Psi(\bbeta)\leq 0\}$. Define
\begin{align*}
    \Mo(\lambda,\vx)&=\min_{\bbeta\in B_\infty^N}\Bigl(\frac{1}{2}\|\vx-\bbeta\|_2^2+\lambda(\|\bbeta\|_1-k)\Bigr),\mbox{ and }
    d_\lambda(\vx)=\frac{1}{2}\dist^2(\vx,T)-\Mo(\lambda,\vx)+\lambda\bigl(\|\hat\bbeta[\vx]\|_1-k\bigr).
\end{align*}Recall that $\lambda_\self=\argmin_{\lambda\geq0}\bE d_\lambda(\vy).$ For every $\lambda\geq0$ and $t\in\bR$, let $\eta_{\lambda,j}(t)=\argmin(\frac{1}{2}(t-s)^2+\lambda|s|:|s|\leq1)=\sign(t)\min\{1,(|t|-\lambda)_+\}$. That is, the only difference from the \(\ell_q\) ellipsoid and the \(\ell_1\) ball is that we change the feasible set of \(\Mo\) from \(\bR^N\) to \(B_\infty^N\). Nevertheless, our analysis still holds. Repeating the proof for the $\ell_q$ ellipsoid, we obtain the following result.
\begin{Theorem}\label{theo:k_interpolated_body}
Let $1\leq k\leq N$ and $T=kB_1^N\cap B_\infty^N$. For every $\sigma_\xi>0$, $\lambda_\self$ exists uniquely and satisfies
\begin{align}\label{eq:eta_k_interpolated_body}
    \sum_{j=1}^N\bE\min\bigl\{1,(|y_j|-\lambda_\self)_+\bigr\}=\bE\min\Bigl\{k,\sum_{j=1}^N\min\{1,|y_j|\}\Bigr\}.
\end{align}
Let
\begin{align*}
    r_\star^2=\sum_{j=1}^N\min\Bigl\{\min\bigl\{|\beta_j^\star|^2,\lambda_\self^2\bigr\} + \sigma_\xi^2,\, \bigl|\beta_j^\star\bigr|^2+\bE\bigl[\min\bigl\{1,\bigl(|\xi_j|-\lambda_\self\bigr)_+^2\bigr\}\bigr]\Bigr\}.
\end{align*}
Then there exist absolute constants $c,C,C'>0$ such that if $r_\star\geq C'\sigma_\xi$, then $cr_\star^2\leq\bE\|\hat\bbeta-\bbeta^\star\|_2^2\leq Cr_\star^2$.
\end{Theorem}

\subsection{Examples}\label{sec:examples_LSE}

This subsection provides some concrete examples to show how \(\lambda_\self\) and \(r_\star^2\) can be explicitly estimated. The computation of the following examples may be found in Section~\ref{sec:computation_LSE_examples}.

\begin{Example}\label{example:ellipsoid_rate}
    Suppose $1<q<\infty$, $\alpha>2-\frac{2}{q}$. Let $\sigma_j = j^{-\alpha}$, $T=\Sigma^{1/2}B_q^N$, and $\bxi\sim\cN(\vzero,\sigma_\xi^2 I_N)$. Let $0<\delta<1$ be some parameter. There exist constants $c=c(\alpha,q,\delta),C=C(\alpha,q,\delta)$, such that for any $0<\sigma_\xi\leq c$, $N\geq C\sigma_\xi^{-\frac{2q}{\alpha q+2}}$, and any $\sum_{j=1}^N j^{\frac{\alpha q}{2}}\bigl|\beta_j^\star\bigr|^q\leq 1-\delta$, there holds $\lambda_\self\asymp r_\star^2\asymp\bE\|\hat\bbeta-\bbeta^\star\|_2^2\asymp \bigl(\sigma_\xi^2\bigr)^{\frac{\alpha q+2-q}{\alpha q+2}}$.
\end{Example}
This example obtains an explicit representation of the sharp rate for a polynomially decaying spectrum under the source condition \(\sum_{j=1}^N j^{\alpha q/2}|\beta_j^\star|<\infty\). The difference between this example and \cite{wei_from_gaussian_2020} is that we allow a general \(q\), and we do not assume that \(\boldsymbol\beta^\star=\vzero\) or that it is close to some coordinate axis. We do not intend to recover the conclusions of \cite{wei_from_gaussian_2020} in these cases---although this can be done.
\begin{Example}\label{example:k_interpolated}
    Let $1\leq k\leq N$, $T=kB_1^N\cap B_\infty^N$, $\bbeta^\star=\vzero$. We say $\xi$ is Gaussian if $\bxi\sim\cN(\vzero,\sigma_\xi^2 I_N)$; and is exponential, if $\bxi$ has i.i.d. coordinate $\xi_j$ that has probability density function $t\mapsto \frac{1}{\sqrt{2}\sigma_\xi}\exp(-\frac{\sqrt{2}}{\sigma_\xi}|t|)$. Then for any $\sigma_\xi>0$, 
    \begin{align*}
        r_\star^2\asymp \min\bigl\{k,N\min\{\sigma_\xi,1\}\bigr\}\begin{cases}
            \min\Bigl\{1,\frac{\sigma_\xi}{\sqrt{\log\bigl(e+\frac{N\min\{\sigma_\xi,1\}}{k}\bigr)}}\Bigr\}, & \mbox{ if } \xi \mbox{ is Gaussian},\\
            \min\{\sigma_\xi, 1\}, & \mbox{ if }\xi\mbox{ is exponential}.
        \end{cases}
    \end{align*}
\end{Example}
To the best of our knowledge, this example has not appeared in mathematical statistics, and we include it only to illustrate that Theorem~\ref{theo:k_interpolated_body} can still yield explicit instance-specific sharp rates for a more complicated index set and for probability measures beyond Gaussian.

\subsection{Comments}\label{sec:comments}

\begin{Remark}\label{rmk:poincare_noise}
    This remark considers the case where \(\bxi\) is not necessarily a Gaussian random vector. Suppose $\bxi$ is symmetric, $\bE[\bxi\otimes\bxi]=\sigma_\xi^2 I_N$; and satisfies Poincar\'e's inequality (see, for instance, \cite[Definition 1.1.1]{artstein-avidanAsymptoticGeometricAnalysis2021}) with some parameter $C_{\mathrm{P}}$, that is, for any $1$-Lipschitz function $h:\bR^N\to\bR$, $\Var(h(\bxi))\leq C_{\mathrm{P}}\sigma_\xi^2\bE\|\nabla h(\bxi)\|_2^2$. Then Theorem~\ref{theo:ell_q}, Theorem~\ref{theo:ell_1_ball} and Theorem~\ref{theo:k_interpolated_body} hold with constants depending further on $C_{\mathrm{P}}$.
\end{Remark}

\begin{Remark}
    In Theorem~\ref{theo:ell_1_ball}, one may show that
    \begin{align*}
        r_\star^2 \asymp \sum_{j=1}^N\min\bigl\{\bigl|\beta_j^\star\bigr|^2,\lambda_\self^2\bigr\} + \frac{N\sigma_\xi^2\exp\bigl(-\frac{\lambda_\self^2}{2\sigma_\xi^2}\bigr)}{\bigl(1+\frac{\lambda_\self}{\sigma_\xi}\bigr)^3},
    \end{align*}see Lemma~\ref{lem:ell_1_equivalence}.
\end{Remark}
Theorem~\ref{theo:ell_q}, Theorem~\ref{theo:ell_1_ball}, and Theorem~\ref{theo:k_interpolated_body} assume \(r_\star^2\gtrsim\sigma_\xi^2\), which we call the ``non-degenerate regime''---the noise level \(\sigma_\xi\) of interest for the general constrained LSE is typically a constant, whereas \(r_\star\) is usually a sequence growing to infinity with \(N\), and hence we do not consider \(\sigma_\xi\gtrsim r_\star\) to be the most interesting regime in this model, which is also the same regime that the theorem (see Theorem~\ref{theo:chatterjee}) of Chatterjee \cite{chatterjeeNewPerspectiveLeast2014} can support. Combining Theorem~\ref{theo:ell_q}, Theorem~\ref{theo:ell_1_ball}, and Theorem~\ref{theo:k_interpolated_body}, we claim that we have developed a mechanism for obtaining explicit instance-specific sharp rates for the expected estimation error of the convex-constrained LSE in the sequence model. We emphasize once again that we do not pursue a master theorem valid for general convex bodies as in Chatterjee \cite{chatterjeeNewPerspectiveLeast2014}---in fact, we do not claim that this approach can be directly applied to arbitrary convex bodies. In fact, depending on the index set under consideration, this approach often needs to be modified. What we emphasize is that the proof strategy itself introduced below may be universal.

Below we provide geometric perspectives on Theorem~\ref{theo:ell_q}, Theorem~\ref{theo:ell_1_ball}, and Theorem~\ref{theo:k_interpolated_body}, and explain the general proof strategy behind these theorems.  This is another aspect of the contributions of this work: the methodological perspective.

We first construct, via the KKT conditions, a transformation of the constrained LSE into an unconstrained RERM with a data-dependent regularization functional. For uniformity of exposition, we let $K\subset\mathbb{R}^N$ be a nonempty closed convex set such that $T = \{\bbeta\in K: \Psi(\bbeta)\leq 0\}$ for some convex function $\Psi:\bR^N\to\bR$. Then the KKT conditions for $\hat\bbeta=\arg\min(\frac{1}{2}\|\vy-\bbeta\|_2^2:\bbeta\in T)$ guarantee the existence of a Lagrange multiplier $\hat\lambda$ for $\Psi$ such that (see Proposition~\ref{prop:implicit_regularization})
\begin{align*}
    \hat\bbeta = \argmin\Bigl(\frac{1}{2}\|\vy-\bbeta\|_2^2+\hat\lambda\Psi(\bbeta):\bbeta\in K\Bigr).
\end{align*}
This is precisely the implicit regularization induced by the index set of the constrained LSE. 

The second step of the proof is to construct a surrogate estimator $\tilde\bbeta$, which replaces the data-dependent regularization functional $\hat\lambda\Psi$ by the deterministic regularization functional $\lambda_\self\Psi$, such that $\mathbb{E}\|\hat\bbeta-\tilde\bbeta\|_2^2\leq c r_\star^2$ for some $c\le1$. Here, the surrogate estimator is precisely the proximal operator of $\lambda_\self\Psi$, see Section~\ref{sec:convex_analysis}, and for the index sets and $\Psi$ we consider, it is well-behaved enough that the proximal operator is separable; that is,
\begin{align*}
    \tilde\bbeta=\argmin\Bigl(\frac{1}{2}\|\vy-\bbeta\|_2^2+\lambda_\self\Psi(\bbeta):\bbeta\in K\Bigr) = \Prox_{\lambda_\self \Psi+\iota_K}[\vy] = \sum_{j=1}^N\eta_{\self,j}(y_j)\ve_j,
\end{align*}where $\iota_K$ is the indicator function, $\vy=(y_j)_{j=1}^N$, and $\eta_{\self,j}(\cdot)$ is the coordinate proximal map associated with $\lambda_\self\Psi(\cdot_j)$ restricted on $K$. The purpose of the second step is to prove that the expected squared Euclidean distance between the surrogate estimator and $\hat\bbeta$ is negligible compared with $r_\star^2$; hence the analysis of $\mathbb{E}\|\hat\bbeta-\bbeta^\star\|_2^2$ can be reduced to that of the estimator $\tilde\bbeta$, which has a very well-behaved deterministic regularization functional (indeed, coordinate-separable). The investigation of the upper bound on \(\bE\|\hat\bbeta-\tilde\bbeta\|_2^2\) precisely illustrates the motivation for defining \(d_\lambda\) and \(\lambda_\self\) in this way. For any \(\lambda\geq0\), note that \(\bbeta\mapsto \frac{1}{2}\|\vy-\bbeta\|_2^2+\lambda\Psi(\bbeta)\) is a \(1\)-strongly convex function (see, for instance, \cite[Definition 10.5]{bauschkeConvexAnalysisMonotone2011}), and \(\Prox_{\lambda\Psi}[\vy]\) is precisely its minimizer; hence, by the definition of strong convexity, we have
\begin{align*}
    \forall\lambda\geq0,\qquad \frac{1}{2}\|\vy-\hat\bbeta\|_2^2 - \frac{1}{2}\|\vy-\Prox_{\lambda\Psi}[\vy]\|_2^2 + \lambda\Psi(\hat\bbeta) - \lambda\Psi(\Prox_{\lambda\Psi}[\vy]) \geq \frac{1}{2}\|\hat\bbeta-\Prox_{\lambda\Psi}[\vy]\|_2^2.
\end{align*}
Later in Lemma~\ref{lem:d_varepsilon} we prove that the left-hand side of the above inequality is exactly \(d_\lambda(\vy)\), that is, \(d_\lambda(\vy)\) characterizes an upper bound on the (square of the) Euclidean distance between \(\hat\bbeta\) and \(\Prox_{\lambda\Psi}[\vy]\); therefore, in order to make this expected approximation error as small as possible, we define \(\lambda_\self\) in \eqref{eq:def_divergence} as the minimizer of \(\bE[d_\lambda(\vy)]\).

The third step is to establish the instance-specific sharp rate of the estimation error of $\tilde\bbeta$, namely to prove that $\mathbb{E}\|\tilde\bbeta-\bbeta^\star\|_2^2\asymp r_\star^2$. Since the $\tilde\bbeta$ we consider here is coordinate-separable, this step becomes extremely simple and does not require the full strength of the Feature Space Decomposition (FSD) method. Even so, its estimation error can still naturally be endowed with geometric meaning under the FSD framework. Here, we refer interested readers to the author's PhD thesis for a comprehensive understanding of the FSD framework, \cite{shangFeatureSpaceDecomposition2026}.

In the squared risk setting of this work, the FSD framework states that the smallness of the estimation error or population excess risk has exactly two sources: 1. on the estimation subspace, the projected estimator estimates the projected signal; 2. on the free subspace, the sum of the squared Euclidean norms of the projected estimator and the projected signal is not too large. We first explain the terminology involved here.
Let us abbreviate $\eta_{\self,j}$ by $\eta_j$. Define
\begin{align}\label{eq:def_J_star}
    J_\star=\Bigl\{j\in[N]:\, \bigl|\beta_j^\star - \eta_j(\beta_j^\star)\bigr|^2 + \sigma_\xi^2 \leq \bigl|\beta_j^\star\bigr|^2 + \bE\bigl|\eta_j(\xi_j)\bigr|^2\Bigr\},\quad\mbox{ and }J_\star^c=[N]\backslash J_\star.
\end{align}Recall that $\{\ve_j\}_{j=1}^N$ is a canonical basis of the feature space $\mathbb{R}^N$. We set $V_{J_\star} := \Span\{\ve_j: j\in J_\star\}$ and $V_{J_\star^c} :=\Span\{\ve_j: j\in J_\star^c\}$, called the \uline{estimation subspace} and the \uline{free subspace}, respectively. We adopt the convention that for any $\bbeta\in\mathbb{R}^N$, $\bbeta_{J_\star}=\Proj_{V_{J_\star}}[\bbeta]$ and $\bbeta_{J_\star^c}=\Proj_{V_{J_\star^c}}[\bbeta]$. We say that estimation occurs on the estimation subspace, namely that $\bE\|\tilde\bbeta_{J_\star}-\bbeta_{J_\star}^\star\|_2^2$ is small; while on the free subspace, the squared Euclidean norms of the projections of the estimator and the signal are relatively small, namely $\|\bbeta_{J_\star^c}^\star\|_2^2 + \mathbb{E}\|\tilde\bbeta_{J_\star^c}\|_2^2$ is small. Below we explain the rationale for these two phenomena. We first define the two bias terms and two variance terms induced by the projections onto the two subspaces. Define
\begin{align*}
\Bias_{J_\star} &= \sum_{j\in J_\star}\bigl|\beta_j^\star - \eta_j(\beta_j^\star)\bigr|^2, \qquad& \Bias_{J_\star^c} &= \|\bbeta_{J_\star^c}^\star\|_2^2,\\ \Var_{J_\star} &= \sigma_\xi^2|J_\star|, \qquad\mbox{and}& \Var_{J_\star^c} &= \sum_{j\in J_\star^c}\bE\bigl|\eta_j(\xi_j)\bigr|^2.
\end{align*}
Let
\begin{align*}
    r^2(V_{J_\star},V_{J_\star^c})=\BIAS^\FSD+\VAR^\FSD,\quad\mbox{where}\quad \BIAS^\FSD=\Bias_{J_\star}+\Bias_{J_\star^c},\mbox{ and }\VAR^\FSD=\Var_{J_\star}+\Var_{J_\star^c}.
\end{align*}
Taking Theorem~\ref{theo:ell_q} as an example, we show that $r^2(V_{J_\star},V_{J_\star^c})$ characterizes $\mathbb{E}\|\tilde\bbeta-\bbeta^\star\|_2^2$. The remaining examples still satisfy similar equivalence relations.
\begin{Proposition}\label{prop:equivalence_r_square_star}
Under the assumptions of Theorem~\ref{theo:ell_q}, we have
    \begin{align*}
        \bE\bigl\|\tilde\bbeta-\bbeta^\star\bigr\|_2^2\asymp r^2(V_{J_\star},V_{J_\star^c}) = \sum_{j=1}^N\min\Bigl\{ \bigl|\beta_j^\star-\eta_j(\beta_j^\star)\bigr|^2+\sigma_\xi^2,\, \bigl|\beta_j^\star\bigr|^2+\bE\bigl|\eta_j(\xi_j)\bigr|^2 \Bigr\}.
    \end{align*}Moreover, there exist constants $c=c(q), C=C(q)$ such that $c(q)r^2(V_{J_\star},V_{J_\star^c})\leq r_\star^2\leq C(q)r^2(V_{J_\star},V_{J_\star^c})$.
\end{Proposition}
\begin{proof}[Proof of Proposition~\ref{prop:equivalence_r_square_star}]
    The claimed equality readily follows; while $\bE\|\tilde\bbeta-\bbeta^\star\|_2^2\asymp r^2(V_{J_\star},V_{J_\star^c})$ is a consequence of \eqref{eq:upper_tilde_bbeta_estimation_error} and of Proposition~\ref{prop:lower_surrogate_tilde_beta}. The ``moreover'' part is a consequence of Lemma~\ref{lem:ell_q_coordinate_estimates}.
\end{proof}
Let us now use Proposition~\ref{prop:equivalence_r_square_star} to explain the two phenomena of the expected estimation error of $\tilde\bbeta$ as explained by the FSD framework: the contribution of $V_{J_\star}$ to the estimation error of $\tilde\bbeta$ being small is due to $\bE\|\tilde\bbeta_{J_\star} - \bbeta_{J_\star}^\star\|_2^2=\sum_{j\in J_\star}\bE\bigl|\eta_j(\beta_j^\star+\xi_j)-\beta_j^\star\bigr|^2$ being small, which is manifested as $\Bias_{J_\star}$ and $\Var_{J_\star}$ being small; whereas the contribution on $V_{J_\star^c}$ is due to $\|\tilde\bbeta_{J_\star^c}\|_2^2 + \|\bbeta_{J_\star^c}^\star\|_2^2$ being small, rather than $\tilde\bbeta_{J_\star^c}$ estimating $\bbeta_{J_\star^c}^\star$. This can be seen from the fact that $\|\bbeta_{J_\star^c}^\star\|_2^2$ still appears in $r^2(V_{J_\star},V_{J_\star^c})\asymp \mathbb{E}\|\tilde\bbeta-\bbeta^\star\|_2^2$ — if $\bbeta_{J_\star^c}^\star$ were estimated, we should expect the bias on $V_{J_\star^c}$ no longer to be this part of the energy itself but rather in a form of cancellation (i.e., $\bigl|\beta_j^\star-\eta_j(\beta_j^\star)\bigr|^2$) that is closer to that in \(\Bias_{J_\star}\) — but the opposite is the case, and hence we say that this part of the signal is not estimated. Instead, the contribution on $V_{J_\star^c}$ is that the sum of $\mathbb{E}\|\tilde\bbeta_{J_\star^c}\|_2^2\leq 2\|\bbeta_{J_\star^c}^\star\|_2^2 + 2\Var_{J_\star^c}$ and $\|\bbeta_{J_\star^c}^\star\|_2^2$ is small.

Summarizing the above discussion, we know that at least for Theorem~\ref{theo:ell_q}, Theorem~\ref{theo:ell_1_ball}, and Theorem~\ref{theo:k_interpolated_body}, these constrained LSEs approximately possess an implicit- and self-regularization property, that is, although the definition of the constrained LSE itself does not explicitly incorporate any regularization, it implicitly performs regularization, making it ``prefer'' \(\tilde\bbeta\) in the sense of the expected estimation error; hence, we obtain the following intuition: the constrained LSE on \(\rho B_1^N\) is approximately a LASSO estimator with \(\lambda_\self\|\cdot\|_1\) as the regularization functional; the constrained LSE on the \(\ell_q\) ellipsoid is approximately an RERM with \(\frac{\lambda_\self}{q}\|\Sigma^{-1/2}\cdot\|_q^q\) as the regularization functional; and the \(k\)-interpolated one is approximately a ``constrained'' LASSO on \(B_\infty^N\) with \(\lambda_\self\|\cdot\|_1\) as the regularization functional. This is also why we call \(\lambda_\self\) the self-regularization strength, and the form of this regularization functional comes from the geometric structure of the index set of the constrained LSE.

\subsection{Isotonic Regression}\label{sec:isotonic_regression}
In this subsection, we record the explicit instance-specific sharp rate for univariate isotonic regression. This is another contribution of this paper in terms of statistical conclusions. Although this problem can still be analyzed by FSD, since its proof is longer than the elementary proof that we will introduce in this subsection, we adopt a more direct, elementary approach. Isotonic regression is one of the most typical examples of shape-constrained regression problems, \cite{ayer_empirical_1955,Grenander01071956,van_de_geer_estimating_1990,meyer_degrees_2000,durot_sharp_2002,durot_Lp_2007,cator_adaptivity_2011,chatterjeeNewPerspectiveLeast2014,chatterjee_risk_2015,luss_bounded_2017,bellecSharpOracleInequalities2018,soloff_distribution-free_2019,yang_contraction_2019,han_isotonic_2019,gao_estimation_2020,hanNoisyLinearInverse2023}; we refer the reader to P. Bellec's PhD thesis (\cite{bellecSharpOracleInequalities2016}) for a comprehensive understanding of shape-constrained regression. We also refer the reader to the special issue \cite{samworth_editorial_2018}, survey \cite{guntuboyina_nonparametric_2018}, and Samworth's ICM talk \cite{samworthNonparametricInferenceShape2026} for recent developments in nonparametric shape-constrained regression. The isotonic regression refers to the following estimator
\begin{align*}
    \hat\bbeta\in\argmin\Bigl(\frac{1}{2}\|\vy-\bbeta\|_2^2:\,\bbeta\in T\Bigr),\quad\mbox{ where }\quad T=\Bigl\{\bbeta = (\beta_j)_{j=1}^N \in\bR^N: \beta_1\leq\beta_2\leq\cdots\leq\beta_N\Bigr\}.
\end{align*}That is, we project $\vy$ onto the ``increasing'' index set in Euclidean distance to obtain $\hat\bbeta$; this kind of restriction on the ``shape'' of $\hat\bbeta$ is what makes this class of problems known as shape-constrained regression. We refer interested readers to \cite{chatterjeeNewPerspectiveLeast2014} and the references therein for a comprehensive understanding of this problem. Isotonic regression is well known for its estimation error having a rate of $N^{1/3}$ in the general case; in particular, Chatterjee, in \cite[Theorem 2.2]{chatterjeeNewPerspectiveLeast2014}, obtained the following result by using Theorem~\ref{theo:chatterjee}.
\begin{Theorem}[\cite{chatterjeeNewPerspectiveLeast2014}, Theorem~2.2]\label{theo:chatterjee_isotonic}
    Suppose $N\geq 2$. Define $A=N\min(\beta_{j+1}^\star-\beta_j^\star: 1\leq j<N)$, $B=N\max(\beta_{j+1}^\star-\beta_j^\star:1\leq j<N)$, and $D=\max\{\beta_N^\star-\beta_1^\star,\sigma_\xi\}$. Then there exist absolute constants $c,C>0$ such that $t_{\bbeta^\star}$ defined in Theorem~\ref{theo:chatterjee} satisfies
    \begin{align*}
        c\sigma_\xi^{2/3}\frac{A^{8/3}}{B^{4/3}D}N^{\frac{1}{6}}\leq t_{\bbeta^\star}\leq C\sigma_\xi^{2/3}D^{1/3}N^{\frac{1}{6}}.
    \end{align*}Moreover, when $t_{\bbeta^\star}\gtrsim\sigma_\xi$, $\bE\|\hat\bbeta-\bbeta^\star\|_2^2\asymp t_{\bbeta^\star}^2$.
\end{Theorem}
In particular, when $\beta_{j+1}^\star-\beta_j^\star\asymp \frac{\sigma_\xi}{N}$, Theorem~\ref{theo:chatterjee_isotonic} gives $\mathbb{E}\|\hat\bbeta-\bbeta^\star\|_2^2\asymp \sigma_\xi^2 N^{1/3}$. However, precisely because of the difficulty of estimating the offset localized Gaussian mean width when applying Theorem~\ref{theo:chatterjee}, the estimate of $t_{\bbeta^\star}$ provided by Theorem~\ref{theo:chatterjee_isotonic} is not always meaningful; we give an example below.
\begin{Example}\label{example:isotonic}
    Let $N\geq4$ be an even integer. Let $\sigma_\xi>0$, and take
    \begin{align*}
        \beta_i^\star = \begin{cases}
            0, & 1\leq i\leq \frac{1}{2}N,\\
            \sigma_\xi, & \frac{1}{2}N<i\leq N.
        \end{cases}
    \end{align*}
    Then $A=0$, $B=N\sigma_\xi$, and $D=\sigma_\xi$, and $0\leq t_{\bbeta^\star}\lesssim \sigma_\xi N^{1/6}$. Theorem~\ref{theo:chatterjee_isotonic} implies that if $t_{\bbeta^\star}\gtrsim\sigma_\xi$, then $0\leq\bE\|\hat\bbeta-\bbeta^\star\|_2^2\leq C \sigma_\xi^2 N^{\frac{1}{3}}$. Therefore, for this particular $\bbeta^\star$, Theorem~\ref{theo:chatterjee_isotonic} cannot provide a meaningful conclusion — it is neither clear how to verify $t_{\bbeta^\star}\gtrsim\sigma_\xi$, nor can we obtain matching upper and lower bounds. In fact, we show below using our sharp rate that for this example, $\mathbb{E}\|\hat\bbeta-\bbeta^\star\|_2^2\asymp \sigma_\xi^2\log(N)$, so that the condition $t_{\bbeta^\star}\gtrsim\sigma_\xi$ in Theorem~\ref{theo:chatterjee_isotonic} could have been satisfied, yet its upper estimate of $t_{\bbeta^\star}$ cannot provide the sharp rate of the expected estimation error. In fact, if there exists \(j\in[N-1]\) such that \(\beta_{j+1}^\star=\beta_j^\star\), then \(A=0\), and Theorem~\ref{theo:chatterjee_isotonic} fails giving meaningful result in this case. This is because Theorem~\ref{theo:chatterjee_isotonic} does not make use of the local structure of $\bbeta^\star$.
\end{Example}
Below we present the main result of this subsection: an explicit, user-friendly instance-specific sharp rate for the expected estimation error of isotonic regression that holds for any $\sigma_\xi>0$, any $\bbeta^\star\in T$, and any $N\in\mathbb{N}_+$.

For any $\bbeta^\star\in T$ and any $1\leq i\leq N$, we define
\begin{align}\label{eq:def_k_star_i_pm}
    \begin{aligned}
        k_{\star,i}^+&=\max\biggl\{ 1\leq k\leq N-i+1:\, \beta_{i+k-1}^\star-\beta_i^\star\leq\frac{\sigma_\xi}{\sqrt{k}} \biggr\},\mbox{ and }
        k_{\star,i}^-=\max\biggl\{ 1\leq k\leq i:\, \beta_i^\star - \beta_{i-k+1}^\star \leq\frac{\sigma_\xi}{\sqrt{k}} \biggr\}.
    \end{aligned}
\end{align}The proof of the following theorem can be found in Section~\ref{sec:proof_theo_isotonic}.
\begin{Theorem}\label{theo:isotonic}
     There exists an absolute constant $c>0$ such that, for any $N\in\bN_+$, any $\sigma_\xi>0$, and any $\bbeta^\star\in T$, 
     \begin{align*}
         c r_\star^2 \leq \bE\|\hat\bbeta-\bbeta^\star\|_2^2\leq 10r_\star^2,\quad\mbox{ where }\quad r_\star^2 = \sigma_\xi^2\sum_{i=1}^N\biggl(\frac{1}{k_{\star,i}^+}+\frac{1}{k_{\star,i}^-}\biggr).
     \end{align*}
\end{Theorem}
Theorem~\ref{theo:isotonic} establishes explicit upper and lower bounds that match for the expected estimation error of isotonic regression in the Gaussian sequence model, for any $\sigma_\xi>0$ and any $\bbeta^\star$. To the best of our knowledge, there has been no result in the existing literature that, like Theorem~\ref{theo:isotonic}, provides matching upper and lower bounds on the expected estimation error for any \(\bbeta^\star\in T\), any \(\sigma_\xi>0\), and any \(N\in\bN_{+}\), and hence we say that Theorem~\ref{theo:isotonic} fully characterizes the expected estimation error of well-specified univariate isotonic regression. Existing results in the literature either give only upper bounds (\cite{zhang_risk_2002,chatterjee_risk_2015,bellecSharpOracleInequalities2018,gao_estimation_2020}), or give matching upper and lower bounds only under fairly special assumptions on the signal and the noise (\cite{chatterjeeNewPerspectiveLeast2014,soloff_distribution-free_2019}). We use the following example to illustrate that this instance-specific sharp rate is also user-friendly.
\begin{Examplebis}{example:isotonic}
    We compute $k_{\star,i}^+$ and $k_{\star,i}^-$. First, fix any $1\leq i\leq\frac{1}{2}N$. For any $1\leq k\leq \frac{1}{2}N-i+1$, $\beta_{i+k-1}^\star-\beta_i^\star=0$; while for any $\frac{1}{2}N-i+2\leq k\leq N-i+1$, we have $k\geq2$ and $\sqrt{k}(\beta_{i+k-1}^\star-\beta_i^\star)=\sqrt{k}\sigma_\xi>\sigma_\xi$. Hence $k_{\star,i}^+=\frac{1}{2}N-i+1$, for all $1\leq i\leq \frac{1}{2}N$. Now fix any $\frac{1}{2}N<i\leq N$, then for any $1\leq k\leq N-i+1$, we have $\frac{1}{2}N<i\leq i+k-1\leq N$, hence $\sqrt{k}(\beta_{i+k-1}^\star-\beta_i^\star)=0$, hence by the definition, $k_{\star,i}^+=N-i+1$. Computing $k_{\star,i}^-$ similarly, we obtain
    \begin{align*}
        \forall 1\leq i\leq N,\qquad k_{\star,i}^+=\begin{cases}
            \frac{1}{2}N-i+1, & 1\leq i\leq \frac{1}{2}N,\\
            N-i+1, & \frac{1}{2}N<i\leq N,
        \end{cases}\quad \mbox{ and }\quad k_{\star,i}^-=\begin{cases}
            i, & 1\leq i\leq \frac{1}{2}N,\\
            i-\frac{1}{2}N, & \frac{1}{2}N<i\leq N.
        \end{cases}
    \end{align*}Therefore, by reindexing the sum,
    \begin{align*}
        \sigma_\xi^{-2}r_\star^2&=\sum_{i=1}^{N/2}\Bigl(\frac{1}{i}+\frac{1}{\frac{1}{2}N-i+1}\Bigr) + \sum_{i>N/2}^N\Bigl(\frac{1}{i-\frac{1}{2}N} + \frac{1}{N-i+1}\Bigr) = 4\sum_{i=1}^{N/2}\frac{1}{i} \asymp \log(N).
    \end{align*}Theorem~\ref{theo:isotonic} therefore implies that $\bE\|\hat\bbeta-\bbeta^\star\|_2^2\asymp \sigma_\xi^2\log(N)$.
\end{Examplebis}
In the following discussion, we take ``isotonic regression'' as the subject in order to reflect that the properties discussed below are not artificial byproducts of the analysis, but rather the intrinsic structure of isotonic regression. Theorem~\ref{theo:isotonic} shows that the sharp rate of isotonic regression depends on the sequence $(k_{\star,i}^\pm)_{i=1}^N$, and this sequence has a natural interpretation, as we show in the following. 

For any \(i\in[N]\), \(k_{\star,i}^\mp\) defines two intervals that intersect at \(i\), namely \(I_i^-:=\{i-k_{\star,i}^-+1,\cdots,i\}\) and \(I_i^+:=\{i,\cdots,i+k_{\star,i}^+-1\}\). The lengths of these two intervals are \(k_{\star,i}^\mp\), respectively. Since \(k_{\star,i}^\mp\) depends on \(\bbeta^\star,\sigma_\xi\), we say that these two intervals are \emph{adaptive}; and since they lie on either side of $i$, we say they are \emph{local}. Below we explain that \(\frac{\sigma_\xi}{\sqrt{k_{\star,i}^\mp}}\) can be regarded as a uniform upper bound on the variation over the two intervals \(I_i^\mp\) at position \(i\), or equivalently, at each position \(i\), isotonic regression filters the indices according to the variation and selects all maximal index sets whose variation does not exceed this threshold. For any \(j,\ell\in I_i^\mp\) with \(j\leq \ell\), we denote \(\Delta_{\ell,j}=\beta_\ell^\star-\beta_j^\star\). Then \(\max(\Delta_{\ell,j}: j,\ell\in I_i^\mp, j\leq \ell) \leq \frac{\sigma_\xi}{\sqrt{k_{\star,i}^\mp}}\), that is, the variation over \(I_i^\mp\) has this uniform upper bound, and hence we may regard it as a kind of variation threshold. The more slowly the signal varies, the smaller this threshold is, that is, the larger \(k_{\star,i}^\mp\) is. Since \(r_\star^2\) is the sum of squares of these thresholds, isotonic regression prefers signals that vary more smoothly, that is, for signals that vary more smoothly, its expected estimation error is smaller. Therefore, for each \(i\in[N]\), isotonic regression identifies the largest interval on which the uniform variation of the signal does not exceed a given threshold. Isotonic regression repeats the above operation for each position \(i\in[N]\) of the signal, and takes the sum of squares of the thresholds in the final estimation error.

Below we provide more examples to illustrate that \(r_\star^2\) is user-friendly.

\begin{Example}\label{example:uniform_increase}
    Let $N\geq2$, $\sigma_\xi>0$, and $v\geq0$. Let $\beta_1^\star\in\bR$, and $\beta_i^\star=\beta_1^\star+v(i-1)$. Then one may compute that $k_{\star,1}^+\asymp \min\bigl\{N,\max\{1,(\frac{\sigma_\xi}{v})^{2/3}\}\bigl\}$, $k_{\star,i}^+=\min\{N-i+1,k_{\star,1}^+\}$, and $k_{\star,i}^-=\min\{i,k_{\star,1}^+\}$. Consequently,
    \begin{align*}
        \bE\|\hat\bbeta-\bbeta^\star\|_2^2\asymp \begin{cases}
            \sigma_\xi^2 \log(N), & 0\leq v\leq \sigma_\xi\Bigl(\frac{\log(N)}{N}\Bigr)^{3/2},\\
            N\sigma_\xi^{4/3}v^{2/3}, & \sigma_\xi\Bigl(\frac{\log(N)}{N}\Bigr)^{3/2} < v \leq \sigma_\xi,\\
            N\sigma_\xi^2, & v>\sigma_\xi.
        \end{cases}
    \end{align*}That is, for a signal \(\bbeta^\star\) growing uniformly at speed \(v\), there exist two thresholds in terms of the speed \(v\), namely \(\sigma_\xi(\frac{\log(N)}{N})^{3/2}\) and \(\sigma_\xi\), such that when \(v\) does not exceed \(\sigma_\xi(\frac{\log(N)}{N})^{3/2}\), the estimation error does not increase as \(v\) increases---that is, for slowly growing signals, the learning difficulty is no different from that for constant signals; when \(v\) lies between the two thresholds, the estimation error increases rapidly at the rate \(v^{2/3}\) as \(v\) increases---the common case \(v=\frac{\sigma_\xi}{N}\) corresponds to \(\sigma_\xi^2 N^{1/3}\), while \(v=\frac{\sigma_\xi}{\sqrt{N}}\) corresponds to \(\sigma_\xi^2 N^{2/3}\); when \(v\) exceeds \(\sigma_\xi\), the growth of the estimation error with respect to the speed \(v\) enters a saturation stage, at which point increasing \(v\) no longer makes the estimation error increase faster --- in this case, \(\frac1N\bE\|\hat\bbeta-\bbeta^\star\|_2^2\) no longer tends to \(0\) as \(N\to\infty\) (given a constant \(\sigma_\xi\)). Qualitatively speaking, according to Theorem~\ref{theo:isotonic}, as \(v\) increases, the uniform variation increases, so that the adaptive threshold \(\frac{\sigma_\xi}{\sqrt{k_{\star,i}^\mp}}\) has to increase accordingly; quantitatively, the rate of increase of the variation is proportional to the difference of the indices, that is, \(\beta_{i+k-1}^\star-\beta_i^\star = v(k-1)\); hence, for it not to exceed the threshold \(\frac{\sigma_\xi}{\sqrt{k}}\), \(k_{\star,i}^\mp\) naturally needs to be proportional to \(\bigl(\frac{\sigma_\xi}{v}\bigr)^{2/3}\), which explains the origin of this exponent \(\frac23\).
\end{Example}

\begin{Example}\label{example:K_plateaus}
    Let $1\leq K\leq N$, and suppose that $K$ divides $N$. Let $\beta_i^\star = \sigma_\xi\lfloor \frac{K(i-1)}{N}\rfloor$. Then $\bbeta^\star$ has $K$ plateaus of equal width, and for any $1\leq \ell\leq K$, such that $\frac{(\ell-1)N}{K}<i\leq \frac{\ell N}{K}$, then $\beta_i^\star=(\ell-1)\sigma_\xi$. For any $1\leq k\leq i-\frac{(\ell-1)N}{K}$, we have $i-k+1\leq i\leq \frac{\ell N}{K}$; while $i-k+1\geq \frac{(\ell-1)N}{K}+1$, so $\beta_{i-k+1}^\star=\beta_i^\star$; whereas for $k>i-\frac{(\ell-1)N}{K}$, $i-k+1\leq\frac{(\ell-1)N}{K}$ and hence $\beta_{i-k+1}^\star=\ell\sigma_\xi$. Therefore $k_{\star,i}^-=i-(\ell-1)\frac{N}{K}$. Similarly, $k_{\star,i}^+=\ell\frac{N}{K}-i+1$. Consequently,
    \begin{align*}
        \bE\|\hat\bbeta-\bbeta^\star\|_2^2 \asymp \sigma_\xi^2\sum_{\ell=1}^K\sum_{i=\frac{(\ell-1)N}{K+1}}^{\frac{\ell N}{K}}\Bigl( \frac{1}{i-(\ell-1)\frac{N}{K}} + \frac{1}{\ell\frac{N}{K}-i+1} \Bigr) = 2\sigma_\xi^2\sum_{j=1}^{\frac{N}{K}}\frac{K}{j} \asymp \sigma_\xi^2 K \log\Bigl(\frac{eN}{K}\Bigr).
    \end{align*}\cite{chatterjee_risk_2015,gao_estimation_2020} proved that the worst-case estimation error lower bound for the set of non-decreasing signals consisting of at most \(K\) constant pieces is \(\sigma_\xi^2 K\log(\frac{eN}{K})\), and hence this example shows that all these equal-width signals are in the worst position.
\end{Example}

\begin{Example}\label{example:plateau_transition_plateau}
    Suppose that $3$ divides $N$, $\sigma_\xi>0$, and $1\leq m\leq \frac{1}{3}N$. Let $\beta_i^\star = \sigma_\xi\min\{1,\frac{1}{m}(i-\frac{N}{3})_+\}$. Then one may compute that
    \begin{align*}
        \bE\|\hat\bbeta-\bbeta^\star\|_2^2 \asymp \sigma_\xi^2\Bigl(\log(eN)+m^{1/3}\Bigr).
    \end{align*}Since the computation in this example is relatively involved, we place it in item~\emph{3.} of Section~\ref{sec:examples_revisited}.
    In dose-response studies, a typical monotone signal exhibits a plateau-transition-plateau structure: the treatment effect remains nearly constant at low doses, increases through a transition region, and eventually saturates at high doses. Such response patterns are commonly modeled by logistic curves; see, for instance, \cite{bornkamp_response-adaptive_2011}. This example is a simplified version of this model.

    This example reveals an interesting phenomenon: since the signal is monotonically increasing and the two plateaus are fixed at $0$ and $\sigma_\xi$, its total variation is always $\beta_N^\star-\beta_1^\star=\sigma_\xi$, independently of the transition width $m$. Nevertheless, as $m$ increases, the transition becomes more gradual, while the estimation error of the isotonic LSE increases. In particular, when $m\lesssim(\log(N))^3$, the estimation error remains of logarithmic order; when $m\asymp N$, it reaches the order $\sigma_\xi^2N^{1/3}$. This sharply contrasts with the classical risk upper bounds depending on the total variation \cite{zhang_risk_2002,chatterjee_risk_2015}, which yield the same upper bound of order $\sigma_\xi^2N^{1/3}$ for all $m$ in this example. In contrast, the instance-specific matching upper and lower bounds in Theorem~\ref{theo:isotonic} reveal that two monotone signals with exactly the same total variation can have expected estimation errors differing by a polynomial factor. This demonstrates that the local variation structure of the signal plays an essential role in precisely characterizing the estimation error of the isotonic LSE.
\end{Example}


\subsection{Explicit Instance-Specific Sharp Rates for Gaussian Linear Regression in the Asymptotic Regime}\label{sec:asymptotic_isotonic}

Recently, Han \cite{hanNoisyLinearInverse2023} established a reduction method from the constrained LSE of Gaussian linear regression to the constrained LSE of a sequence model, defining a fixed-point equation via the expected estimation error of the sequence model with an \emph{effective} noise level, and proved that the solution of this equation is, in the asymptotic regime, the leading term of the estimation error of the original Gaussian linear regression. This result establishes an instance-specific connection between the estimation errors of the two setups under fairly general assumptions; however, solving this fixed point is usually difficult. For example, \cite[Corollary 3.4]{hanNoisyLinearInverse2023} still yields only an upper bound on the estimation error, rather than matching upper and lower bounds. The instance-specific explicit sharp rates provided by Theorem~\ref{theo:isotonic} precisely complement this point. In this subsection, for the example of isotonic regression, we present a corollary combining Theorem~\ref{theo:isotonic} with the conclusions in \cite{hanNoisyLinearInverse2023}, which characterizes the explicit instance-specific sharp convergence rate of isotonic regression in the Gaussian linear regression setup in the limiting regime.

To this end, we first specify the setup---these setups are needed only in this subsection. We let \(n\in\bN_{+}\) be the sample size, while the original \(N\in\bN_{+}\) is the dimension of the space in which the design lies, that is, let \(X_1^{(N)},\cdots,X_n^{(N)}\) be \(n\) independent copies of a standard Gaussian random vector in \(\bR^N\). Let \(\bxi=(\xi_i^{(N)})_{i=1}^n\) be a centered random vector independent of \(X_1^{(N)},\cdots,X_n^{(N)}\) such that $\bE[(\xi_i^{(N)})^2]=\sigma_\xi^2$ and such that its coordinates are independent. We still let \(T = \{\bbeta\in\bR^N:\beta_1\leq\cdots\leq\beta_N\}\). Let \(\bbeta^{\star,(N)}\in T\) be the signal. Since the results in \cite{hanNoisyLinearInverse2023} hold only in the asymptotic regime, we define a sequence of problems indexed by the dimension \(N\geq 2\), namely \(\{(n(N),\bbeta^{\star,(N)}, (X_i^{(N)},\xi_i^{(N)})_{i=1}^{n(N)})\}_{N\geq 2}\). The limits considered below are all as \(N\to\infty\) while $\sigma_\xi^2$ being fixed. For convenience, we omit the superscript \((N)\) in the following writing. Let \(Y_i=\langle \bbeta^\star,X_i\rangle + \xi_i\) be the response variable, where \(\xi_i\) is the \(i\)-th component of \(\bxi\). Let $\hat\bbeta\in\argmin\Bigl(\frac{1}{2n}\sum_{i=1}^{n}\bigl(Y_i -\langle \bbeta,X_i\rangle\bigr)^2: \bbeta\in T\Bigr)$ be the isotonic regression. \cite{hanNoisyLinearInverse2023} establishes a connection between the estimation error of the convex-constrained LSE under such a Gaussian linear regression setup and the expected estimation error of the constrained LSE of a Gaussian sequence model on \(\bR^N\) with \emph{effective} noise level \(\frac{\sigma_\xi}{\sqrt n}\), which we call the auxiliary problem. By Theorem~\ref{theo:isotonic}, we know that the expected estimation error of such an auxiliary problem depends on \(k_{\star,i}^\pm\), and we recall here its definition under this effective noise level as
\begin{align*}
    k_{\star,i}^+ = \max\Bigl\{1\leq k\leq N-i+1:\, \beta_{i+1-1}^\star-\beta_i^\star \leq \frac{\sigma_\xi}{\sqrt{nk}} \Bigr\},\mbox{ and } k_{\star,i}^- = \max\Bigl\{1\leq k\leq i:\, \beta_i^\star - \beta_{i-k+1}^\star \leq \frac{\sigma_\xi}{\sqrt{nk}} \Bigr\}.
\end{align*}Let $r_\star^2 = \frac{\sigma_\xi^2}{n}\sum_{i=1}^N\Bigl(\frac{1}{k_{\star,i}^+} + \frac{1}{k_{\star,i}^-}\Bigr)$. Han's reduction relies on the statistical dimension of the tangent cone of \(T\) (\cite{amelunxen_living_2014}). To this end, we define \(K = \bigl|\bigl\{\beta_i^\star:1\leq i\leq N\bigr\}\bigr|\), that is, the number of distinct components of \(\bbeta^\star\), and number its distinct coordinates as \(0=i_0<i_1<\cdots<i_K=N\), so that \(\beta_{i_{\ell-1}+1}^\star = \cdots = \beta_{i_\ell}^\star\), and for any \(\ell<K\), \(\beta_{i_\ell}^\star<\beta_{i_\ell+1}^\star\). From this we can introduce the following assumption.
\begin{Assumption}\label{assumption:han}
    As $N\to\infty$, the following hold:
    \begin{enumerate}
        \item $r_\star^2\to 0$;
        \item $\log(eK)+\log\log(16N) = o(nr_\star^2)$;
        \item either $r_\star^2\inf(\|\bbeta^\star-a\one_N\|_2^2:a\in\bR)\to0$, or $\limsup\frac{1}{n}\sum_{\ell=1}^K\sum_{j=1}^{i_\ell-i_{\ell-1}}\frac{1}{j}<2$.
    \end{enumerate}
\end{Assumption}The following corollary is the main result of this subsection, and its proof can be found in Section~\ref{sec:proof_coro_han}.
\begin{Corollary}\label{coro:asymptotic_Gaussian}
    Grant Assumption~\ref{assumption:han}, then there exists an absolute constant $c>0$ such that for any fixed $\varepsilon>0$, there exists $N_0\in\bN_+$ such that
    \begin{align*}
        \forall N\geq N_0,\qquad \bP\Bigl(cr_\star^2\leq \bigl\|\hat\bbeta-\bbeta^\star\bigr\|_2^2 \leq 20r_\star^2\Bigr) \geq 1-\varepsilon.
    \end{align*}
\end{Corollary}
To the best of our knowledge, Corollary~\ref{coro:asymptotic_Gaussian} is the first result that gives explicit, user-friendly matching upper and lower bounds on the estimation error under such broad conditions. The works in the literature either give such conclusions only for very special signals (\cite{durot_Lp_2007}), or are restricted to the small-noise limit (\cite{soloff_distribution-free_2019}); or hold for general signals and noise but are neither explicit nor user-friendly (\cite{hanNoisyLinearInverse2023,chatterjeeNewPerspectiveLeast2014}). Below we present some examples.
\begin{Examplebis}{example:uniform_increase}
    One may compute that $\inf(\|\bbeta^\star-a\one_N\|_2^2: a\in\bR) = \frac{1}{12}v^2N(N^2-1)$. Suppose $\log(N)+N\min\{1,(v\sqrt{n}/\sigma_\xi)^{2/3}\}=o(n)$; $\frac{v\sqrt{n}}{\sigma_\xi}(\frac{N}{\log(N)})^{3/2}\to\infty$; and either $r_\star^2v^2N(N^2-1)\to0$ or $\limsup\frac{N}{n}<2$, then Assumption~\ref{assumption:han} holds. In particular, when $v=\frac{\sigma_\xi}{N}$, then $\limsup \frac{N}{n}<2$ is sufficient for all the assumptions and $\bigl\|\hat\bbeta-\bbeta^\star\bigr\|_2^2\asymp \sigma_\xi^2\min\{\frac{N}{n},\, \frac{N^{1/3}}{n^{2/3}}\}$ tending to probability $1$; when $v=\frac{\sigma_\xi}{\sqrt{N}}$, then $N=o(n)$ suffices and $\bigl\|\hat\bbeta-\bbeta^\star\bigr\|_2^2\asymp \sigma_\xi^2\frac{N}{n}$, tending to probability $1$.
\end{Examplebis}
\begin{Examplebis}{example:K_plateaus}
    Suppose further that $K\log(\frac{eN}{K})=o(n)$, and $n=\omega(\log(N))$. Then $\bigl\|\hat\bbeta-\bbeta^\star\bigr\|_2^2\asymp \sigma_\xi^2\frac{K\log(\frac{eN}{K})}{n}$ happens with probability tending to $1$.
\end{Examplebis}
\begin{Examplebis}{example:plateau_transition_plateau}
    Suppose further $\log(N)=o(n)$, and $\limsup\frac{m}{n}<2$, then Assumption~\ref{assumption:han} holds, and $\|\hat\bbeta-\bbeta^\star\|_2^2 \asymp \frac{\sigma_\xi^2}{n}\bigl(\log(N)+\min\{m,(nm)^{1/3}\}\bigr)$ tending to probability $1$.
\end{Examplebis}
At the end of this section, we point out that other conclusions of this paper, such as those for the \(\ell_q\) ellipsoid, the \(\ell_1\) ball, and the \(k\)-interpolated body, can also be combined with the results of \cite{hanNoisyLinearInverse2023} to obtain instance-specific sharp rates for Gaussian linear regression in the asymptotic regime. Due to length constraints, we do not demonstrate this here.

\section{Exemplar of our analysis: \texorpdfstring{$\ell_q$}{lq} ellipsoid}\label{sec:exemplar_ell_q}
We prove only the Gaussian version of the conclusion. For the conclusion under noise satisfying a Poincar\'e inequality, we briefly explain at the end (see Remark~\ref{rmk:Poincare_noise_proof}) how to modify the proof.
\subsection{Preliminaries}\label{sec:preliminary}

We first establish the existence and uniqueness of $\lambda_\self$ defined by \eqref{eq:general_self_calibration}.
\begin{Lemmabis}{lem:calibration_equation_l_q}
    Under the assumptions of Theorem~\ref{theo:ell_q}, $\lambda_\self>0$ exists uniquely, and satisfies
    \begin{align}\tag{\ref{eq:general_self_calibration}}
    \sum_{j=1}^N\sigma_j^{-q/2}\bE|\eta_{\self,j}(y_j)|^q=\bE\min\Bigl\{1,\sum_{j=1}^N\sigma_j^{-q/2}|y_j|^q\Bigr\},\mbox{ where }\eta_{\self,j}(\cdot) = \argmin\Bigl(\frac{1}{2}(\cdot-s)^2+\frac{\lambda_\self}{q}\sigma_j^{-q/2}|s|^q: s\in\bR\Bigr).
    \end{align}
\end{Lemmabis}
In this proof, we write $T$ in the following form
\begin{align*}
 T=\{\bbeta\in\mathbb R^N:\Psi(\bbeta)\leq0\},\quad \mbox{ where }\quad 
 \Psi(\bbeta)=\frac{1}{q}(\|\Sigma^{-1/2}\bbeta\|_q^q-1).
\end{align*}
Write $\Prox_{\lambda }$ for the proximal operator of $\lambda \Psi$, and abbreviate  $d_{\lambda_\self}$ as $d_\self$, where $d_\lambda$ is defined in \eqref{eq:def_divergence}. For any $\vx\in\bR^N$, let
\begin{align}\label{eq:def_implicit_regularization}
 \bar\lambda[\vx]:=\min\left\{\lambda\geq0:
 \Mo(\lambda,\vx)=\frac12\|\vx-\Proj_T[\vx]\|_2^2\right\}.
\end{align}The $\bar\lambda$ defined above always exists uniquely and is a continuous mapping.
\begin{Lemmabis}{lem:existence_l_q}
    We have $\bar\lambda[\bullet]=\|\Sigma^{1/2}(\bullet-\Proj_T[\bullet])\|_{\frac{q}{q-1}}$, hence is unique and is continuous.
\end{Lemmabis}

For the constrained LSE studied in this work, we introduce the following two novel ingredients: 
\begin{enumerate}
    \item using implicit regularization to transform the constrained LSE into an unconstrained RERM with a data-dependent regularization functional; and 
    \item constructing a surrogate estimator $\tilde\bbeta$ that has a deterministic regularization functional and is an unconstrained RERM, so that it can be conveniently handled, and whose approximation error to $\hat\bbeta$ can be shown to satisfy $\mathbb{E}\|\hat\bbeta-\tilde\bbeta\|_2^2\leq  \sigma_\xi^2$. 
\end{enumerate}

\subsection{Implicit regularization: from constrained to unconstrained problems}\label{sec:implicit_regularization}

Below, we abbreviate $\hat\bbeta[\vy]$ and $\bar\lambda[\vy]$ as $\hat\bbeta$ and $\hat\lambda$, respectively.

In this subsection, we show how to transform the constrained LSE into an unconstrained RERM with a data-dependent regularization functional by exploring the implicit regularization of the constrained LSE. The basic analytical tools are the KKT conditions from convex analysis. In this subsection, we prove the following proposition. 
\begin{Proposition}\label{prop:implicit_regularization}
    Under the conditions of Theorem~\ref{theo:ell_q}, the following hold:
    \begin{enumerate}
        \item For any $\vx\in\bR^N$ and any $\lambda\geq0$, 
        \begin{align*}
     d_\lambda(\vx)=\Mo(\bar\lambda[\vx],\vx)-\Mo(\lambda,\vx)
     +(\lambda-\bar\lambda[\vx])\Psi(\hat\bbeta[\vx])\in\bR_+,
    \end{align*}
    
    \item For any $\vy\in\bR^N$, recall that $\hat\lambda=\bar\lambda[\vy]$. Then
    \begin{align}\label{eq:implicit_regularization}
        \hat\bbeta = \Proj_T[\vy] = \argmin\Bigl( \frac{1}{2}\|\vy-\bbeta\|_2^2 + \hat\lambda \Psi(\bbeta): \bbeta\in\bR^N \Bigr),\quad\mbox{ and }\quad \hat\lambda \Psi(\hat\bbeta)=0.
    \end{align}
    
    \item Moreover, for any $\lambda>0$, $\Prox_\lambda[\bullet] = \sum_{j=1}^N \eta_{\lambda,j}(\bullet_j)\ve_j$, and each $\eta_j$ is non-decreasing, $1$-Lipschitz, and satisfies $\eta_j(0)=0$; where $$\eta_{\lambda,j}(\cdot)=\langle\ve_j,\Prox_\lambda[\cdot\ve_j]\rangle = \argmin\Bigl(\frac{1}{2}(\cdot-s)^2+\frac{\lambda}{q}\sigma_j^{-q/2}|s|^q:s\in\bR\Bigr).$$
    \end{enumerate}
\end{Proposition}
    \begin{proof}[Proof of Proposition~\ref{prop:implicit_regularization}]
    Suppose item~\emph{2.} is verified, we first prove item~\emph{1.}. The identity for $d_\lambda(\vx)$ in item~\emph{1.} is a direct consequence of the definition \eqref{eq:def_implicit_regularization} of $\bar\lambda[\vx]$, $\hat\lambda=\bar\lambda[\vy]$, $\hat\lambda\Psi(\hat\bbeta)=0$ from \eqref{eq:implicit_regularization}, and the following identity
    \begin{align}\tag{\ref{eq:identity_Moreau}}
        \Mo\bigl(\bar\lambda[\bullet],\bullet\bigr) = \frac{1}{2}\bigl\|\bullet-\Proj_T[\bullet]\bigr\|_2^2 = \frac{1}{2}\dist^2(\bullet,T).
    \end{align}Moreover, by Danskin's theorem, see, for instance, Lemma~\ref{lem:danskin}, $d_\lambda(\vx)$ is the Bregman divergence of the convex function $\lambda\mapsto -\Mo(\lambda,\vx)$ evaluated at $\bar\lambda[\vx]$, hence $d_\lambda(\vx)\geq 0$ for any $\vx\in\bR^N$.
    
    We now prove item~\emph{2.} By \eqref{eq:identity_Moreau}, item~\emph{1.} and item~\emph{2.} of Lemma~\ref{lem:moreau_proximal_equi}, $\Prox_{\hat\lambda}[\vy]=\Proj_T[\vy]=\hat\bbeta$, and $\hat\lambda \Psi(\hat\bbeta)=0$. By the definition of $\Prox_{\hat\lambda}[\vy]$ and $\hat\lambda=\bar\lambda[\vy]$, \eqref{eq:implicit_regularization} holds.

    We next prove the ``moreover'' part. The definition of $\Psi$ gives
    \begin{align*}
      \forall\bbeta\in\bR^N,\qquad \Psi(\bbeta)\geq \Psi(\vzero),\quad\mbox{ and }\quad \lambda \Psi(\bbeta)
      =\sum_{j=1}^N\lambda \Psi(\beta_j\ve_j)-(N-1)\lambda \Psi(\vzero).
    \end{align*}
    By the definition of $\Prox_\lambda$,
    \begin{align*}
        \Prox_\lambda[\bullet] = \argmin\Bigl(\frac{1}{2}\|\bullet-\bbeta\|_2^2 + \lambda \Psi(\bbeta):\bbeta\in\bR^N\Bigr) = \argmin\Bigl(\sum_{j=1}^N\bigl( \frac{1}{2}(\bullet_j-\beta_j)^2 + \lambda \Psi(\beta_j\ve_j) \bigr):\bbeta\in\bR^N\Bigr),
    \end{align*}where we have used the fact that $\lambda \Psi(\vzero)$ does not depend on $\bbeta$. As a result,
    \begin{align*}
        \langle\ve_j,\Prox_\lambda[\bullet]\rangle = \argmin\Bigl(\frac{1}{2}(\bullet_j-t)^2 + \lambda \Psi(t\ve_j) :t\in\bR\Bigr) = \langle\ve_j,\Prox_\lambda[\bullet_j\ve_j]\rangle = \eta_j(\bullet_j).
    \end{align*}From the definition of $\Psi$, we know $\vzero\in\argmin(\lambda \Psi(\bbeta):\bbeta\in\bR^N)$, hence \cite[Proposition 12.28]{bauschkeConvexAnalysisMonotone2011} implies that $\Prox_\lambda(\vzero)=\vzero$, and consequently $\eta_j(0)=0$; while the other properties of $\eta_j$ follow from the fact that it is a proximal operator, see, for instance, Section~\ref{sec:convex_analysis}.
\end{proof}
Proposition~\ref{prop:implicit_regularization} states that the constrained LSE $\hat\bbeta$ can be written as an unconstrained RERM with a data-dependent regularization term $\hat\lambda \Psi$.

\subsection{Surrogate estimator}

In this subsection, we construct an RERM $\tilde\bbeta$ with a deterministic regularization functional such that $\mathbb{E}\|\hat\bbeta-\tilde\bbeta\|_2^2 \leq \sigma_\xi^2$. In this subsection, we prove the following proposition. Let
\begin{align}\label{eq:def_surrogate_tilde_beta}
    \tilde\bbeta = \argmin\Bigl(\frac{1}{2}\|\vy-\bbeta\|_2^2+\lambda_\self \Psi(\bbeta):\bbeta\in\bR^N\Bigr) = \Prox_\self[\vy].
\end{align}
\begin{Lemma}\label{lem:coordinate_surrogate}
    For $\Psi(\bbeta)=\frac{1}{q}(\|\Sigma^{-1/2}\bbeta\|_q^q-1)$, $\forall 1\leq j\leq N$, we have $\tilde\beta_j = \eta_j(\beta_j^\star+\xi_j).$
\end{Lemma}
\begin{proof}[Proof of Lemma~\ref{lem:coordinate_surrogate}]
    From Proposition~\ref{prop:implicit_regularization}, $\Prox_\self[\vy]=\sum_{j=1}^N\eta_j(y_j)\ve_j$. Since $\tilde\bbeta = \Prox_\self[\vy]$, the claimed identity follows.
\end{proof}The aim of this subsection is to prove the following proposition.
\begin{Proposition}\label{prop:surrogate}
    Under the conditions of Theorem~\ref{theo:ell_q}, $\bE\|\hat\bbeta-\tilde\bbeta\|_2^2\leq \sigma_\xi^2$.
\end{Proposition}
We first prove the following two lemmas.
\begin{Lemma}\label{lem:d_varepsilon}
    For any $\lambda\in\bR_+$, $\frac{1}{2}\|\hat\bbeta-\Prox_\lambda[\vy]\|_2^2\leq d_\lambda(\vy)$.
\end{Lemma}
\begin{proof}[Proof of Lemma~\ref{lem:d_varepsilon}]
    By \eqref{eq:identity_Moreau} and $\hat\lambda=\bar\lambda[\vy]$, $\Mo(\hat\lambda,\vy)=\frac{1}{2}\|\vy-\hat\bbeta\|_2^2$. Recall that by definition, for $\bbeta\in\mathbb R^N$, $\lambda \Psi(\bbeta)$ is the regularization functional.
    By \eqref{eq:def_divergence} and Proposition~\ref{prop:implicit_regularization}, $\hat\lambda \Psi(\hat\bbeta)=0$, and we have
    \begin{align}\label{eq:d_varepsilon_y_proof}
        d_\lambda(\vy) &= \frac{1}{2}\|\vy - \hat\bbeta\|_2^2 - \Mo(\lambda,\vy) + \lambda \Psi(\hat\bbeta).
    \end{align}By \eqref{eq:identity_Moreau_Proximal} applied to $\Psi=\lambda \Psi(\cdot)$, together with the fact that for any $\lambda\in\bR_+$, $\Mo(\lambda,\cdot)=\min(\frac{1}{2}\|\cdot-\bbeta\|_2^2+\lambda \Psi(\bbeta):\bbeta\in\bR^N)$ is the Moreau envelope of $\lambda \Psi(\cdot)$; we obtain
    \begin{align*}
        \Mo(\lambda,\vy) = \frac{1}{2}\bigl\| \vy-\Prox_\lambda[\vy] \bigr\|_2^2 + \lambda \Psi\bigl(\Prox_\lambda[\vy]\bigr).
    \end{align*}Therefore,
    \begin{align*}
        d_\lambda(\vy) = \frac{1}{2}\|\vy-\hat\bbeta\|_2^2 - \frac{1}{2}\|\vy-\Prox_\lambda[\vy]\|_2^2 + \lambda \Psi(\hat\bbeta) - \lambda \Psi(\Prox_\lambda[\vy]).
    \end{align*}From \eqref{eq:gradient_Moreau}, $\vy-\Prox_\lambda[\vy]\in\partial^-(\lambda \Psi)(\Prox_\lambda[\vy])$, and consequently, $\lambda \Psi(\hat\bbeta)-\lambda \Psi(\Prox_\lambda[\vy])\geq\langle\vy-\Prox_\lambda[\vy],\hat\bbeta-\Prox_\lambda[\vy]\rangle$. In summary,
    \begin{align*}
        d_\lambda(\vy)&\geq \frac{1}{2}\|\vy-\hat\bbeta\|_2^2 - \frac{1}{2}\|\vy-\Prox_\lambda[\vy]\|_2^2 + \langle\vy-\Prox_\lambda[\vy],\hat\bbeta-\Prox_\lambda[\vy]\rangle = \frac{1}{2}\bigl\|\hat\bbeta - \Prox_\lambda[\vy]\bigr\|_2^2.
    \end{align*}
\end{proof}
\begin{Lemma}\label{lem:distance_waist}
    For any $\lambda\in\bR_+$, if $\bar\lambda^{-1}[\lambda]:=\{\vw\in\bR^N:\, \bar\lambda[\vw]=\lambda\}$ is non-empty, then $d_\lambda(\vy)\leq\frac12\dist^2(\vy,\bar\lambda^{-1}[\lambda])$.
\end{Lemma}
\begin{proof}[Proof of Lemma~\ref{lem:distance_waist}]
    By $\lambda \Psi(\hat\bbeta)\leq0$, \eqref{eq:d_varepsilon_y_proof} together with the fact that $\frac{1}{2}\|\vy-\hat\bbeta\|_2^2=\frac{1}{2}\dist^2(\vy,T)$, we have $d_\lambda(\vy)\leq\frac{1}{2}\dist^2(\vy,T)-\Mo(\lambda,\vy)$. We only need to prove that for any $\vw\in\bR^N$ such that $\bar\lambda[\vw]=\lambda$, we have
    \begin{align}\label{eq:objective_distance}
        \frac{1}{2}\dist^2(\vy,T)-\Mo(\lambda,\vy)\leq\frac{1}{2}\|\vy-\vw\|_2^2.
    \end{align}Notice that for any $\lambda\in\bR_+$, $\vx\in\bR^N\mapsto \Mo(\lambda,\vx)$ is convex, hence
    \begin{align*}
        \Mo(\lambda,\vy)-\Mo(\lambda,\vw)\geq \langle\vy-\vw,\nabla_{\bullet_2}\Mo(\lambda,\vw)\rangle.
    \end{align*}
    By \eqref{eq:derivative_Moreau}, $\frac{1}{2}\dist^2(\cdot,T)$ is Fr\'echet differentiable and its Fr\'echet derivative is $1$-Lipschitz. By the descent lemma, see, for instance, \cite[Theorem 18.15, (iii)]{bauschkeConvexAnalysisMonotone2011} applied to $\frac{1}{2}\dist^2(\cdot,T)$, together with  $\nabla_\vw\frac{1}{2}\dist^2(\vw,T) = \nabla_{\bullet_2}\Mo(\lambda,\vw)$ from \eqref{eq:gradient_Moreau_1} and \eqref{eq:derivative_Moreau}; we obtain
    \begin{align*}
        \frac{1}{2}\dist^2(\vy,T)&\leq\frac{1}{2}\dist^2(\vw,T) + \langle\vy-\vw,\nabla_\vw\frac{1}{2}\dist^2(\vw,T)\rangle+\frac{1}{2}\|\vy-\vw\|_2^2\\
        &\leq\frac{1}{2}\dist^2(\vw,T) + \Mo(\lambda,\vy)-\Mo(\lambda,\vw) + \frac{1}{2}\|\vy-\vw\|_2^2.
    \end{align*}Recall that $\vw\in\bR^N$ satisfies $\bar\lambda[\vw]=\lambda$, hence by \eqref{eq:identity_Moreau}, $\frac{1}{2}\dist^2(\vw,T) = \Mo(\bar\lambda[\vw],\vw)=\Mo(\lambda,\vw)$, which implies \eqref{eq:objective_distance}.
\end{proof}
Now we prove Proposition~\ref{prop:surrogate}.
\begin{proof}[Proof of Proposition~\ref{prop:surrogate}]
We divide this proof into the stochastic argument and the deterministic argument.

\medskip
\noindent\emph{Stochastic argument.}
By Lemma~\ref{lem:existence_l_q}, the map $\vu\in\bR^N\mapsto\bar\lambda[\bbeta^\star+\sigma_\xi\vu]\in\bR_+$ is continuous. Apply Lemma~\ref{lem:iso-gaussian-waist} to this map with $m=1$. Since $\bxi/\sigma_\xi\sim\cN(\vzero,I_N)$, there exists a deterministic $\lambda_0\geq0$ such that, for every $u>0$, the event $\Omega_u=\{\dist(\vy,\{\vw\in\bR^N:\bar\lambda[\vw]=\lambda_0\})\leq u\}$ satisfies $\bP(\Omega_u)\geq\bP(|\xi_1|\leq u)$. Indeed, the change of variables $\vw=\bbeta^\star+\sigma_\xi\vu$ gives
\begin{align*}
\dist\bigl(\bxi/\sigma_\xi,\{\vu\in\bR^N:\bar\lambda[\bbeta^\star+\sigma_\xi\vu]=\lambda_0\}\bigr)=\frac{1}{\sigma_\xi}\dist\bigl(\vy,\{\vw\in\bR^N:\bar\lambda[\vw]=\lambda_0\}\bigr).
\end{align*}
The level set is nonempty because $\lambda_0$ belongs to the image of the map. Tonelli's theorem consequently gives
\begin{align*}
\bE\dist^2\bigl(\vy,\{\vw\in\bR^N:\bar\lambda[\vw]=\lambda_0\}\bigr)=\int_0^\infty2u\bP(\Omega_u^c)\mathrm d u\leq\int_0^\infty2u\bP(|\xi_1|>u)\mathrm d u=\sigma_\xi^2.
\end{align*}

\medskip
\noindent\emph{Deterministic argument.}
By \eqref{eq:def_surrogate_tilde_beta}, $\tilde\bbeta=\Prox_{\lambda_\self}[\vy]$. Lemma~\ref{lem:d_varepsilon}, the minimizing property of $\lambda_\self$ in \eqref{eq:def_divergence}, and Lemma~\ref{lem:distance_waist}, applied in this order, yield
\begin{align*}
\frac{1}{2}\bE\|\hat\bbeta-\tilde\bbeta\|_2^2\leq\bE d_{\lambda_\self}(\vy)\leq\bE d_{\lambda_0}(\vy)\leq\frac{1}{2}\bE\dist^2\bigl(\vy,\{\vw\in\bR^N:\bar\lambda[\vw]=\lambda_0\}\bigr)\leq\frac{1}{2}\sigma_\xi^2.
\end{align*}
\end{proof}

\subsection{End of the proof of the upper side of Theorem~\texorpdfstring{\ref{theo:ell_q}}{lq} in the nondegenerate case}

By the triangle inequality and Proposition~\ref{prop:surrogate}, we have $\bE\|\hat\bbeta-\bbeta^\star\|_2^2\leq 2\bE\|\hat\bbeta-\tilde\bbeta\|_2^2+2\bE\|\tilde\bbeta-\bbeta^\star\|_2^2\leq 2\sigma_\xi^2 + 2\bE\|\tilde\bbeta-\bbeta^\star\|_2^2$. We are left with the upper bound for $\bE\|\tilde\bbeta-\bbeta^\star\|_2^2$. Fix $1\leq j\leq N$. From Proposition~\ref{prop:implicit_regularization}, $\eta_j$ is $1$-Lipschitz, hence $\bigl|\eta_j(\beta_j^\star+\xi_j)-\eta_j(\beta_j^\star)\bigr|\leq\bigl|\xi_j\bigr|$; combining with $(a+b)^2\leq 2(a^2+b^2)$, we obtain 
    \begin{align}\label{eq:upper_1}
        \bE\bigl|\eta_j(\beta_j^\star+\xi_j) - \beta_j^\star\bigr|^2\leq 2\bigl|\eta_j(\beta_j^\star)-\beta_j^\star\bigr|^2 + 2\bE\bigl|\eta_j(\beta_j^\star+\xi_j)-\eta_j(\beta_j^\star)\bigr|^2\leq 2\bigl|\eta_j(\beta_j^\star)-\beta_j^\star\bigr|^2 + 2\sigma_\xi^2.
    \end{align}On the other hand, notice that $\bigl|\eta_j(\beta_j^\star+\xi_j)-\eta_j(\xi_j)-\beta_j^\star\bigr|\leq\bigl|\beta_j^\star\bigr|$. Indeed, when $\beta_j^\star\geq0$, since $\eta_j$ is increasing, $\eta_j(\beta_j^\star+\xi_j)-\eta_j(\xi_j)\geq 0$, and consequently by $1$-Lipschitzness, $\eta_j(\beta_j^\star+\xi_j)-\eta_j(\xi_j)\leq \beta_j^\star$. When $\beta_j^\star<0$, similarly we have $\beta_j^\star\leq\eta_j(\beta_j^\star+\xi_j)-\eta_j(\xi_j)\leq 0$. Combining the above two cases, we obtain $\bigl|\eta_j(\beta_j^\star+\xi_j)-\eta_j(\xi_j)-\beta_j^\star\bigr|\leq\bigl|\beta_j^\star\bigr|$. As a result,
    \begin{align}\label{eq:upper_2}
        \bE\bigl|\eta_j(\beta_j^\star+\xi_j) - \beta_j^\star\bigr|^2\leq 2\bE\bigl|\eta_j(\xi_j)\bigr|^2 + 2\bE\bigl|\eta_j(\beta_j^\star+\xi_j)-\eta_j(\xi_j)-\beta_j^\star\bigr|^2\leq  2\bE\bigl|\eta_j(\xi_j)\bigr|^2 + 2\bigl|\beta_j^\star\bigr|^2.
    \end{align}Combining \eqref{eq:upper_1}, \eqref{eq:upper_2}, $\tilde\beta_j=\eta_j(\beta_j^\star+\xi_j)$ from Lemma~\ref{lem:coordinate_surrogate}, and summing over $j\in[N]$, we obtain
    \begin{align}\label{eq:upper_tilde_bbeta_estimation_error}
        \bE\bigl\|\tilde\bbeta - \bbeta^\star\bigr\|_2^2\leq 2\sum_{j=1}^N\min\Bigl\{ \bigl|\eta_j(\beta_j^\star)-\beta_j^\star\bigr|^2 + \sigma_\xi^2 ,\,\bE\bigl|\eta_j(\xi_j)\bigr|^2 + \bigl|\beta_j^\star\bigr|^2\Bigr\}.
    \end{align}Using the identity for $r^2(V_{J_\star},V_{J_\star^c})$ from Proposition~\ref{prop:equivalence_r_square_star}, we conclude the proof of the upper side of Theorem~\ref{theo:ell_q}.

\begin{Remark}\label{rmk:Poincare_noise_proof}
Under the Poincar\'e assumption in Theorem~\ref{theo:ell_q}, apply Lemma~\ref{lem:poincare_fiber} to $h=\bar\lambda$ and a median $\lambda_0$ of $\bar\lambda[\vy]$, in place of Lemma~\ref{lem:iso-gaussian-waist}. The deterministic argument in Proposition~\ref{prop:surrogate} then gives $\bE\|\hat\bbeta-\tilde\bbeta\|_2^2\leq2C_{\mathrm P}\sigma_\xi^2$. The Poincar\'e inequality also implies $\bE|\xi_j|^4\leq(1+4C_{\mathrm P})\sigma_\xi^4$; by symmetry and the Paley--Zygmund inequality, $\bP(\xi_j\geq\frac{\sigma_\xi}{\sqrt2})=\bP(\xi_j\leq-\frac{\sigma_\xi}{\sqrt2})\geq\frac{1}{8(1+4C_{\mathrm P})}$. Use these events in place of the Gaussian one-sided events in Lemmas~\ref{lem:upper_min_noise_signal} and~\ref{lem:lower_each_coordinate}. The coordinate moment equivalence $\|\xi_j\|_{L_p}\asymp_{p,C_{\mathrm P}}\sigma_\xi$, for every $p>0$ (see, for instance, \cite[Theorem 1.1.4]{artstein-avidanAsymptoticGeometricAnalysis2021}), together with Lemma~\ref{lem:ell_q_coordinate_estimates}, preserves~\eqref{eq:explicit_ellipsoid_rate}. Thus the same upper and lower bounds hold, with $c,C,C'$ depending on $q$ and $C_{\mathrm P}$.
\end{Remark}

\subsection{The lower bound on the estimation error}
The proof of the lower bound follows the argument in \cite{lecueSharpConvergenceRates2025,lecuesharp2026}.
In this subsection, we prove the following proposition.
\begin{Proposition}\label{prop:lower_surrogate_tilde_beta}
    Under the conditions of Theorem~\ref{theo:ell_q}, for Gaussian noise $\bxi\sim\cN(\vzero,\sigma_\xi^2I_N)$, $\bE\|\tilde\bbeta-\bbeta^\star\|_2^2\geq \frac{1}{704}r^2(V_{J_\star},V_{J_\star^c})$.
\end{Proposition}Once this proposition is proved, combining it with Proposition~\ref{prop:surrogate}, we immediately have
\begin{align}\label{eq:lower_hat_bbeta_small_noise}
\begin{aligned}
\bE\|\tilde\bbeta-\bbeta^\star\|_2^2&\leq2\bE\|\hat\bbeta-\bbeta^\star\|_2^2+2\bE\|\tilde\bbeta-\hat\bbeta\|_2^2,\mbox{ hence }
\bE\|\hat\bbeta-\bbeta^\star\|_2^2\geq\frac{1}{1408}r^2(V_{J_\star},V_{J_\star^c})-\sigma_\xi^2\geq\frac{1}{2816}r^2(V_{J_\star},V_{J_\star^c}),
\end{aligned}
\end{align}
where the last inequality holds whenever $r^2(V_{J_\star},V_{J_\star^c})\geq2816\sigma_\xi^2$.

We first prove the following lemma.
\begin{Lemma}\label{lem:upper_min_noise_signal}
    For any $1\leq j\leq N$, $\min\{|\beta_j^\star|^2,\sigma_\xi^2\}\leq 10(|\beta_j^\star-\eta_j(\beta_j^\star)|^2+\bE|\eta_j(\xi_j)|^2)$.
\end{Lemma}
\begin{proof}[Proof of Lemma~\ref{lem:upper_min_noise_signal}]
    Fix $1\leq j\leq N$. By symmetry of $\xi_j$, replacing $\beta_j^\star$ by $-\beta_j^\star$ and $\eta_j(t)$ by $-\eta_j(-t)$ if necessary, it suffices to consider $\beta_j^\star\geq0$.
    From Proposition~\ref{prop:implicit_regularization}, we know that $\eta_j$ is non-decreasing and $1$-Lipschitz, hence $\bR\ni t\mapsto t-\eta_j(t)\in\bR$ is non-decreasing.
    \begin{enumerate}
        \item When $\sigma_\xi\leq\beta_j^\star$, applying $t\mapsto t-\eta_j(t)$ with $t=\sigma_\xi$ and $t=\beta_j^\star$ respectively yields that $\sigma_\xi-\eta_j(\sigma_\xi)\leq \beta_j^\star-\eta_j(\beta_j^\star)$, hence $\eta_j(\sigma_\xi)\geq \sigma_\xi -|\eta_j(\beta_j^\star)-\beta_j^\star|$.
        \item When $\sigma_\xi>\beta_j^\star$, since $\eta_j$ is non-decreasing, we have $\eta_j(\sigma_\xi)\geq\eta_j(\beta_j^\star)\geq \beta_j^\star - |\eta_j(\beta_j^\star)-\beta_j^\star|$.
    \end{enumerate}Together with the fact that $\eta_j(\sigma_\xi)\geq0$, we obtain
    \begin{align*}
        \eta_j(\sigma_\xi) \geq \Bigl(\min\{\sigma_\xi,\beta_j^\star\} - \bigl|\eta_j(\beta_j^\star)-\beta_j^\star\bigr|\Bigr)_+.
    \end{align*}Let $\Omega_{j,+}=\bigl\{\xi_j\geq \sigma_\xi\bigr\}$. Then $\bP(\Omega_{j,+})\geq\frac19$, and on $\Omega_{j,+}$, $\eta_j(\xi_j)\geq\eta_j(\sigma_\xi)\geq\Bigl(\min\{\sigma_\xi,\beta_j^\star\} - \bigl|\eta_j(\beta_j^\star)-\beta_j^\star\bigr|\Bigr)_+.$ Multiplying $\one_{\{\Omega_{j,+}\}}$ on both sides, and taking expectations, we obtain
    \begin{align*}
        \bE\bigl[\bigl|\eta_j(\xi_j)\bigr|^2\bigr]\geq\Bigl(\min\{\sigma_\xi,\beta_j^\star\} - \bigl|\eta_j(\beta_j^\star)-\beta_j^\star\bigr|\Bigr)_+^2\bE[\one_{\{\Omega_{j,+}\}}] \geq \frac{1}{9}\Bigl(\min\{\sigma_\xi,\beta_j^\star\} - \bigl|\eta_j(\beta_j^\star)-\beta_j^\star\bigr|\Bigr)_+^2.
    \end{align*}Now since
    \begin{align*}
        &\min\bigl\{\sigma_\xi,\beta_j^\star\bigr\} \leq \bigl|\beta_j^\star-\eta_j(\beta_j^\star)\bigr| + \Bigl( \min\bigl\{\sigma_\xi,\beta_j^\star\bigr\} - \bigl| \beta_j^\star - \eta_j(\beta_j^\star) \bigr| \Bigr)_+,
    \end{align*}and since $(a+b)^2\leq10(a^2+\frac{1}{9}b^2)$; squaring both sides, we obtain
    \begin{align*}
        \min\bigl\{\sigma_\xi^2,(\beta_j^\star)^2\bigr\} \leq 10\bigl|\beta_j^\star - \eta_j(\beta_j^\star)\bigr|^2 + \frac{10}{9}\Bigl(\min\{\sigma_\xi,\beta_j^\star\} - \bigl|\eta_j(\beta_j^\star)-\beta_j^\star\bigr|\Bigr)_+^2\leq 10\bigl|\beta_j^\star - \eta_j(\beta_j^\star)\bigr|^2 + 10\bE\bigl|\eta_j(\xi_j)\bigr|^2.
    \end{align*}
\end{proof}

Through Lemma~\ref{lem:upper_min_noise_signal}, we can establish an equivalent expression (up to an absolute multiplicative constant) for $r^2(V_{J_\star},V_{J_\star^c})$.
\begin{Lemma}\label{lem:equivalence_r_square_proof}
    We have $\sum_{j=1}^N\bigl( \bigl|\beta_j^\star-\eta_j(\beta_j^\star)\bigr|^2 + \bE\bigl|\eta_j(\xi_j)\bigr|^2 \bigr)\leq r^2(V_{J_\star},V_{J_\star^c})\leq 16\sum_{j=1}^N\bigl( \bigl|\beta_j^\star-\eta_j(\beta_j^\star)\bigr|^2 + \bE\bigl|\eta_j(\xi_j)\bigr|^2 \bigr).$
\end{Lemma}
\begin{proof}[Proof of Lemma~\ref{lem:equivalence_r_square_proof}]
     Recall that from Proposition~\ref{prop:implicit_regularization}, $\eta_j(0)=0$, and $\eta_j$ is non-decreasing, and is $1$-Lipschitz. Consequently, $\beta_j^\star\eta_j(\beta_j^\star)\geq 0$;  $|\eta_j(\beta_j^\star)|=|\eta_j(\beta_j^\star)-\eta_j(0)|\leq|\beta_j^\star|$; and hence
     \begin{align*}
         |\beta_j^\star-\eta_j(\beta_j^\star)|=\Bigl|\bigl|\beta_j^\star\bigr|-\bigl|\eta_j(\beta_j^\star)\bigr|\Bigr|=\bigl|\beta_j^\star\bigr| - \bigl|\eta_j(\beta_j^\star)\bigr|\leq|\beta_j^\star|.
     \end{align*}Similarly, $\bE\bigl|\eta_j(\xi_j)\bigr|^2 = \bE\bigl|\eta_j(\xi_j)-\eta_j(0)\bigr|^2\leq\bE|\xi_j|^2=\sigma_\xi^2$. Combining, we obtain
     \begin{align*}
         &\bigl|\beta_j^\star-\eta_j(\beta_j^\star)\bigr|^2 + \bE\bigl|\eta_j(\xi_j)\bigr|^2 \leq \min\Big\{ |\beta_j^\star|^2 + \bE\bigl|\eta_j(\xi_j)\bigr|^2,\,   |\beta_j^\star-\eta_j(\beta_j^\star)|^2 + \sigma_\xi^2 \Big\}\\
         &\leq \min\bigl\{\sigma_\xi^2,|\beta_j^\star|^2\bigr\} + \bE\bigl|\eta_j(\xi_j)\bigr|^2 + \bigl|\beta_j^\star-\eta_j(\beta_j^\star)\bigr|^2,
     \end{align*}where we have used $\min\{a+c,b+d\}\leq a+d+\min\{b,c\}$ for any $a,b,c,d\geq0$. Combining Lemma~\ref{lem:upper_min_noise_signal}, we obtain
     \begin{align*}
         \bigl|\beta_j^\star-\eta_j(\beta_j^\star)\bigr|^2 + \bE\bigl|\eta_j(\xi_j)\bigr|^2\leq \min\Bigl\{ \bigl|\beta_j^\star\bigr|^2+\bE\bigl|\eta_j(\xi_j)\bigr|^2, \bigl|\beta_j^\star-\eta_j(\beta_j^\star)\bigr|^2+\sigma_\xi^2 \Bigr\} \leq 11\Bigl(\bE\bigl|\eta_j(\xi_j)\bigr|^2 + \bigl|\beta_j^\star-\eta_j(\beta_j^\star)\bigr|^2\Bigr).
     \end{align*}Summing over $j\in[N]$ and using Proposition~\ref{prop:equivalence_r_square_star} completes the proof.
\end{proof}
The last ingredient needed to prove Proposition~\ref{prop:lower_surrogate_tilde_beta} is the following lemma.
\begin{Lemma}\label{lem:lower_each_coordinate}
    Let $\eta:\bR\to\bR$ be a non-decreasing $1$-Lipschitz function such that $\eta(0)=0$. Then for any $1\leq j\leq N$ such that $\beta_j^\star\geq0$, we have $\bigl|\beta_j^\star-\eta(\beta_j^\star)\bigr|^2 + \bE\bigl|\eta(\xi_j)\bigr|^2\leq 44\bE\bigl|\eta(\beta_j^\star+\xi_j)-\beta_j^\star\bigr|^2$.
\end{Lemma}
\begin{proof}[Proof of Lemma~\ref{lem:lower_each_coordinate}]
    We separate the proof into two steps.
    \begin{enumerate}
        \item We prove that $\bE\bigl|\eta(\beta_j^\star+\xi_j)-\beta_j^\star\bigr|^2\geq\frac{1}{2}\bigl|\beta_j^\star-\eta(\beta_j^\star)\bigr|^2.$

        Let $\Omega_{j,0}=\{\xi_j\leq 0\}$, then $\bP(\Omega_{j,0})=\frac{1}{2}$. On $\Omega_{j,0}$, $\xi_j\leq 0$, hence $\beta_j^\star+\xi_j\leq\beta_j^\star$. Since $\eta$ is non-decreasing and $1$-Lipschitz, $\eta(\beta_j^\star+\xi_j)\leq\eta(\beta_j^\star)\leq\beta_j^\star$. Therefore, $\bigl|\eta(\beta_j^\star+\xi_j)-\beta_j^\star\bigr|\geq\bigl|\eta(\beta_j^\star)-\beta_j^\star\bigr|$. Multiplying $\one_{\{\Omega_{j,0}\}}$ on both sides and taking expectations give
        \begin{align*}
            \bE\bigl|\eta(\beta_j^\star+\xi_j)-\beta_j^\star\bigr|^2\geq \bE\Bigl[\bigl|\eta(\beta_j^\star+\xi_j)-\beta_j^\star\bigr|^2\one_{\{\Omega_{j,0}\}}\Bigr] \geq \bigl|\eta(\beta_j^\star)-\beta_j^\star\bigr|^2\bP(\Omega_{j,0}) = \frac{1}{2}\bigl|\eta(\beta_j^\star)-\beta_j^\star\bigr|^2.
        \end{align*}

        \item We prove that $\bE\bigl|\eta(\beta_j^\star+\xi_j)-\beta_j^\star\bigr|^2\geq \frac{1}{20}\bE\bigl|\eta(\xi_j)\bigr|^2.$

        We first show that 
        \begin{align}\label{eq:pre_lower_coordinate}
            \bE\bigl|\eta(\beta_j^\star+\xi_j)-\beta_j^\star\bigr|^2\geq \frac{1}{9}\min\{\sigma_\xi^2,(\beta_j^\star)^2\}.
        \end{align}
        Let $\Omega_{j,1}=\{\xi_j\leq-\sigma_\xi\}$, then $\bP(\Omega_{j,1})\geq\frac19$. On $\Omega_{j,1}$, $\beta_j^\star+\xi_j\leq\beta_j^\star-\sigma_\xi$, hence $\eta(\beta_j^\star+\xi_j)\leq \eta(\beta_j^\star-\sigma_\xi)\leq(\beta_j^\star-\sigma_\xi)_+$, where we have used the fact that $\eta$ is $1$-Lipschitz, fact that $\eta(0)=0$, and the fact that $\eta$ is non-decreasing. Consequently, $\beta_j^\star - \eta(\beta_j^\star+\xi_j)\geq \beta_j^\star-(\beta_j^\star-\sigma_\xi)_+=\min\{\sigma_\xi,\beta_j^\star\}$ on $\Omega_{j,1}$. Since $\beta_j^\star\geq0$, squaring, multiplying $\one_{\{\Omega_{j,1}\}}$, and taking expectations gives $\bE\bigl|\eta(\beta_j^\star+\xi_j)-\beta_j^\star\bigr|^2\geq \frac{1}{9}\min\{\sigma_\xi^2,(\beta_j^\star)^2\}$, that is, \eqref{eq:pre_lower_coordinate}.

        Now we separate into two cases.
        \begin{enumerate}
            \item When $0\leq\beta_j^\star\leq\sigma_\xi$. Again, using the fact that $\eta$ is $1$-Lipschitz, we obtain that for any $t\in\bR$, $0\leq\eta(\beta_j^\star+t)-\eta(t)\leq\beta_j^\star$, hence $0\leq \beta_j^\star - \bigl( \eta(\beta_j^\star+t)-\eta(t) \bigr)\leq \beta_j^\star$. From $\eta(t) = \bigl(\eta(\beta_j^\star+t)-\beta_j^\star\bigr) + \Bigl( \beta_j^\star - \bigl( \eta(\beta_j^\star+t)-\eta(t) \bigr) \Bigr)$, we obtain
            \begin{align*}
                \bigl|\eta(t)\bigr|&\leq \bigl|\eta(\beta_j^\star+t)-\beta_j^\star\bigr| + \Bigl|\beta_j^\star - \bigl( \eta(\beta_j^\star+t)-\eta(t) \bigr) \Bigr| = \bigl|\eta(\beta_j^\star+t)-\beta_j^\star\bigr| +  \beta_j^\star - \bigl( \eta(\beta_j^\star+t)-\eta(t) \bigr)\\
                &\leq \bigl|\eta(\beta_j^\star+t)-\beta_j^\star\bigr| + \beta_j^\star.
            \end{align*}Setting $t=\xi_j$, squaring and taking expectations, we obtain
            \begin{align*}
                \bE\bigl|\eta(\xi_j)\bigr|^2\leq 2\bE\bigl|\eta(\beta_j^\star+\xi_j)-\beta_j^\star\bigr|^2 + 2(\beta_j^\star)^2.
            \end{align*}By \eqref{eq:pre_lower_coordinate}, together with the assumption $\beta_j^\star\leq\sigma_\xi$, we have $2(\beta_j^\star)^2= 2\min\{\sigma_\xi^2,(\beta_j^\star)^2\}\leq 18\bE\bigl|\eta(\beta_j^\star+\xi_j)-\beta_j^\star\bigr|^2$, that is,
            \begin{align*}
                \bE\bigl|\eta(\xi_j)\bigr|^2\leq 20\bE\bigl|\eta(\beta_j^\star+\xi_j)-\beta_j^\star\bigr|^2.
            \end{align*}

            \item When $\beta_j^\star>\sigma_\xi$. By the $1$-Lipschitz property of $\eta$, we have $|\eta(\xi_j)|\leq|\xi_j|$. Squaring and taking expectations, $\bE|\eta(\xi_j)|^2\leq\sigma_\xi^2$. Therefore, by \eqref{eq:pre_lower_coordinate}, for $\sigma_\xi<\beta_j^\star$,
            \begin{align*}
                \bE|\eta(\xi_j)|^2\leq\sigma_\xi^2=\min\{\sigma_\xi^2,(\beta_j^\star)^2\}\leq 9\bE\bigl|\eta(\beta_j^\star+\xi_j)-\beta_j^\star\bigr|^2.
            \end{align*}
        \end{enumerate}
    \end{enumerate}The two steps give the claimed inequality, since $2+20\leq44$.
\end{proof}
We are now ready to prove Proposition~\ref{prop:lower_surrogate_tilde_beta}.
\begin{proof}[Proof of Proposition~\ref{prop:lower_surrogate_tilde_beta}]
    Recall that $\tilde\bbeta=\Prox_\self[\vy] = \sum_{j=1}^N\eta_j(\beta_j^\star+\xi_j)\ve_j$, see Lemma~\ref{lem:coordinate_surrogate}. Then $\bE\|\tilde\bbeta-\bbeta^\star\|_2^2=\sum_{j=1}^N\bE|\eta_j(\beta_j^\star+\xi_j)-\beta_j^\star|^2$. We only deal with the $j$'s such that $\beta_j^\star\geq0$. Lemma~\ref{lem:lower_each_coordinate} applied to $\eta=\eta_j$ gives
    \begin{align*}
        \forall 1\leq j\leq N,\mbox{ such that }\beta_j^\star\geq0,\qquad \bigl|\beta_j^\star-\eta_j(\beta_j^\star)\bigr|^2 + \bE\bigl|\eta_j(\xi_j)\bigr|^2\leq 44\bE\bigl|\eta_j(\beta_j^\star+\xi_j)-\beta_j^\star\bigr|^2.
    \end{align*}For the indices $j$ such that $\beta_j^\star<0$, we apply the above displayed inequality with $\beta_j^\star$ replaced by $-\beta_j^\star$, $\eta_j(\cdot)$ replaced by $-\eta_j(-\cdot)$. Then $-\beta_j^\star>0$, $-\eta_j(-\cdot)$ is non-decreasing, $1$-Lipschitz and satisfies $-\eta_j(-0)=0$. Then Lemma~\ref{lem:lower_each_coordinate} gives
    \begin{align*}
        \forall 1\leq j\leq N,\mbox{ such that }\beta_j^\star<0,\qquad \bigl|-\beta_j^\star+\eta_j(\beta_j^\star)\bigr|^2 + \bE\bigl|\eta_j(-\xi_j)\bigr|^2\leq 44\bE\bigl|-\eta_j(\beta_j^\star-\xi_j)+\beta_j^\star\bigr|^2.
    \end{align*}Notice that $\xi_j$ is symmetric, hence
    \begin{align*}
        \forall 1\leq j\leq N,\mbox{ such that }\beta_j^\star<0,\qquad \bigl|\beta_j^\star-\eta_j(\beta_j^\star)\bigr|^2 + \bE\bigl|\eta_j(\xi_j)\bigr|^2\leq 44\bE\bigl|\eta_j(\beta_j^\star+\xi_j)-\beta_j^\star\bigr|^2.
    \end{align*}Summing over all $1\leq j\leq N$ together with Lemma~\ref{lem:equivalence_r_square_proof}, we obtain that
    \begin{align*}
        \bE\|\tilde\bbeta-\bbeta^\star\|_2^2\geq\frac{1}{44}\sum_{j=1}^N\Bigl(\bigl|\beta_j^\star-\eta_j(\beta_j^\star)\bigr|^2 + \bE\bigl|\eta_j(\xi_j)\bigr|^2\Bigr) \geq \frac{1}{704}r^2(V_{J_\star},V_{J_\star^c}).
    \end{align*}
\end{proof}

\section{Summary and Outlook}

This work proposes a methodology based on the Feature Space Decomposition method for establishing explicit, instance-specific sharp rates for convex-constrained least squares under some sequence models, and this methodology can, under certain conditions, be extended beyond Gaussian noise. We illustrate this method with the \(\ell_q\) ellipsoid, the \(\ell_1\) ball, and \(kB_1^N\cap B_\infty^N\) as examples. We also provide an elementary proof of the sharp rate for arbitrary signals and arbitrary noise for isotonic regression in the well-specified regime. By combining the reduction method of Han \cite{hanNoisyLinearInverse2023}, we establish explicit, user-friendly matching upper and lower bounds in Gaussian linear regression in the asymptotic regime for a fairly broad range of signals and noise levels. An interesting direction for future research is to extend the sequence model to the random design setting and establish non-asymptotic sharp convergence rates for Gaussian random-design regression analog to \cite{hanNoisyLinearInverse2023}; another interesting direction for future research is to extend the elementary proof for isotonic regression in this paper to multivariate isotonic regression.

\section*{Acknowledgment}
\textbf{LLM Usage.} The author developed the proofs presented in this paper with the assistance of AI-based tools. The AI model used was ChatGPT 6 Astra. The author takes full responsibility for the correctness of the proofs.

\appendix
\section{Proof of Theorem~\texorpdfstring{\ref{theo:isotonic}}{isotonic}}\label{sec:proof_theo_isotonic}

We use, without proof, the following well-known max-min/min-max representation of isotonic regression; this result can be found in, e.g., \cite{barlow_statistical_1972,dai_bias_2020}:
\begin{align}\label{eq:max_min_isotonic_min_max}
    \forall i\in[N],\qquad \hat\beta_i = \max_{1\leq a\leq i}\min_{i\leq b\leq N}\frac{1}{b-a+1}\sum_{j=a}^b y_j = \min_{i\leq b\leq N}\max_{1\leq a\leq i}\frac{1}{b-a+1}\sum_{j=a}^by_j.
\end{align}

\subsection{Upper side}

We first prove the upper side, which is the easy part of the proof of Theorem~\ref{theo:isotonic}. Fix an arbitrary $i\in[N]$; we use the min-max representation \eqref{eq:max_min_isotonic_min_max} and take $b = i+k_{\star,i}^+-1$. By the definition of $k_{\star,i}^+$ provided in \eqref{eq:def_k_star_i_pm}, we know that $i\leq b = i+k_{\star,i}^+-1 \leq N$, hence \eqref{eq:max_min_isotonic_min_max} together with $y_j=\beta_j^\star+\xi_j$ imply that
\begin{align*}
    &\hat\beta_i - \beta_i^\star &\leq \max_{1\leq a\leq i}\frac{1}{i+k_{\star,i}^+-a}\sum_{j=a}^{i+k_{\star,i}^+-1}y_j-\beta_i^\star \leq \max_{1\leq a\leq i}\frac{1}{i+k_{\star,i}^+-a}\sum_{j=a}^{i+k_{\star,i}^+-1}\Bigl(\beta_j^\star-\beta_i^\star\Bigr) + \max_{1\leq a\leq i}\frac{1}{i+k_{\star,i}^+-a}\sum_{j=a}^{i+k_{\star,i}^+-1}\xi_j.
\end{align*}
Now we use the ``local smoothing'' provided by $k_{\star,i}^+$; that is, by the definition of $k_{\star,i}^+$ provided by \eqref{eq:def_k_star_i_pm} and $\beta_1^\star\leq\beta_2^\star\leq\cdots\leq\beta_N^\star$, we have
\begin{align*}
    \forall i\in[N],\, \forall 1\leq a\leq i,\, \forall a\leq j\leq i+k_{\star,i}^+-1,\qquad \beta_j^\star-\beta_i^\star\leq \beta_{i+k_{\star,i}^+-1}^\star-\beta_i^\star\leq \frac{\sigma_\xi}{\sqrt{k_{\star,i}^+}}.
\end{align*}Therefore,
\begin{align*}
    \forall i\in[N],\qquad \hat\beta_i-\beta_i^\star\leq \frac{\sigma_\xi}{\sqrt{k_{\star,i}^+}} + \max_{1\leq a\leq i}\frac{1}{i+k_{\star,i}^+-a}\sum_{j=a}^{i+k_{\star,i}^+-1}\xi_j.
\end{align*}Taking $(\cdot)_+:=\max\{\cdot,0\}$ on both sides together with $(u+v)_+\leq u+|v|$ for any $u\geq 0,\,v\in\bR$; and letting $k=i+k_{\star,i}^+-a$ (then $1\leq a\leq i$ is equivalent to $k_{\star,i}^+\leq k\leq i+k_{\star,i}^+-1$) yield that $\forall i\in[N]$,
\begin{align*}
    (\hat\beta_i-\beta_i^\star)_+\leq \frac{\sigma_\xi}{\sqrt{k_{\star,i}^+}} + \max_{1\leq a\leq i}\bigg|\frac{1}{i+k_{\star,i}^+-a}\sum_{j=a}^{i+k_{\star,i}^+-1}\xi_j\bigg| = \frac{\sigma_\xi}{\sqrt{k_{\star,i}^+}} + \max_{k_{\star,i}^+\leq k\leq i+k_{\star,i}^+-1}\bigg|\frac{1}{k}\sum_{j=i+k_{\star,i}^+-k}^{i+k_{\star,i}^+-1}\xi_j\bigg|.
\end{align*}Since $\xi_1,\cdots,\xi_N$ are independent, the distribution of the above maximum satisfies that
\begin{align*}
    \mathrm{Law}\Biggl(\max_{k_{\star,i}^+\leq k\leq i+k_{\star,i}^+-1}\bigg|\frac{1}{k}\sum_{j=i+k_{\star,i}^+-k}^{i+k_{\star,i}^+-1}\xi_j\bigg| \Biggr) = \mathrm{Law}\Biggl(\max_{k_{\star,i}^+\leq k\leq i+k_{\star,i}^+-1}\bigg|\frac{1}{k}\sum_{j=1}^{k}\xi_j\bigg| \Biggr).
\end{align*}By Lemma~\ref{lem:martingale} applied to $h=k_{\star,i}^+$ and $m=i+k_{\star,i}^+-1$; together with squaring and taking expectations, we obtain that
\begin{align*}
    \forall i\in[N],\qquad \bE\bigl[(\hat\beta_i-\beta_i^\star)_+^2\bigr] \leq 2\frac{\sigma_\xi^2}{k_{\star,i}^+} + \frac{8\sigma_\xi^2}{k_{\star,i}^+} = 10\frac{\sigma_\xi^2}{k_{\star,i}^+}.
\end{align*}Similarly, by taking $a=i-k_{\star,i}^-+1$ in the max-min representation provided by \eqref{eq:max_min_isotonic_min_max} and repeating the above steps, we obtain
\begin{align*}
    \forall i\in[N],\qquad \hat\beta_i-\beta_i^\star\geq -\frac{\sigma_\xi}{\sqrt{k_{\star,i}^-}} - \max_{i\leq b\leq N}\biggl|\frac{1}{b-i+k_{\star,i}^-}\sum_{j=i-k_{\star,i}^-+1}^b\xi_j\biggr|
\end{align*}Taking $(\cdot)_-$ on both sides, reversing the inequality, and repeating the above steps, we obtain
\begin{align*}
    \forall i\in[N],\qquad \bE\bigl[(\hat\beta_i-\beta_i^\star)_-^2\bigr] \leq  10\frac{\sigma_\xi^2}{k_{\star,i}^-}.
\end{align*}Combining the above, and summing over $1\leq i\leq N$, we obtain
\begin{align*}
    \bE\|\hat\bbeta-\bbeta^\star\|_2^2 \leq 10\sigma_\xi^2\sum_{i=1}^N\Bigl(\frac{1}{k_{\star,i}^+} + \frac{1}{k_{\star,i}^-}\Bigr) = 10r_\star^2.
\end{align*}

\subsection{Lower side}
The proof of the lower side is the difficult part of the proof of Theorem~\ref{theo:isotonic}; we divide it into two parts: the stochastic argument and the deterministic argument.

\subsubsection{Stochastic argument}
The main result of this subsubsection is the following lemma.
\begin{Lemma}\label{lem:stochastic_lower_isotonic}
    For any $N\in\bN_+$, any $\sigma_\xi>0$, and any $k_{\star,i}^+ \leq i\leq N$, define
    \begin{align*}
        \Omega_{i,1}:=\Bigl\{ \sum_{j=i-k_{\star,i}^++1}^i\xi_j\geq 3\sigma_\xi\sqrt{k_{\star,i}^+} \Bigr\},\quad\mbox{ and }\quad \Omega_{i,2}:=\bigcap_{t=1}^{N-i}\Bigl\{ \sum_{j=i+1}^{i+t}\xi_j\geq -\sigma_\xi\sqrt{k_{\star,i}^+} - \frac{\sigma_\xi t}{4\sqrt{k_{\star,i}^+}} \Bigr\}.
    \end{align*}Then there exists an absolute constant $c>0$ such that $\min_{1\leq i\leq N}\bP(\Omega_{i,1}\cap\Omega_{i,2})>c$.
\end{Lemma}
\begin{proof}[Proof of Lemma~\ref{lem:stochastic_lower_isotonic}]
    Notice that $\mathrm{Law}(\sum_{j=i-k_{\star,i}^++1}^i\xi_j)=\cN(0,\sigma_\xi^2k_{\star,i}^+)$, then $\bP(\Omega_{i,1})>c_1$ for some absolute constant $c_1$. The second random event holds with probability at least $c_2$ for some $c_2>0$ as a consequence of Lemma~\ref{lem:martingale} with $k=t$, $m=N-i$, $G_j=\frac{\xi_{i+j}}{\sigma_\xi}$, $a=\sqrt{k_{\star,i}^+}$, and $b=\frac{1}{4\sqrt{k_{\star,i}^+}}$. In particular, when $i=N$, $\bP(\Omega_{i,2})=1$. Now since $\Omega_{i,1}$ depends only on $(\xi_j)_{j=i-k_{\star,i}^++1}^i$; while $\Omega_{i,2}$ depends only on $(\xi_j)_{j=i+1}^N$, the two random events are independent, hence $\bP(\Omega_{i,1}\cap\Omega_{i,2})=\bP(\Omega_{i,1})\bP(\Omega_{i,2})\geq c_1c_2>c$ for some absolute constant $c$. Since the above argument is independent of $i$, we obtain $\min_{1\leq i\leq N}\bP(\Omega_{i,1}\cap\Omega_{i,2})>c$.
\end{proof}

\subsubsection{Deterministic argument}

Without loss of generality, we assume that $k_{\star,i}^+\leq k_{\star,i}^-$. In this subsubsection, we fix $i\in[N]$ and work on $\Omega_{i,1}\cap\Omega_{i,2}$. Since we ultimately want a lower bound in expectation, we do not need a union bound.

We first prove the following inequality:
\begin{align}\label{eq:proof_isotonic_inequ}
    \forall 0\leq t\leq N-i,\qquad \sum_{j=i-k_{\star,i}^++1}^{i+t} \bigl(y_j-\beta_i^\star\bigr) \geq \frac{\sigma_\xi(k_{\star,i}^++t)}{4\sqrt{k_{\star,i}^+}}.
\end{align}
\textbf{Proof of \eqref{eq:proof_isotonic_inequ}.} Since we have assumed $k_{\star,i}^+\leq k_{\star,i}^-$, by the definition of $k_{\star,i}^-$ (see \eqref{eq:def_k_star_i_pm}), we have
\begin{align*}
    \forall i-k_{\star,i}^++1\leq j\leq i,\qquad \beta_i^\star-\beta_j^\star \leq \beta_i^\star - \beta_{i-k_{\star,i}^++1}^\star \leq \beta_i^\star - \beta_{i-k_{\star,i}^-+1}^\star \leq \frac{\sigma_\xi}{\sqrt{k_{\star,i}^-}} \leq \frac{\sigma_\xi}{\sqrt{k_{\star,i}^+}}.
\end{align*}
Combining the above inequalities with $y_j=\beta_j^\star+\xi_j$ and $\beta_1^\star\leq\cdots\leq\beta_N^\star$, we obtain: for any $0\leq t\le N-i$,
\begin{align*}
    \sum_{j=i-k_{\star,i}^++1}^{i+t} \bigl(y_j-\beta_i^\star\bigr) &= \sum_{j=i-k_{\star,i}^++1}^{i} \bigl(y_j-\beta_i^\star\bigr) + \sum_{j=i+1}^{i+t} \bigl(y_j-\beta_i^\star\bigr)\\
    &= \sum_{j=i-k_{\star,i}^++1}^{i} \bigl(\beta_j^\star-\beta_i^\star\bigr) + \sum_{j=i+1}^{i+t} \bigl(\beta_j^\star-\beta_i^\star\bigr) + \sum_{j=i-k_{\star,i}^++1}^{i}\xi_j + \sum_{j=i+1}^{i+t}\xi_j\\
    &\geq \sigma_\xi\sqrt{k_{\star,i}^+} - \frac{\sigma_\xi t}{4\sqrt{k_{\star,i}^+}} + \sum_{j=i+1}^{i+t} \bigl(\beta_j^\star-\beta_i^\star\bigr),
\end{align*}where we have used the definition of $\Omega_{i,1}$ and $\Omega_{i,2}$ in the last inequality.
\begin{enumerate}
    \item When $0\leq t<k_{\star,i}^+$. By $\beta_j^\star\geq\beta_i^\star$, we obtain that
    \begin{align*}
        \sum_{j=i-k_{\star,i}^++1}^{i+t} \bigl(y_j-\beta_i^\star\bigr) &\geq \sigma_\xi\sqrt{k_{\star,i}^+} - \frac{\sigma_\xi t}{4\sqrt{k_{\star,i}^+}} = \frac{\sigma_\xi}{4\sqrt{k_{\star,i}^+}}\bigl( 4k_{\star,i}^+-t \bigr) \geq \frac{\sigma_\xi}{4\sqrt{k_{\star,i}^+}}\bigl( k_{\star,i}^++t \bigr),
    \end{align*}where we have used $t<k_{\star,i}^+$ in the last inequality to obtain $4k_{\star,i}^+-t > k_{\star,i}^++t$.

    \item In the following we deal with the case when $k_{\star,i}^+\leq t\leq N-i$. By the definition of $k_{\star,i}^+$ provided in \eqref{eq:def_k_star_i_pm}, we know
    \begin{align*}
        k_{\star,i}^+\leq N-i+1,\quad\mbox{ and }\quad \beta_{i+k_{\star,i}^+}^\star - \beta_i^\star > \frac{\sigma_\xi}{\sqrt{k_{\star,i}^++1}}.
    \end{align*}Therefore, for any $j\geq k_{\star,i}^++i$, by $\beta_1^\star\leq\cdots\leq\beta_N^\star$, we have $\beta_j^\star-\beta_i^\star\geq\beta_{k_{\star,i}^++i}^\star-\beta_i^\star>\frac{\sigma_\xi}{\sqrt{k_{\star,i}^++1}}$. Consequently, when $k_{\star,i}^+\leq t\leq N-i$, by using $\beta_j^\star-\beta_i^\star\geq0$ for any $i+1\leq j\leq i+k_{\star,i}^+-1$, we obtain
    \begin{align*}
        &\sum_{j=i+1}^{i+t} \bigl(\beta_j^\star-\beta_i^\star\bigr) = \Bigl(\sum_{j=i+1}^{i+k_{\star,i}^+-1} + \sum_{j=i+k_{\star,i}^+}^{i+t}\Bigr) \bigl(\beta_j^\star-\beta_i^\star\bigr) \geq \sum_{j=i+k_{\star,i}^+}^{i+t}\bigl(\beta_j^\star-\beta_i^\star\bigr) \geq \frac{\sigma_\xi(t-k_{\star,i}^++1)}{\sqrt{k_{\star,i}^++1}}.
    \end{align*}Consequently,
    \begin{align*}
        \sum_{j=i-k_{\star,i}^++1}^{i+t} \bigl(y_j-\beta_i^\star\bigr) &\geq \sigma_\xi\sqrt{k_{\star,i}^+} - \frac{\sigma_\xi t}{4\sqrt{k_{\star,i}^+}} + \frac{\sigma_\xi(t-k_{\star,i}^++1)}{\sqrt{k_{\star,i}^++1}} \geq \frac{\sigma_\xi(k_{\star,i}^++t)}{4\sqrt{k_{\star,i}^+}},
    \end{align*}where we have used $t\geq k_{\star,i}^+-1$ in the last inequality.
\end{enumerate}
\textbf{The lower bound for $\bE[(\hat\beta_i-\beta_i^\star)_+^2]$.} By the assumption $k_{\star,i}^+\leq k_{\star,i}^-$ and the definition of $k_{\star,i}^-$, we know $1\leq i-k_{\star,i}^++1\leq i$ and $k_{\star,i}^+\leq k_{\star,i}^-\leq i$. By the max-min representation of $\hat\beta_i$ provided by \eqref{eq:max_min_isotonic_min_max}, and taking $a=i-k_{\star,i}^++1$, we obtain $1\leq a\leq i$ and consequently, for any $i\in[N]$ such that $k_{\star,i}^+\leq k_{\star,i}^-$,
\begin{align}\label{eq:lower_difference_isotonic}
    \hat\beta_i-\beta_i^\star \geq \min_{i\leq b\leq N}\frac{1}{b-i+k_{\star,i}^+}\sum_{j=i-k_{\star,i}^++1}^b \bigl( y_j - \beta_i^\star\bigr) = \min_{0\leq t\leq N-i}\frac{1}{t+k_{\star,i}^+}\sum_{j=i-k_{\star,i}^++1}^{i+t}\bigl(y_j-\beta_i^\star\bigr) \geq \frac{\sigma_\xi}{4\sqrt{k_{\star,i}^+}},
\end{align}where we have let $t=b-i$ in the equality, and \eqref{eq:proof_isotonic_inequ} in the last inequality. Similarly, for any $i\in[N]$ such that $k_{\star,i}^-< k_{\star,i}^+$, there exist random events
\begin{align*}
    \Omega_{i,1}^-=\Bigl\{\sum_{j=i}^{i+k_{\star,i}^--1}\xi_j\leq -3\sigma_\xi\sqrt{k_{\star,i}^-}\Bigr\},\quad \Omega_{i,2}^-=\bigcap_{t=1}^{i-1}\Bigl\{ \sum_{j=i-t}^{i-1}\xi_j\leq \sigma_\xi\sqrt{k_{\star,i}^-} + \frac{\sigma_\xi t}{4\sqrt{k_{\star,i}^-}} \Bigr\},
\end{align*}
with $\Omega_{i,1}^-\cap \Omega_{i,2}^-$ strictly positive probability, on which $\hat\beta_i-\beta_i^\star\leq-\frac{\sigma_\xi}{4\sqrt{k_{\star,i}^-}}$.

\textbf{End of the proof of the lower side.} Now square both sides and use \eqref{eq:lower_difference_isotonic} on $\Omega_{i,1}\cap\Omega_{i,2}$; use nonnegativity otherwise. Summing over $1\leq i\leq N$ and separating into $i$ such that $k_{\star,i}^+\leq k_{\star,i}^-$ or the reverse, and using Lemma~\ref{lem:stochastic_lower_isotonic}, we obtain $\bE\|\hat\bbeta-\bbeta^\star\|_2^2\geq c^2\sigma_\xi^2\sum_{i=1}^N(\frac{1}{k_{\star,i}^+}+\frac{1}{k_{\star,i}^-})=c^2r_\star^2$.
\section{Auxiliary lemmas}

\subsection{Convex Analysis}\label{sec:convex_analysis}

The following facts are standard in convex analysis, and may be found, for instance, in \cite[Chapter 12]{bauschkeConvexAnalysisMonotone2011}. For a lower semicontinuous convex function $g:\mathbb R^N\to\bR$ bounded below, its Moreau envelope is defined as $$\Mo_g(\cdot):=\min\Bigl\{\frac{1}{2}\|\cdot-\vx\|_2^2+g(\vx):\vx\in\mathbb R^N\Bigr\},$$ a real-valued function on $\bR^N$; while we recall the definition of the proximal operator of $g$ as
\begin{align*}
    \Prox_g[\bullet]\in\argmin\Bigl(\frac{1}{2}\|\bullet-\vx\|_2^2+g(\vx):\,\vx\in\bR^N\Bigr).
\end{align*}
All subdifferentials of functions defined on $\mathbb R^N$ are understood on $\mathbb R^N$: $\vu\in\partial^-g(\bbeta)$ means $g(\vgamma)-g(\bbeta)\geq\langle\vu,\vgamma-\bbeta\rangle$ for every $\vgamma\in\mathbb R^N$. For the function $\Mo(\bullet_1,\bullet_2)$ used in this work, when $\lambda\in\bR_+$ is fixed, $\Mo(\lambda,\bullet)$ is precisely the Moreau envelope of $\lambda \Psi$ on $\mathbb R^N$.  For each $\vx$, $\lambda\mapsto\Mo(\lambda,\vx)$ is concave; for each $\lambda\geq0$, $\vx\mapsto\Mo(\lambda,\vx)$ is convex. It is easy to see that
\begin{align}\label{eq:identity_Moreau_Proximal}
    \Mo_g(\bullet) = \frac{1}{2}\Bigl\|\bullet-\Prox_g[\bullet]\Bigr\|_2^2 + g\Bigl(\Prox_g[\bullet]\Bigr).
\end{align}It is also known that $\Prox_g$ and $\mathrm{Id}-\Prox_g$ are firmly non-expansive, hence $1$-Lipschitz, \cite[Proposition 12.27 and Proposition 4.2]{bauschkeConvexAnalysisMonotone2011}. In the sense of monotone operators on $\mathbb{R}^N$, $\Prox_g$ is also a nondecreasing operator, \cite[Section 20.1]{bauschkeConvexAnalysisMonotone2011}. Moreover, $\Mo_g$ is Fr\'echet differentiable on $\bR^N$, and $\nabla\Mo_g(\bullet) = \bullet - \Prox_g[\bullet]$, see \cite[Proposition 12.29]{bauschkeConvexAnalysisMonotone2011} (this is in fact a consequence of Danskin's theorem, see Lemma~\ref{lem:danskin}). In particular, for the function $\Mo(\bullet_1,\bullet_2)$ used in this work,
\begin{align}\label{eq:gradient_Moreau_1}
    \forall\lambda\in\bR_+,\quad\forall \vw\in\bR^N,\qquad \nabla_{\bullet_2}\Mo(\lambda,\vw) = \vw-\Prox_\lambda[\vw].
\end{align}
For the squared distance to $T$, these facts give (see \cite[Proposition 12.27, Corollary 12.30]{bauschkeConvexAnalysisMonotone2011})
\begin{align}\label{eq:derivative_Moreau}
    \nabla_\vw\frac{1}{2}\dist^2(\vw,T) = \vw-\Proj_T[\vw],\qquad \forall\vw\in\bR^N,\qquad\mbox{ and }\nabla_\bullet\frac{1}{2}\dist^2(\bullet,T):\bR^N\to\bR^N \mbox{ is }1\mbox{-Lipschitz}.
\end{align}
For any lower semicontinuous convex function $g$ on $\mathbb R^N$, \cite[Proposition 12.26]{bauschkeConvexAnalysisMonotone2011} gives that for any $\vy\in\bR^N$ and $\bbeta\in\mathbb R^N$,
\begin{align}\label{eq:gradient_Moreau}
    \bbeta=\Prox_g[\vy] \quad\Leftrightarrow\quad \forall \vu\in\mathbb R^N,\, g(\vu)-g(\bbeta)\geq\langle \vy-\bbeta,\vu-\bbeta\rangle,\quad\Leftrightarrow\quad \vy-\bbeta\in\partial^-g(\bbeta).
\end{align}In particular, $\bbeta\mapsto\frac{1}{2}\|\vx-\bbeta\|_2^2+\lambda \Psi(\bbeta)$ is $1$-strongly convex on $\mathbb R^N$ for any $\vx\in\bR^N$ and $\lambda\geq0$.
\begin{Lemma}\label{lem:moreau_proximal_equi}
    Under the conditions of Theorem~\ref{theo:ell_q}, for any $\vx\in\bR^N$ and $\lambda\in\bR_+$, the following are equivalent:
    \begin{enumerate}
        \item $\Mo(\lambda,\vx)=\frac{1}{2}\dist^2(\vx,T)$;
        \item $\Prox_\lambda[\vx] = \Proj_T[\vx]$, and $\lambda \Psi(\Proj_T[\vx])=0$;
        \item $\vx-\Proj_T[\vx]\in\partial^-(\lambda \Psi)(\Proj_T[\vx])$, and $\lambda \Psi(\Proj_T[\vx])=0$.
    \end{enumerate}Moreover, for any $\vx\in\bR^N$, at least one $\lambda\geq0$ satisfies these conditions. The minimum-norm choice $\bar\lambda[\vx]$ is unique, and
    \begin{align}\label{eq:identity_Moreau}
        \Mo\bigl(\bar\lambda[\bullet],\bullet\bigr) = \frac{1}{2}\bigl\|\bullet-\Proj_T[\bullet]\bigr\|_2^2 = \frac{1}{2}\dist^2(\bullet,T).
    \end{align}
\end{Lemma}
\begin{proof}[Proof of Lemma~\ref{lem:moreau_proximal_equi}]
    \begin{enumerate}
        \item \emph{1.} $\Rightarrow$ \emph{2.}. Since $\Proj_T[\vx]\in T\subset\mathbb R^N$ and $\lambda \Psi(\Proj_T[\vx])\leq0$,
        \begin{align*}
            \Mo(\lambda,\vx)&\leq\frac{1}{2}\|\vx-\Proj_T[\vx]\|_2^2+\lambda \Psi(\Proj_T[\vx])\leq\frac{1}{2}\|\vx-\Proj_T[\vx]\|_2^2=\frac{1}{2}\dist^2(\vx,T).
        \end{align*}
        Under \emph{1.}, both inequalities are equalities. Thus $\lambda \Psi(\Proj_T[\vx])=0$, and $\Proj_T[\vx]$ minimizes the strongly convex objective defining $\Mo(\lambda,\vx)$. Its unique minimizer is $\Prox_\lambda[\vx]$.

        \item \emph{2.} $\Rightarrow$ \emph{3.}. Apply \eqref{eq:gradient_Moreau} with $g=\lambda \Psi$ and $\vy=\vx$.

        \item \emph{3.} $\Rightarrow$ \emph{1.}. By the definition of subgradient on $\mathbb R^N$, for every $\bbeta\in\mathbb R^N$,
        \begin{align*}
            \frac{1}{2}\|\vx-\bbeta\|_2^2+\lambda \Psi(\bbeta)&\geq\frac{1}{2}\|\vx-\Proj_T[\vx]\|_2^2+\lambda \Psi(\Proj_T[\vx])+\frac{1}{2}\|\bbeta-\Proj_T[\vx]\|_2^2\\
            &=\frac{1}{2}\dist^2(\vx,T)+\frac{1}{2}\|\bbeta-\Proj_T[\vx]\|_2^2.
        \end{align*}
        Taking the minimum over $\bbeta\in\mathbb R^N$, with equality at $\Proj_T[\vx]$, proves \emph{1.}.
    \end{enumerate}
    Since $\Psi$ is continuous and $\Psi(\vzero)<0$, $\vzero$ is strictly feasible. Slater's condition and the KKT theorem therefore provide $\lambda\geq0$ satisfying \emph{3.}. Finally, $\lambda\mapsto\Mo(\lambda,\vx)$ is concave and upper semicontinuous, and is bounded above by $\frac12\dist^2(\vx,T)$. Its nonempty set of maximizers in $\bR_+$ is closed and convex. It consequently has a unique minimum-norm element $\bar\lambda[\vx]$, by \cite[Theorem 3.14]{bauschkeConvexAnalysisMonotone2011}.
\end{proof}
The following Danskin's lemma may be found, for instance, in \cite[Theorem 4.13, p. 273]{bonnansPerturbationAnalysisOptimization2000}.
\begin{Lemma}[Danskin's theorem]\label{lem:danskin}
Let \(E\) be a finite-dimensional Euclidean space, let \(K\) be a nonempty compact set, and let \(\Psi:E\times K\to\mathbb R\) be continuous. Assume that, for every \(\boldsymbol a\in K\), \(\Psi(\cdot,\boldsymbol a)\) is differentiable and that \((\boldsymbol u,\boldsymbol a)\mapsto\nabla_{\boldsymbol u}\Psi(\boldsymbol u,\boldsymbol a)\) is continuous. Then, for every \(\boldsymbol u,\boldsymbol v\in E\),
\begin{align*}
\mathbf D\Big(\boldsymbol u\mapsto\max_{\boldsymbol a\in K}\Psi(\boldsymbol u,\boldsymbol a)\Big)(\boldsymbol u)[\boldsymbol v]=\max_{\boldsymbol a\in\operatorname*{arg\,max}_{\boldsymbol b\in K}\Psi(\boldsymbol u,\boldsymbol b)}\big\langle\nabla_{\boldsymbol u}\Psi(\boldsymbol u,\boldsymbol a),\boldsymbol v\big\rangle,
\end{align*}
and
\begin{align*}
\mathbf D\Big(\boldsymbol u\mapsto\min_{\boldsymbol a\in K}\Psi(\boldsymbol u,\boldsymbol a)\Big)(\boldsymbol u)[\boldsymbol v]=\min_{\boldsymbol a\in\operatorname*{arg\,min}_{\boldsymbol b\in K}\Psi(\boldsymbol u,\boldsymbol b)}\big\langle\nabla_{\boldsymbol u}\Psi(\boldsymbol u,\boldsymbol a),\boldsymbol v\big\rangle.
\end{align*}
If the maximizer is unique, then
\begin{align*}
\nabla_{\boldsymbol u}\max_{\boldsymbol a\in K}\Psi(\boldsymbol u,\boldsymbol a)=\nabla_{\boldsymbol u}\Psi(\boldsymbol u,\boldsymbol a),\qquad \boldsymbol a\in\operatorname*{arg\,max}_{\boldsymbol b\in K}\Psi(\boldsymbol u,\boldsymbol b),
\end{align*}
and if the minimizer is unique, then
\begin{align*}
\nabla_{\boldsymbol u}\min_{\boldsymbol a\in K}\Psi(\boldsymbol u,\boldsymbol a)=\nabla_{\boldsymbol u}\Psi(\boldsymbol u,\boldsymbol a),\qquad \boldsymbol a\in\operatorname*{arg\,min}_{\boldsymbol b\in K}\Psi(\boldsymbol u,\boldsymbol b).
\end{align*}
\end{Lemma}

\subsection{Probabilistic tools}\label{sec:probabilistic_tools}
\begin{Lemma}[\cite{gromovIsoperimetryWaistsConcentration2003}, Theorem 3.1; \cite{klartag_convex_2017}, Theorem 1.1]\label{lem:iso-gaussian-waist}
Let $1\leq m\leq N$, let $\Psi:\mathbb R^N\to\mathbb R^m$ be continuous, and let $\vg,\vh$ be independent, standard Gaussian random vectors on $\bR^N$ and $\bR^m$ respectively. There exists $\boldsymbol t\in \Psi(\mathbb R^N)$ such that $\mathbb P(\operatorname{dist}(\boldsymbol g,\Psi^{-1}(\boldsymbol t))\leq u)\geq\mathbb P(\|\boldsymbol h\|_2\leq u)$ for every $u>0$.
\end{Lemma}
\begin{Lemma}\label{lem:poincare_fiber}
Let $h:\bR^N\to\bR$ be continuous, and let $a$ be a median of $h(\vy)$. If $\bxi$ satisfies the Poincar\'e inequality with parameter $C_{\mathrm P}\sigma_\xi^2$, then
\begin{align*}
    \bE\dist^2\bigl(\vy,h^{-1}(a)\bigr)\leq 2C_{\mathrm P}\sigma_\xi^2.
\end{align*}
\end{Lemma}
\begin{proof}
For $a=\Med(h(\vy))$, let
\begin{align*}
    s_a:\vx\in\bR^N\mapsto s_a(\vx)=\begin{cases}
        \dist(\vx,h^{-1}(a)), & h(\vx)> a,\\
        0, & h(\vx)=a,\\
        -\dist(\vx,h^{-1}(a)), & h(\vx)<a.
    \end{cases}
\end{align*}Then $|s_a(\bullet)|=\dist(\bullet,h^{-1}(a))$. We show that $s_a$ is $1$-Lipschitz. Indeed, for any $\vx,\vx'\in\bR^N$, if $\sign(h(\vx)-a)=\sign(h(\vx')-a)$, then $|s_a(\vx)-s_a(\vx')|=|\dist(\vx,h^{-1}(a))-\dist(\vx',h^{-1}(a))|\leq\|\vx-\vx'\|_2$. When $\sign(h(\vx)-a)\neq\sign(h(\vx')-a)$, without loss of generality, we assume that $h(\vx)>a>h(\vx')$. Since $h$ is continuous, there exists $\vz\in[\vx,\vx']:=\{\alpha\vx+(1-\alpha)\vx':0\leq \alpha\leq 1\}$, such that $h(\vz)=a$, hence $|s_a(\vx)-s_a(\vx')|=\dist(\vx,h^{-1}(a))+\dist(\vx',h^{-1}(a))\leq\|\vx-\vz\|_2+\|\vx'-\vz\|_2=(1-\alpha)\|\vx-\vx'\|_2+\alpha\|\vx-\vx'\|_2=\|\vx-\vx'\|_2$.

Moreover, $\bP(s_a(\vy)\geq0)=\bP(h(\vy)\geq a)\geq\frac12$; while $\bP(s_a(\vy)\leq0)=\bP(h(\vy)\leq a)\geq\frac12$, indicating that $0$ is the median of $s_a(\vy)$. Now, applying the Poincar\'e's inequality to $s_a(\bullet)$, we obtain
\begin{align*}
    \bE\dist^2(\vy,h^{-1}(a)) = \bigl(\bE[s_a(\vy)]\bigr)^2 + \Var\bigl(s_a(\vy)\bigr) \leq 2\Var(s_a(\vy))\leq 2C_{\mathrm{P}}\sigma_\xi^2,
\end{align*}where we have used the mean-median inequality to obtain $\bigl(\bE[s_a(\vy)]\bigr)^2\leq \Var(s_a(\vy))$, that is, applying \cite[Problem 7.11.9, p. 355]{grimmettProbabilityRandomProcesses2020} to $X=\bE[s_a(\vy)]-s_a(\vy)$, $t=\bE[s_a(\vy)]-\varepsilon$ and sending $\varepsilon\to 0^+$ to obtain $\frac{1}{2}\leq\bP(s_a(\vy)\leq 0)\leq \frac{\Var(s_a(\vy))}{\Var(s_a(\vy))+(\bE[s_a(\vy)])^2}$.
\end{proof}

\begin{Lemma}\label{lem:martingale}
    Let $m\in\bN_+$, and let $G_1,\cdots,G_m$ be independent standard Gaussian random variables. For any $1\leq h\leq m$, and any $a,b>0$, the following hold
    \begin{align*}
        \bE\max_{h\leq k\leq m}\Bigl|\frac{1}{k}\sum_{j=1}^k G_j \Bigr|^2 \leq \frac{4}{h},\mbox{ and } \bP\biggl(\bigcap_{k=1}^m\Bigl\{\sum_{j=1}^k G_j + bk>-a \Bigr\}\biggr)\geq1-\exp(-2ab).
    \end{align*}
\end{Lemma}
\begin{proof}[Proof of Lemma~\ref{lem:martingale}]
    It is easy to verify that $\mathsf M_k:=\frac{1}{k}\sum_{j=1}^k G_j$ is a martingale for $(\mathsf M_k)_{k=m}^h$ (in fact, $\mathsf M_k$ is a backward martingale for $k$ from $h$ to $m$; see, e.g., \cite[Section 12.7, p. 499, Example 2]{grimmettProbabilityRandomProcesses2020}), and hence by Doob's maximal inequality in $L^2$ form, as in \cite[Section 12.6]{grimmettProbabilityRandomProcesses2020}, $\mathbb{E}[\max_{h\leq k\leq m}|\mathsf M_k|^2] \leq 4\mathbb{E}|\mathsf M_h|^2=\frac4h$.

    For the second inequality, we consider a one-dimensional Brownian motion $(B_t)_{t\geq 0}$ starting from $0$. Then, by the increment property of Brownian motion, $\mathrm{Law}((B_k)_{k=1}^m) = \mathrm{Law}((\sum_{j=1}^k G_j)_{k=1}^m)$. On the other hand, by the first passage time property of Brownian motion, see \cite[Section 13.4, Corollary 14]{grimmettProbabilityRandomProcesses2020}, $\mathbb{P}(\sup(B_t-bt : t\geq 0)\geq a)=\exp(-2ab)$; combined with $\min(B_k+bk : 1\leq k\leq m) \geq \inf(B_t+bt : t\geq 0)$, we have
    \begin{align*}
        \bP\biggl( \bigcup_{k=1}^m\Bigl\{ \sum_{j=1}^k G_j + bk \leq -a \Bigr\}\biggr) = \bP\biggl( \min_{1\leq k\leq m}\Bigl(B_k+bk\Bigr)\leq-a \biggr) \leq \bP\biggl( \inf_{t\geq 0}\Bigl(B_t+bt\Bigr)\leq-a \biggr)\leq \exp(-2ab).
    \end{align*}
\end{proof}

\subsection{Other lemmas}

\begin{Lemma}\label{lem:existence_l_q}
    Under the assumptions of Theorem~\ref{theo:ell_q}, $\bar\lambda[\bullet]=\|\Sigma^{1/2}(\bullet-\Proj_T[\bullet])\|_{\frac{q}{q-1}}$, hence is unique and $\bar\lambda:\bR^N\to\bR$ is continuous.
\end{Lemma}
\begin{proof}[Proof of Lemma~\ref{lem:existence_l_q}]
    Since $\Psi(\vzero)<0$, Slater's condition and the KKT conditions  give a nonnegative multiplier satisfying 
    \begin{align*}
        \forall j\in [N],\qquad x_j-\hat\beta_j[\vx]=\lambda\sigma_j^{-q/2}|\hat\beta_j[\vx]|^{q-2}\hat\beta_j[\vx],\quad\mbox{ and }\quad \frac{\lambda}{q}\Bigl(\|\Sigma^{-1/2}\hat\bbeta[\vx]\|_q^q-1\Bigr)=0.
    \end{align*}
    Multiplying by $\Sigma^{1/2}$ and taking the $\ell_{q/(q-1)}$ norm gives 
    \begin{align*}
        \|\Sigma^{1/2}(\vx-\hat\bbeta[\vx])\|_{\frac{q}{q-1}}=\lambda\|\Sigma^{-1/2}\hat\bbeta[\vx]\|_q^{q-1}.
    \end{align*}
    If $\lambda>0$, then $\frac{\lambda}{q}\bigl(\|\Sigma^{-1/2}\hat\bbeta[\vx]\|_q^q-1\bigr)=0$ implies that $\|\Sigma^{-1/2}\hat\bbeta[\vx]\|_q^q=1$; if $\lambda=0$, then $\hat\bbeta[\vx]=\vx$. Thus the Lagrangian multiplier $\lambda$ is unique and equals
    \begin{align}\label{eq:ellipsoid_multiplier}
    \bar\lambda[\vx]=\bigl\|\Sigma^{1/2}(\vx-\Proj_T[\vx])\bigr\|_{\frac{q}{q-1}}.
    \end{align}
    It is continuous because $\Proj_T$ is continuous. Moreover, $\bE\hat\lambda<\infty$, since $\|\vy-\Proj_T[\vy]\|_2\leq\|\vy\|_2$ and $\bE\|\vy\|_2^2<\infty$.
\end{proof}
\begin{Lemma}\label{lem:calibration_equation_l_q}
    Under the assumptions of Theorem~\ref{theo:ell_q}, $\lambda_\self>0$ exists uniquely, and satisfies
    \begin{align}\tag{\ref{eq:general_self_calibration}}
    \sum_{j=1}^N\sigma_j^{-q/2}\bE|\eta_j(y_j)|^q=\bE\min\Bigl\{1,\sum_{j=1}^N\sigma_j^{-q/2}|y_j|^q\Bigr\}.
    \end{align}
\end{Lemma}
\begin{proof}[Proof of Lemma~\ref{lem:calibration_equation_l_q}]
We first compute $\frac{\mathrm{d}}{\mathrm{d}\lambda}\Mo(\lambda,\vy)$. Apply Lemma~\ref{lem:danskin} to $E=\bR$, $\vu=\lambda$, $K=B_2(\vzero;\|\vy\|_2)$ (where we have used that $\|\Prox_\lambda[\vy]\|_2\leq\|\vy\|_2$), and $\Psi(u,a)=\frac{1}{2}\|\vy-\bbeta\|_2^2+\lambda \Psi(\bbeta)$; since the minimizer of $\Mo(\lambda,\vy)=\min(\frac{1}{2}\|\vy-\bbeta\|_2^2+\lambda \Psi(\bbeta):\bbeta\in\bR^N)$ is unique and equals $\Prox_\lambda[\vy]$, we have
\begin{align*}
    \frac{\mathrm{d}}{\mathrm{d}\lambda}\Mo(\lambda,\vy)=\frac{\mathrm{d}}{\mathrm{d}\lambda}\Bigl(\frac{1}{2}\|\vy-\bbeta\|_2^2+\lambda \Psi(\bbeta)\Bigr)\Big|_{\bbeta=\Prox_\lambda[\vy]} = \Psi\Bigl(\Prox_\lambda[\vy]\Bigr).
\end{align*}
Recall that
\begin{align}\tag{\ref{eq:def_divergence}}
    d_\lambda(\vy)=\frac12\dist^2(\vy,T)-\Mo(\lambda,\vy)+\lambda \Psi(\hat\bbeta).
\end{align}Then by dominated convergence theorem,
\begin{align*}
    \frac{\mathrm{d}}{\mathrm{d}\lambda}\bE[d_\lambda(\vy)] = \bE\Bigl[\frac{\mathrm{d}}{\mathrm{d}\lambda} d_\lambda(\vy)\Bigr] = \bE \Psi(\hat\bbeta) - \bE\bigl[\Psi(\Prox_\lambda[\vy])\bigr].
\end{align*}Moreover by \eqref{eq:calibration_ell_q}, $\lambda\mapsto \bigl|\eta_{\lambda,j}(y_j)\bigr|$ is almost surely decreasing, hence $\lambda\mapsto\bE\bigl[\Psi(\Prox_\lambda[\vy])\bigr]$ is decreasing, and $\lambda\mapsto \bE[d_\lambda(\vy)]$ is strictly convex, hence $\lambda_\self$ is the unique real number that solves $\frac{\mathrm{d}}{\mathrm{d}\lambda}\bE[d_\lambda(\vy)]=0$. Combining this with $\Prox_{\lambda_\self}[\vy]=\sum_{j=1}^N\eta_j(y_j)\ve_j$ from Proposition~\ref{prop:implicit_regularization}, we obtain
\begin{align*}
    \frac{\mathrm d}{\mathrm d\lambda}\bE d_\lambda(\vy)\bigg|_{\lambda=\lambda_\self}=\frac1q\biggl(\bE\min\Bigl\{1,\sum_{j=1}^N\sigma_j^{-q/2}|y_j|^q\Bigr\}-\sum_{j=1}^N\sigma_j^{-q/2}\bE|\eta_j(y_j)|^q\biggr)=0,
\end{align*}that is
\begin{align}\tag{\ref{eq:general_self_calibration}}
    \sum_{j=1}^N\sigma_j^{-q/2}\bE|\eta_j(y_j)|^q=\bE\min\Bigl\{1,\sum_{j=1}^N\sigma_j^{-q/2}|y_j|^q\Bigr\}.
\end{align}
\end{proof}

\begin{Lemma}\label{lem:ell_q_coordinate_estimates}
Under the conditions of Theorem~\ref{theo:ell_q}.
For any $\lambda\geq0$, $j\in[N]$, and $t\in\bR$, define
    \begin{align*}
        \eta_{\lambda,j}(t):=\argmin\Bigl(\frac{1}{2}(t-s)^2+\frac{\lambda}{q}\sigma_j^{-q/2}|s|^q:\, s\in\bR\Bigr).
    \end{align*}In particular, let $\eta_j=\eta_{\lambda_\self,j}$.
Then for every $j\in[N]$ and $t\in\bR$,
\begin{align}\label{eq:estimate_ell_q_ellipsoid}
\begin{aligned}
\min\bigl\{\frac12,2^{1/(1-q)}\bigr\}\min\Bigl\{|t|,\bigl(\frac{\sigma_j^{q/2}|t|}{\lambda_\self}\bigr)^{1/(q-1)}\Bigr\}&\leq|\eta_j(t)|\leq\min\Bigl\{|t|,\bigl(\frac{\sigma_j^{q/2}|t|}{\lambda_\self}\bigr)^{1/(q-1)}\Bigr\},\\
\min\bigl\{\frac12,2^{1-q}\bigr\}\min\Bigl\{|t|,\frac{\lambda_\self}{\sigma_j^{q/2}}|t|^{q-1}\Bigr\}&\leq|t-\eta_j(t)|\leq\min\Bigl\{|t|,\frac{\lambda_\self}{\sigma_j^{q/2}}|t|^{q-1}\Bigr\}.
\end{aligned}
\end{align}
\end{Lemma}
\begin{proof}
The KKT condition gives
    \begin{align}\label{eq:calibration_ell_q}
    \eta_{\lambda,j}(t)+\lambda\sigma_j^{-q/2}|\eta_{\lambda,j}(t)|^{q-2}\eta_{\lambda,j}(t)=t,\qquad\eta_j=\eta_{\lambda_\self,j}.
    \end{align}
    Since the objective in $\argmin(\frac{1}{2}(t-s)^2+\frac{\lambda}{q}\sigma_j^{-q/2}|s|^q:s\in\bR)$ is strongly convex, $\eta_{\lambda,j}(t)$ is uniquely defined. Moreover, from \eqref{eq:calibration_ell_q}, we have $t = \eta_{\lambda,j}(t)\bigl(1+\lambda\sigma_j^{-q/2}|\eta_{\lambda,j}(t)|^{q-2}\bigr)$, hence $\sign(\eta_j(t))=\sign(t)$, and
\begin{align*}
|t|=|\eta_j(t)|+\lambda_\self\sigma_j^{-q/2}|\eta_j(t)|^{q-1},\qquad|t-\eta_j(t)|=\frac{\lambda_\self}{\sigma_j^{q/2}}|\eta_j(t)|^{q-1}.
\end{align*}
The first identity gives $|\eta_j(t)|\leq|t|$ and $|\eta_j(t)|\leq(\frac{\sigma_j^{q/2}|t|}{\lambda_\self})^{1/(q-1)}$. If $|\eta_j(t)|\geq\frac12|t|$, the lower bound for $|\eta_j(t)|$ follows immediately. Otherwise, $\lambda_\self\sigma_j^{-q/2}|\eta_j(t)|^{q-1}=|t|-|\eta_j(t)|>\frac12|t|$, giving $|\eta_j(t)|\geq2^{1/(1-q)}(\frac{\sigma_j^{q/2}|t|}{\lambda_\self})^{1/(q-1)}$. This proves the first row of \eqref{eq:estimate_ell_q_ellipsoid}.

For the second row of \eqref{eq:estimate_ell_q_ellipsoid}, the identities give $|t-\eta_j(t)|\leq|t|$ and $|t-\eta_j(t)|\leq\frac{\lambda_\self}{\sigma_j^{q/2}}|t|^{q-1}$. If $|\eta_j(t)|\leq\frac12|t|$, then $|t-\eta_j(t)|\geq\frac12|t|$. Otherwise, $|t-\eta_j(t)|=\frac{\lambda_\self}{\sigma_j^{q/2}}|\eta_j(t)|^{q-1}\geq2^{1-q}\frac{\lambda_\self}{\sigma_j^{q/2}}|t|^{q-1}$. The two cases give the required lower bound.
\end{proof}

\begin{Lemma}\label{lem:ell_1_equivalence}
    In Theorem~\ref{theo:ell_1_ball}, we have
    \begin{align*}
        \sum_{j=1}^N\min\Bigl\{\min\bigl\{\bigl|\beta_j^\star\bigr|^2,\lambda_\self^2\bigr\}+\sigma_\xi^2,\, \bigl|\beta_j^\star\bigr|^2+\bE\bigl(|\xi_j|-\lambda_\self\bigr)_+^2\Bigr\} \asymp \sum_{j=1}^N\min\bigl\{\bigl|\beta_j^\star\bigr|^2,\lambda_\self^2\bigr\} + \frac{N\sigma_\xi^2\exp\bigl(-\frac{\lambda_\self^2}{2\sigma_\xi^2}\bigr)}{\bigl(1+\frac{\lambda_\self}{\sigma_\xi}\bigr)^3}.
    \end{align*}
\end{Lemma}
\begin{proof}[Proof of Lemma~\ref{lem:ell_1_equivalence}]
    Applying \cite[Equation (119), (122), (123)]{lecuesharp2026}, one may check that $\bE\bigl(|\xi_j|-\lambda_\self\bigr)_+^2 \asymp \frac{\sigma_\xi^2\exp\bigl(-\frac{\lambda_\self^2}{2\sigma_\xi^2}\bigr)}{\bigl(1+\frac{\lambda_\self}{\sigma_\xi}\bigr)^3}$.
    \begin{enumerate}
        \item When $|\beta_j^\star|\leq\lambda_\self$. Notice that $\bE(|\xi_j|-\lambda_\self)_+^2\leq\sigma_\xi^2$. There holds $\min\{|\beta_j^\star|^2,\lambda_\self^2\}+\sigma_\xi^2 = \bigl|\beta_j^\star\bigr|^2+\sigma_\xi^2\geq\bigl|\beta_j^\star\bigr|^2 + \bE(|\xi_j|-\lambda_\self)_+^2$. Therefore,
        \begin{align*}
            \min\Bigl\{\min\bigl\{\bigl|\beta_j^\star\bigr|^2,\lambda_\self^2\bigr\}+\sigma_\xi^2,\, \bigl|\beta_j^\star\bigr|^2+\bE\bigl(|\xi_j|-\lambda_\self\bigr)_+^2\Bigr\} = \bigl|\beta_j^\star\bigr|^2 + \bE(|\xi_j|-\lambda_\self)_+^2.
        \end{align*}
        
        \item When $|\beta_j^\star|>\lambda_\self$. Then $\min\{|\beta_j^\star|^2,\lambda_\self^2\}=\lambda_\self^2$. From $|\xi_j|\leq\lambda_\self + (|\xi_j|-\lambda_\self)_+$, we have $\sigma_\xi^2\leq 2\lambda_\self^2 + 2\bE(|\xi_j|-\lambda_\self)_+^2$. Then
        \begin{align*}
            \lambda_\self^2 + \bE\bigl(|\xi_j|-\lambda_\self\bigr)_+^2 \leq \min\Bigl\{\lambda_\self^2+\sigma_\xi^2,\, \bigl|\beta_j^\star\bigr|^2+\bE\bigl(|\xi_j|-\lambda_\self\bigr)_+^2\Bigr\}\leq \lambda_\self^2+\sigma_\xi^2 \leq 3\bigl(\lambda_\self^2 + \bE\bigl(|\xi_j|-\lambda_\self\bigr)_+^2\bigr).
        \end{align*}
    \end{enumerate}Summing over $1\leq j\leq N$ together with the above estimate for $\bE(|\xi_j|-\lambda_\self)_+^2$, we finish the proof. 
\end{proof}

\subsection{Proof of Corollary~\texorpdfstring{\ref{coro:asymptotic_Gaussian}}{Han}}\label{sec:proof_coro_han}

Take $\{\vG^{(N)}\}_{N\geq 2}$ where $\vG^{(N)}\sim\cN(\vzero,I_N)$. In the following we omit the superscript $(N)$. Let
\begin{align*}
    \cT(T,\bbeta^\star) = \Bigl\{ \vh\in\bR^N:\, h_i\leq h_{i+1},\quad \forall i\mbox{ such that } \beta_i^\star=\beta_{i+1}^\star  \Bigr\}.
\end{align*}It is shown in \cite[Proposition 3.1]{bellecSharpOracleInequalities2018} that $\cT(T,\bbeta^\star)$ is precisely the tangent cone of $T$ evaluated at $\bbeta^\star$, and
\begin{align}\label{eq:statistical_dimension_isotonic}
    \mathbb E\|\mathrm{Proj}_{\mathcal T(T,\boldsymbol\beta^\star)}(\boldsymbol G)\|_2^2 =\sum_{\ell=1}^K\sum_{j=1}^{i_\ell-i_{\ell-1}}\frac{1}{j}.
\end{align}
To prove Corollary~\ref{coro:asymptotic_Gaussian}, we only need to use \cite[Theorem 2.7]{hanNoisyLinearInverse2023}\footnote{To obtain this normalization, apply Han's result to the design vectors $X_i/\sqrt N$, the signal $\sqrt N\boldsymbol\beta^\star$, and the constraint $\sqrt N T=T$.}. In the following we use the notation from \cite{hanNoisyLinearInverse2023} in our language. Let $\mathrm{err}(t)=\big\|\mathrm{Proj}_T(\sqrt N\boldsymbol\beta^\star+t\boldsymbol G)-\sqrt N\boldsymbol\beta^\star\big\|_2^2$, and
\begin{align*}
\mathrm{lrt}(t) = \bigl\|\sqrt{N}\bbeta^\star+t\vG-\sqrt{N}\bbeta^\star\bigr\|_2^2 - \bigl\|\sqrt{N}\bbeta^\star + t\vG - \Proj_T\bigl[ \sqrt{N}\bbeta^\star+t\vG \bigr]\bigr\|_2^2 =t^2\|\boldsymbol G\|_2^2-\mathrm{dist}^2(\sqrt N\boldsymbol\beta^\star+t\boldsymbol G,T).
\end{align*}The noise transformation of \cite{hanNoisyLinearInverse2023} is $\omega_{n/N}(r)=\sqrt{\frac{N(r^2+\sigma_\xi^2)}{n}}$. Let $r_N>0$ be the solution of the fixed point equation 
\begin{align}\label{eq:fixed_point_equation_Han}
    \frac1N\mathbb E\,\mathrm{err}(\omega_{n/N}(r_N))=\frac{1}{N}\bE\mathrm{err}\Bigl(\sqrt{\frac{N(r_N^2+\sigma_\xi^2)}{n}}\Bigr)=r_N^2,
\end{align}
and let
\begin{align*}
    \omega_N=\omega_{n/N}(r_N) = \sqrt{\frac{N(r_N^2+\sigma_\xi^2)}{n}},\quad\mbox{ then }\quad \bE\mathrm{err}(\omega_N) = Nr_N^2.
\end{align*}
Furthermore, let
\begin{align*}
\bar r_N^2=\frac1N\mathbb E\,\mathrm{err}\big(\omega_{n/N}(0)\big) = \frac{1}{N}\bE\mathrm{err}\Bigl(\sigma_\xi\sqrt{\frac{N}{n}}\Bigr),\quad\mbox{ and }\quad \mathfrak L_N=\log\Big(1+\sum_{\ell=1}^K\sum_{j=1}^{i_\ell-i_{\ell-1}}\frac1j\Big)+\log\log(16N).
\end{align*}

\begin{enumerate}
\item We compute Equation 2.14 of \cite{hanNoisyLinearInverse2023} and verify that it tends to zero. By positive homogeneity of projection onto a closed convex cone,
\begin{align*}
\bar r_N^2&=\frac1N\mathbb E\big\|\mathrm{Proj}_T(\sqrt N\boldsymbol\beta^\star+\sigma_\xi\sqrt{N/n}\boldsymbol G)-\sqrt N\boldsymbol\beta^\star\big\|_2^2 =\mathbb E\big\|\mathrm{Proj}_T(\boldsymbol\beta^\star+\frac{\sigma_\xi}{\sqrt n}\boldsymbol G)-\boldsymbol\beta^\star\big\|_2^2,
\end{align*}which is precisely the expected estimation error of the isotonic regression at noise level $\frac{\sigma_\xi}{\sqrt{n}}$.
By Theorem~\ref{theo:isotonic}, there exists an absolute constant $c>0$ such that $2cr_\star^2\leq \bE\|\Proj_T[\vy]-\bbeta^\star\|_2^2\leq 10r_\star^2$, hence $\bar r_N^2\asymp r_\star^2$. Thus item~\emph{1.} of Assumption~\ref{assumption:han}, that is, $r_\star^2\to0$, verifies Equation 2.14 of \cite{hanNoisyLinearInverse2023}. Moreover, $\mathfrak L_N\asymp\log(eK)+\log\log(16N)$. Item~\emph{2.} of Assumption~\ref{assumption:han}, that is, $\log(eK)+\log\log(16N)=o(nr_\star^2)$, yields both $n\gg \mathfrak L_N$ and $\mathfrak L_N/n=o(\bar r_N^2)$.

\item We verify (R2) of \cite{hanNoisyLinearInverse2023} (see Equation 2.9 therein). By \eqref{eq:fixed_point_equation_Han}, together with \cite[Proposition 4.2]{hanNoisyLinearInverse2023}, for any $b\geq 1$ and $t>0$, $\mathrm{err}(t)\leq \mathrm{err}(bt)\leq b^2\mathrm{err}(t)$. Applying this equivalence together with the definition of $\bar r_N^2$ and that of $r_N^2$, we obtain $r_N^2\asymp r_\star^2\to0$. We consider the two alternatives in item~\emph{3.} of Assumption~\ref{assumption:han} separately.

\begin{enumerate}
\item Suppose the first condition of item~\emph{3.} of Assumption~\ref{assumption:han} holds, that is, $r_\star^2\inf(\|\bbeta^\star-a\one_N\|_2^2:a\in\bR)\to0$. 
Recall that by definition,
\begin{align*}
    \dist^2(\sqrt{N}\bbeta^\star+\omega_N\vG,T) = \Bigl\|\sqrt{N}\bbeta^\star + \omega_N\vG\Bigr\|_2^2 - \Bigl\|\Proj_T\bigl[\sqrt{N}\bbeta^\star+\omega_N\vG\bigr]\Bigr\|_2^2.
\end{align*}Then opening the square and using $\bE\vG=\vzero$, we obtain
\begin{align*}
    &\bE\mathrm{lrt}(\omega_N) = \bE\biggl[\omega_N^2\|\vG\|_2^2 - \Bigl\|\sqrt{N}\bbeta^\star + \omega_N\vG\Bigr\|_2^2 + \Bigl\|\Proj_T\bigl[\sqrt{N}\bbeta^\star+\omega_N\vG\bigr]\Bigr\|_2^2\biggr]\\
    &=\bE\Bigl\|\Proj_T\bigl[\sqrt{N}\bbeta^\star+\omega_N\vG\bigr]\Bigr\|_2^2 - N\bigl\|\bbeta^\star\bigr\|_2^2.
\end{align*}Moreover, opening the square again,
\begin{align*}
    &\bE\mathrm{err}(\omega_N) = \bE\Bigl\| \Proj_T\bigl[ \sqrt{N}\bbeta^\star + \omega_N\vG \bigr] \Bigr\|_2^2 - 2\sqrt{N}\Bigl\langle\bE\Proj_T\bigl[\sqrt{N}\bbeta^\star+\omega_N\vG\bigr], \bbeta^\star\Bigr\rangle + N\bigl\|\bbeta^\star\bigr\|_2^2.
\end{align*}
Notice that $a\one_N\in T$, hence the projection residual is orthogonal to $a\one_N$ (see, for instance, \cite[Chapter 28]{bauschkeConvexAnalysisMonotone2011}), that is,
\begin{align*}
    \Bigl\langle \bigl(\sqrt{N}\bbeta^\star + \omega_N\vG\bigr) - \Proj_T\bigl[\sqrt{N}\bbeta^\star + \omega_N\vG\bigr],\, \one_N \Bigr\rangle=0.
\end{align*}Taking expectation on both sides together with $\bE\vG=\vzero$, we obtain
\begin{align*}
    \Bigl\langle \bE\Proj_T\bigl[ \sqrt{N}\bbeta^\star + \omega_N\vG \bigr] - \sqrt{N}\bbeta^\star,\,\one_N\Bigr\rangle=0.
\end{align*}
Notice that $-2N\|\bbeta^\star\|_2^2=2\langle -\sqrt{N}\bbeta^\star,\sqrt{N}\bbeta^\star\rangle$, we obtain
\begin{align*}
&\mathbb E\,\mathrm{lrt}(\omega_N)-\mathbb E\,\mathrm{err}(\omega_N) =2\big\langle\mathbb E\,\mathrm{Proj}_T(\sqrt N\boldsymbol\beta^\star+\omega_N\boldsymbol G)-\sqrt N\boldsymbol\beta^\star,\sqrt N(\boldsymbol\beta^\star-a\boldsymbol1_N)\big\rangle.
\end{align*}
Consequently, Cauchy's inequality applied twice together with $\mathbb E\,\mathrm{err}(\omega_N)=Nr_N^2$ from \eqref{eq:fixed_point_equation_Han} imply
\begin{align*}
    &\frac{1}{2}\Bigl( \bE\mathrm{lrt}(\omega_N) - \bE\mathrm{err}(\omega_N) \Bigr) \leq \bE\bigl\| \Proj_T\bigl(\sqrt{N}\bbeta^\star+\omega_N\vG\bigr) - \sqrt{N}\bbeta^\star \bigr\|_2 \sqrt{N}\bigl\|\bbeta^\star-a\one_N\bigr\|_2\\
    &\leq \Bigl(\bE\bigl\| \Proj_T\bigl(\sqrt{N}\bbeta^\star+\omega_N\vG\bigr) - \sqrt{N}\bbeta^\star \bigr\|_2^2\Bigr)^{1/2} \sqrt{N}\bigl\|\bbeta^\star-a\one_N\bigr\|_2 = \sqrt{N}r_N\sqrt{N}\bigl\|\bbeta^\star-a\one_N\bigr\|_2.
\end{align*}
Therefore, $\frac{\mathbb E\,\mathrm{lrt}(\omega_N)-\mathbb E\,\mathrm{err}(\omega_N)}{2N\sigma_\xi^2}\le\frac{r_N}{\sigma_\xi^2}\inf_{a\in\mathbb R}\|\boldsymbol\beta^\star-a\boldsymbol1_N\|_2\to0,$ where we have used that $r_N^2\leq 20r_\star^2$, together with the assumption that $r_\star^2\inf_{a\in\mathbb R}\|\boldsymbol\beta^\star-a\boldsymbol1_N\|_2^2\to0$.

\item Suppose the second condition of item~\emph{3.} of Assumption~\ref{assumption:han} holds, that is, $\limsup\frac{1}{n}\sum_{\ell=1}^K\sum_{j=1}^{i_\ell-i_{\ell-1}}\frac{1}{j}<2$. Then by definition and opening square,
\begin{align*}
    \mathrm{lrt}(\omega_N)&= \omega_N^2\|\vG\|_2^2 - \inf_{\bbeta\in T}\bigl\|\sqrt{N}\bbeta^\star + \omega_N\vG - \bbeta\bigr\|_2^2 = \sup_{\bbeta\in T}\Bigl( 2\omega_N\langle\boldsymbol G,\boldsymbol\beta-\sqrt N\boldsymbol\beta^\star\rangle-\|\boldsymbol\beta-\sqrt N\boldsymbol\beta^\star\|_2^2 \Bigr).
\end{align*}
Since $T-\sqrt N\boldsymbol\beta^\star\subset\mathcal T(T,\boldsymbol\beta^\star)$ (we have used that for any $i$ such that $\beta_i^\star=\beta_{i+1}^\star$, $\beta_i-\sqrt{N}\beta_i^\star\leq\beta_{i+1}-\sqrt{N}\beta_i^\star=\beta_{i+1}-\sqrt{N}\beta_{i+1}^\star$), the supreme over $\{\bbeta-\sqrt{N}\bbeta^\star:\bbeta\in T\}$ may be enlarged to $\cT(T,\bbeta^\star)$, hence
\begin{align*}
\mathrm{lrt}(\omega_N)&\leq \sup\Bigl( 2\omega_N\langle\vG,\vh\rangle - \|\vh\|_2^2:\, \vh\in\cT(T,\bbeta^\star) \Bigr) = \omega_N^2\|\vG\|_2^2 - \dist^2\Bigl( \omega_N\vG, \cT(T,\bbeta^\star) \Bigr) \\
&=\bigl\|\Proj_{\cT(T,\bbeta^\star)}[\omega_N\vG]\bigr\|_2^2\le\omega_N^2\|\mathrm{Proj}_{\mathcal T(T,\boldsymbol\beta^\star)}(\boldsymbol G)\|_2^2.
\end{align*}
Combining with the definition $\omega_N^2=\frac{N(r_N^2+\sigma_\xi^2)}n$, and $\mathbb E\,\mathrm{err}(\omega_N)=Nr_N^2$ from \eqref{eq:fixed_point_equation_Han}, we obtain
\begin{align*}
\frac{\mathbb E\,\mathrm{lrt}(\omega_N)-\mathbb E\,\mathrm{err}(\omega_N)}{2N\sigma_\xi^2}\le\frac{r_N^2+\sigma_\xi^2}{2n\sigma_\xi^2}\sum_{\ell=1}^K\sum_{j=1}^{i_\ell-i_{\ell-1}}\frac1j-\frac{r_N^2}{2\sigma_\xi^2}.
\end{align*}
Since $r_N^2\to0$, the second assumption ensures that the $\limsup$ of the right-hand side is strictly smaller than $1$. Thus (R2) holds.
\end{enumerate}

\item We substitute the sequence-model bounds into the conclusion of the proposition. Since $\mathfrak L_N/n=o(\bar r_N^2)$, $\frac{\|\widehat{\boldsymbol\beta}-\boldsymbol\beta^\star\|_2^2}{\bar r_N^2}\to 1$ in probability and hence the corollary is proved.
\end{enumerate}

\subsubsection{Examples}\label{sec:examples_revisited}
\begin{enumerate}

\item Example~\ref{example:uniform_increase} (\textbf{revisited}). In this case, $K=N$. Repeating similar computation as in Example~\ref{example:uniform_increase}, we obtain that
\begin{align*}
    k_{\star,1}^+\asymp\min\Bigl\{N,\max\bigl\{1,\bigl(\frac{\sigma_\xi}{v\sqrt{n}}\bigr)^{2/3}\bigr\}\Bigr\},\quad k_{\star,i}^-=\min\{i,k_{\star,1}^+\},\quad\mbox{ and }\quad k_{\star,i}^+=\min\{N-i+1,k_{\star,1}^+\}.
\end{align*}
Therefore, $r_\star^2\asymp\frac{\sigma_\xi^2}{n}\Big[\log N+N\min\Big\{1,\Big(\frac{v\sqrt n}{\sigma_\xi}\Big)^{2/3}\Big\}\Big]$. Moreover, since $\|\bbeta^\star-a\one_N\|_2^2=\sum_{i=1}^N(\beta_1^\star+v(i-1)-a)^2$, one may compute that
\[
\inf_{a\in\mathbb R}\|\boldsymbol\beta^\star-a\boldsymbol{1}_N\|_2^2=\frac{v^2N(N^2-1)}{12}.
\]
As a result, item~\emph{1.} and item~\emph{2.} of Assumption~\ref{assumption:han} are equivalent to
\[
\frac{\log N+N\min\{1,(v\sqrt n/\sigma_\xi)^{2/3}\}}{n}\to0,\quad\mbox{ and }\quad \frac{v\sqrt n}{\sigma_\xi}\Big(\frac{N}{\log N}\Big)^{3/2}\to\infty.
\]
Item~\emph{3.} of Assumption~\ref{assumption:han} is equivalent to either $r_\star^2v^2N(N^2-1)/12\to0$ or $\limsup N/n<2$.

For $v=\sigma_\xi/N$, the condition $\limsup N/n<2$ suffices for all the assumptions, and
\[
r_\star^2\asymp\sigma_\xi^2\min\Big\{\frac{N}{n},\frac{N^{1/3}}{n^{2/3}}\Big\}.
\]
For $v=\sigma_\xi/\sqrt N$, the condition $N=o(n)$ suffices, and $r_\star^2\asymp\sigma_\xi^2N/n$.

\item Example~\ref{example:K_plateaus} (\textbf{revisited}). Suppose $K$ divides $N$ and $\beta_i^\star=\sigma_\xi\lfloor K(i-1)/N\rfloor$. By \eqref{eq:statistical_dimension_isotonic} and the same computation as in Example~\ref{example:K_plateaus}, $k_{\star,i}^-=i-\frac{(\ell-1)N}{K}$, and $k_{\star,i}^+=\frac{\ell N}{K}-i+1$. Therefore,
\[
\mathbb E\|\mathrm{Proj}_{\mathcal T(T,\boldsymbol\beta^\star)}(\boldsymbol G)\|_2^2=K\sum_{j=1}^{N/K}\frac{1}{j},\quad\mbox{ and }\quad  r_\star^2=\frac{2\sigma_\xi^2K}{n}\sum_{j=1}^{N/K}\frac{1}{j}\asymp\frac{\sigma_\xi^2K}{n}\log\Big(\frac{eN}{K}\Big).
\]
Thus $K\log(eN/K)=o(n)$ verifies item~\emph{1.}, $\limsup\frac{1}{n}\sum_{\ell=1}^K\sum_{j=1}^{i_\ell-i_{\ell-1}}\frac{1}{j} = \limsup \frac{1}{n}\sum_{\ell=1}^K \sum_{j=1}^{N/K}\frac{1}{j}=\frac{K}{n}\sum_{j=1}^{N/K}\frac{1}{j}\asymp \frac{K}{n}\log(\frac{eN}{K})<2$. Lastly, item~\emph{2.} follows from $\frac{\log(eK)+\log\log(16N)}{K\log(eN/K)}\to0$.

\item Example~\ref{example:plateau_transition_plateau} (\textbf{revisited}). In this case, $K=m+1$. The leftest plateau has length $N/3$, the middle has length $m-1$; while the rightest has $\frac{2}{3}N-m+1$. Hence $\limsup\frac1n\sum_{\ell=1}^K\sum_{j=1}^{i_\ell-i_{\ell-1}}\frac{1}{j}=\frac1n\bigl( \sum_{j=1}^{N/3}\frac{1}{j} + (m-1) + \sum_{j=1}^{(2/3)N-m+1}\frac{1}{j} \bigr)\asymp \frac{m}{n}+O(\frac{\log(N)}{n})$. By assumption, $\log(N)=o(N)$, and $\limsup\frac{m}{n}<2$, hence item~\emph{3.} of Assumption~\ref{assumption:han} is verified. In the following we repeat the computation of Example~\ref{example:plateau_transition_plateau}. (i) when $1\leq i\leq \frac{1}{3}N$, then $\beta_i^\star=0$, $k_{\star,i}^-=i$, and for any $1\leq k\leq \frac{1}{3}N-i+1$, $i\leq i+k-1\leq\frac{1}{3}N$, hence $\beta_{i+k-1}^\star-\beta_i^\star=0$; hence $k_{\star,i}^+\geq \frac{1}{3}N-i+1$. Therefore, $\sum_{i=1}^{N/3}\frac{1}{i}\leq \sum_{i=1}^{N/3}\bigl(\frac{1}{k_{\star,i}^+}+\frac{1}{k_{\star,i}^-}\bigr) \leq \sum_{i=1}^{N/3}\bigl(\frac{1}{i}+\frac{1}{\frac{1}{3}N-i+1}\bigr)=2\sum_{i=1}^{N/3}\frac{1}{i}\asymp \log(N)$. (ii) when $\frac{1}{3}N+m\leq i\leq N$, similarly, $\sum_{i=\frac{1}{3}N+m}^N\bigl(\frac{1}{k_{\star,i}^+}+\frac{1}{k_{\star,i}^-}\bigr)\asymp\log(N)$. (iii) when $\frac{1}{3}N<i<\frac{1}{3}N+m$, we estimate $k_{\star,i}^\mp$ in the following.
\begin{enumerate}
    \item We estimate $k_{\star,i}^-$. Fix any $1\leq k\leq i$. If $i-k+1>\frac{1}{3}N$, then $\beta_i^\star-\beta_{i-k+1}^\star=\frac{\sigma_\xi}{m}(k-1)$; if $i-k+1\leq\frac{1}{3}N$, then $\beta_i^\star-\beta_{i-k+1}^\star=\frac{\sigma_\xi}{m}(i-\frac{1}{3}N)$. In summary, $\beta_i^\star-\beta_{i-k+1}^\star=\frac{\sigma_\xi}{m}\min\bigl\{k-1,i-\frac{1}{3}N\bigr\}$. By the definition of $k_{\star,i}^-$ from \eqref{eq:def_k_star_i_pm},
    \begin{align}\label{eq:k_star_i_minus_plateau_transition_plateau}
        k_{\star,i}^-=\max\Bigl\{ 1\leq k\leq i:\, \min\bigl\{ k-1,i-\frac{1}{3}N \bigr\} \sqrt{nk} \leq m \Bigr\}.
    \end{align}

    \item We estimate $k_{\star,i}^+$. Fix any $1\leq k\leq N-i+1$. One may repeat the above analysis for $i+k-1<\frac{1}{3}N+m$ and $i+k-1\geq \frac{1}{3}N+m$ separately, to obtain  $\beta_{i+k-1}^\star-\beta_i^\star=\frac{\sigma_\xi}{m}\min\bigl\{k-1,\frac{1}{3}N+m-i\bigr\}$. Again, by the definition of $k_{\star,i}^+$ from \eqref{eq:def_k_star_i_pm},
    \begin{align*}
        k_{\star,i}^+ = \max\Bigl\{ 1\leq k\leq N-i+1:\, \min\bigl\{ k-1,\frac{1}{3}N+m-i \bigr\} \sqrt{nk}\leq m \Bigr\}.
    \end{align*}
\end{enumerate}Since $(k-1)\sqrt{nk}\leq \sqrt{n}k^{3/2}$, we obtain that
\begin{align*}
    k_{\star,i}^\pm\geq\min\{1,\lfloor (\frac{m^2}{n})^{1/3}\rfloor\}.
\end{align*}
We next estimate the upper bounds for $k_{\star,i}^\pm$. From $(i-\frac{1}{3}N) + (\frac{1}{3}N+m-i)=m$, we know $\min\{i-\frac{1}{3}N,\frac{1}{3}N+m-i\}\geq\frac{1}{2}m$. We separate into two cases in the following.
\begin{enumerate}
    \item When $i-\frac{1}{3}N\geq\frac12m$. We estimate $k_{\star,i}^-$ from above. Recall that constraint in \eqref{eq:k_star_i_minus_plateau_transition_plateau}, that is, $\min\{k-1,i-\frac{1}{3}N\}\sqrt{nk}\leq m$.
    \begin{enumerate}
        \item If $k-1\leq i-\frac{1}{3}N$ and $k\geq 2$, then the constraint becomes $(k-1)\sqrt{nk}\leq m$, and $k-1\geq k/2$. Therefore, $\frac{\sqrt{n}}{2}k^{3/2}\leq (k-1)\sqrt{nk}\leq m$, indicating that $k\leq\bigl(\frac{2m}{\sqrt{n}}\bigr)^{2/3}$.

        \item If $k-1>i-\frac{1}{3}N$. Then the constraint becomes $(i-\frac{1}{3}N)\sqrt{nk}\leq m$, hence $\frac{m}{2}\sqrt{nk}\leq (i-\frac{1}{3}N)\sqrt{nk}\leq m$, indicating that $k\leq 4$.
    \end{enumerate}In summary, any $k$ satisfying the constraint in the definition of $k_{\star,i}^-$ in \eqref{eq:k_star_i_minus_plateau_transition_plateau} must satisfy $k\leq \max\{4,(\frac{4m^2}{n})^{1/3}\}\leq 4\max\{1,(\frac{m^2}{n})^{1/3}\}$. Therefore,
    \begin{align*}
        k_{\star,i}^-\leq 4\max\Bigl\{ 1, \Bigl(\frac{m^2}{n}\Bigr)^{1/3} \Bigr\}.
    \end{align*}

    \item When $\frac{1}{3}N+m-1\geq\frac m2$. In this case, $\frac{1}{3}N+m-i>\frac{1}{2}m$. Repeating the computation as in the previous case, we obtain
    \begin{align*}
        k_{\star,i}^+\leq 4\max\Bigl\{ 1, \Bigl(\frac{m^2}{n}\Bigr)^{1/3} \Bigr\}.
    \end{align*}
\end{enumerate}
In summary, for any $\frac{1}{3}N<i<\frac{1}{3}N+m$,
\begin{align*}
    \min\Bigl\{k_{\star,i}^-,k_{\star,i}^+\Bigr\}\leq 4\max\Bigl\{ 1, \Bigl(\frac{m^2}{n}\Bigr)^{1/3} \Bigr\}.
\end{align*}
Therefore, for any such $i$, $\frac{1}{k_{\star,i}^+}+\frac{1}{k_{\star,i}^-}\asymp \min\{1,(\frac{n}{m^2})^{1/3}\}$. Taking sum over $\frac{1}{3}N<i<\frac{1}{3}N+m$, we obtain
\begin{align*}
    \frac{\sigma_\xi^2}{n}\sum_{i=\frac{1}{3}N+1}^{\frac{1}{3}N+m-1} \Bigl( \frac{1}{k_{\star,i}^+}+\frac{1}{k_{\star,i}^-} \Bigr) \asymp \frac{\sigma_\xi^2}{n}(m-1)  \min\{1,(\frac{n}{m^2})^{1/3}\} \asymp \frac{\sigma_\xi^2}{n}\min\{m,(nm)^{1/3}\}.
\end{align*}
In summary, one may compute that
\[
r_\star^2\asymp\frac{\sigma_\xi^2}{n}\Big[\log N+\min\{m,(nm)^{1/3}\}\Big].
\]
Item~\emph{2.} follows from $\frac{\log(e(m+1))+\log\log(16N)}{\log N+\min\{m,(nm)^{1/3}\}}\to0,$ using $\min\{m,(nm)^{1/3}\}\ge m^{1/3}$.

\end{enumerate}

\subsection{Computation of examples in Section~\texorpdfstring{\ref{sec:examples_LSE}}{constrained LSE}}\label{sec:computation_LSE_examples}

\paragraph{Computation of Example~\texorpdfstring{\ref{example:ellipsoid_rate}}{lq rate}.} Notice that for any $s>0$, $\bE|\xi_j|^s=\frac{2^{s/2}\Gamma(\frac{s+1}{2})}{\sqrt{\pi}}\sigma_\xi^s$,
By \eqref{eq:estimate_ell_q_ellipsoid} of Lemma~\ref{lem:ell_q_coordinate_estimates} applied to $\sigma_j = j^{-\alpha}$, $\bE\min\{A,B\}\leq\min\{\bE[A],\bE[B]\}$, together with the aforementioned bound on $\bE|\xi_j|^s$, we obtain
\begin{align*}
    &\bE\bigl|\eta_{\lambda,j}(\xi_j)\bigr|^q\leq \min\bigl\{\bE|\xi_j|^q, (\lambda j^{\alpha q /2 })^{-\frac{q}{q-1}}\bE|\xi_j|^{\frac{q}{q-1}}\bigr\} \lesssim_q \min\Bigl\{\sigma_\xi^q,\, \bigl( \frac{\sigma_\xi}{\lambda j^{\alpha q/2}} \bigr)^{\frac{q}{q-1}} \Bigr\}.
\end{align*}Notice that $\bP(\sigma_\xi\leq |\xi_j|\leq 2\sigma_\xi) = \frac{2}{\sqrt{2\pi}}\int_1^2\exp(-\frac{1}{2}t^2)\mathrm{d}t\geq \frac{2\exp(-2)}{\sqrt{2\pi}}$. Therefore, applying \eqref{eq:estimate_ell_q_ellipsoid} again, the above displayed inequality is reversed, that is,
\begin{align*}
    \bE\bigl|\eta_{\lambda,j}(\xi_j)\bigr|^q &\geq \bE\Bigl[\Bigl|\eta_{\lambda,j}(\xi_j)\Bigr|^q\one_{\{\sigma_\xi\leq|\xi_j|\leq2\sigma_\xi\}}\Bigr] \geq \min\Bigl\{\frac{1}{2^q}, 2^{-\frac{q}{q-1}}\Bigr\}\bP\Bigl(\sigma_\xi\leq \bigl|\xi_j\bigr|\leq 2\sigma_\xi\Bigr)\min\Bigl\{\sigma_\xi^q, \Bigl(\frac{\sigma_\xi}{\lambda j^{\alpha q/2}}\Bigr)^{\frac{q}{q-1}}\Bigr\}\\
    &\gtrsim_q \min\Bigl\{\sigma_\xi^q, \Bigl(\frac{\sigma_\xi}{\lambda j^{\alpha q/2}}\Bigr)^{\frac{q}{q-1}}\Bigr\}.
\end{align*}
As a result,
\begin{align*}
    &\bE\bigl|\eta_{\lambda,j}(\xi_j)\bigr|^q\asymp_q \min\Bigl\{\sigma_\xi^q,\, \bigl( \frac{\sigma_\xi}{\lambda j^{\alpha q/2}} \bigr)^{\frac{q}{q-1}} \Bigr\}.
\end{align*}We first notice that for any $2\leq (\lambda\sigma_\xi^{q-2})^{-\frac{2}{\alpha q}}\leq N$, separating the sum at $j=\lfloor (\lambda\sigma_\xi^{q-2})^{-\frac{2}{\alpha q}} \rfloor$ yields that
\begin{align*}
    \sum_{j=1}^N j^{\alpha q/2}\bE\bigl|\eta_{\lambda,j}(\xi_j)\bigr|^q &\asymp_q \sigma_\xi^q\sum_{j=1}^{\lfloor (\lambda\sigma_\xi^{q-2})^{-\frac{2}{\alpha q}} \rfloor} j^{\alpha q/2} + \Bigl(\frac{\sigma_\xi}{\lambda}\Bigr)^{\frac{q}{q-1}} \sum_{j=\lfloor (\lambda\sigma_\xi^{q-2})^{-\frac{2}{\alpha q}} \rfloor+1}^N j^{-\frac{\alpha q}{2(q-1)}}.
\end{align*}Notice that $\frac{\sigma_\xi}{\lambda} = \sigma_\xi^{q-1}(\lambda\sigma_\xi^{q-2})^{-1}$.
Using $\sigma_\xi^q\sum_{j=1}^{\lfloor (\lambda\sigma_\xi^{q-2})^{-\frac{2}{\alpha q}} \rfloor} j^{\alpha q/2}\asymp_q \sigma_\xi^q\Bigl(\bigl(\lambda\sigma_\xi^{q-2}\bigr)^{-\frac{2}{\alpha q}}\Bigr)^{1+\frac{\alpha q}{2}}$ together with $\Bigl(\frac{\sigma_\xi}{\lambda}\Bigr)^{\frac{q}{q-1}} \sum_{j=\lfloor (\lambda\sigma_\xi^{q-2})^{-\frac{2}{\alpha q}} \rfloor+1}^N j^{-\frac{\alpha q}{2(q-1)}}\lesssim_{\alpha,q} \bigl(\frac{\sigma_\xi}{\lambda}\bigr)^{\frac{q}{q-1}}\Bigl( \bigl(\lambda \sigma_\xi^{q-2}\bigr)^{-\frac{2}{\alpha q}} \Bigr)^{1-\frac{\alpha q}{2(q-1)}} = \sigma_\xi^q \bigl(\lambda\sigma_\xi^{q-2}\bigr)^{-1-\frac{2}{\alpha q}}$, we obtain
\begin{align*}
    \sum_{j=1}^N j^{\alpha q/2}\bE\bigl|\eta_{\lambda,j}(\xi_j)\bigr|^q &\asymp_q \sigma_\xi^q\bigl(\lambda\sigma_\xi^{q-2}\bigr)^{-1-\frac{2}{\alpha q}}.
\end{align*}
We write $c_0,C_0$ such that
\begin{align}\label{eq:two_sides_estimate_sum_eta_j_xi}
    \forall 2\leq (\lambda\sigma_\xi^{q-2})^{-\frac{2}{\alpha q}}\leq N,\qquad c_0 \sigma_\xi^q\bigl(\lambda\sigma_\xi^{q-2}\bigr)^{-1-\frac{2}{\alpha q}}\leq \sum_{j=1}^N j^{\alpha q/2}\bE\bigl|\eta_{\lambda,j}(\xi_j)\bigr|^q \le C_0\sigma_\xi^q\bigl(\lambda\sigma_\xi^{q-2}\bigr)^{-1-\frac{2}{\alpha q}}.
\end{align}
Let $c_1\leq 1$ such that $c_0c_1^{-1-\frac{2}{\alpha q}}\geq 2$, and let $C_1\geq1$ such that $C_0C_1^{-1-\frac{2}{\alpha q}}\leq (\frac{\delta}{2q})^q$. Let 
\begin{align*}
    C\geq \max\Bigl\{1,\,c_1^{-\frac{2}{\alpha q}},\, \Bigl(\frac{(1+\frac{\alpha q}{2})2^{q+1}\sqrt{\pi}}{2^{q/2}\Gamma(\frac{q+1}{2})}\Bigr)^{\frac{2}{2+\alpha q}},\,\frac{4(\alpha q+2)^2\sqrt{\pi}\Gamma(q+\frac{1}{2})}{\delta\Gamma(\frac{q+1}{2})^2}\Bigr\}.
\end{align*}Let
\begin{align*}
    c\leq \min\Bigl\{2^{-\frac{\alpha q+2}{2q}}C_1^{-\frac{\alpha q+2}{\alpha q^2}},\, \Bigl( \frac{C_1^{-\frac{2}{\alpha q}}}{2(C')^2} \Bigr)^{\frac{\alpha q+2}{2q}}\Bigr\},
\end{align*}where $C'$ is the constant in Theorem~\ref{theo:ell_q}.\

Let $\Omega=\bigl\{\sum_{j=1}^N j^{\alpha q/2}|\xi_j|^q \geq 2^q\bigr\}$. Notice that $\sum_{j=1}^Nj^{\alpha q/2}\geq \frac{1}{1+\alpha q/2}N^{1+\frac{\alpha q}{2}}$. Then by the definition of $C$ and the assumption that $N\geq C\sigma_\xi^{-\frac{2q}{\alpha q+2}}$,
\begin{align*}
    \bE\sum_{j=1}^Nj^{\alpha q/2}\bigl|\xi_j\bigr|^q&\geq\frac{2^{q/2}\Gamma(\frac{q+1}{2})}{\sqrt{\pi}(1+\frac{\alpha q}{2})}\sigma_\xi^qN^{1+\frac{\alpha q}{2}}\geq 2^{q+1}.
\end{align*}
Hence $\frac{1}{2}\bE|\xi_1|^q\sum_{j=1}^N j^{\alpha q/2}\geq 2^q$, and on $\Omega^c$, $\sum_{j=1}^Nj^{\alpha q/2}\bigl|\xi_j\bigr|^q<\frac{1}{2}\bE\sum_{j=1}^N j^{\alpha q/2}\bigl|\xi_j\bigr|^q$, and
\begin{align*}
    \Omega^c\subseteq\Bigl\{ \Bigl| \sum_{j=1}^N j^{\alpha q/2}\bigl|\xi_j\bigr|^q - \bE \sum_{j=1}^N j^{\alpha q/2}\bigl|\xi_j\bigr|^q \Bigr|\geq \frac{1}{2}\bE \sum_{j=1}^N j^{\alpha q/2}\bigl|\xi_j\bigr|^q  \Bigr\}.
\end{align*}
Chebyshev's inequality, $\sum_{j=1}^N j^{\alpha q}\leq N^{1+\alpha q}$, $\sum_{j=1}^N j^{\alpha q/2}\geq \frac{2}{\alpha q+2}N^{1+\alpha q/2}$ together with $\frac{\bE|\xi_1|^{2q}}{(\bE|\xi_1|^q)^2}=\frac{\sqrt{\pi}\Gamma(q+\frac{1}{2})}{\Gamma(\frac{q+1}{2})^2}$ imply
\begin{align*}
    \bP\bigl(\Omega^c\bigr) &\leq \frac{4\Var\bigl(\sum_{j=1}^N j^{\alpha q/2}|\xi_j|^q\bigr)}{\bigl(\bE|\xi_1|^q\bigr)^2\bigl(\sum_{j=1}^N j^{\alpha q/2}\bigr)^2} \leq \frac{4\bE\bigl|\xi_1\bigr|^{2q}\sum_{j=1}^N j^{\alpha q}}{\bigl(\bE|\xi_1|^q\bigr)^2\bigl(\sum_{j=1}^N j^{\alpha q/2}\bigr)^2}\leq \frac{(\alpha q+2)^2}{N}\frac{\sqrt{\pi}\Gamma(q+\frac{1}{2})}{\Gamma(\frac{q+1}{2})^2} \leq \frac{\delta}{4},
\end{align*}where we have used the definition of $C$ in the last inequality. On $\Omega$, by the triangular inequality,
\begin{align*}
    \Bigl(\sum_{j=1}^N j^{\alpha q/2}\bigl|y_j\bigr|^q\Bigr)^{1/q} \geq \Bigl(\sum_{j=1}^N j^{\alpha q/2}\bigl|\xi_j\bigr|^q\Bigr)^{1/q} - \Bigl(\sum_{j=1}^N j^{\alpha q/2}\bigl|\beta_j^\star\bigr|^q\Bigr)^{1/q} \geq 2 - (1-\delta)^{1/q}\geq 1.
\end{align*}Moreover, since $\one_{\Omega}\leq\min\{1,\sum_{j=1}^Nj^{\alpha q/2}|y_j|^q\}$, there holds $1-\frac{\delta}{4}\leq\bP(\Omega)\leq\bE\min\{1,\sum_{j=1}^N j^{\alpha q/2}\bigl|y_j\bigr|^q\}\leq 1$.

Now taking $\lambda = c_1(\sigma_\xi^2)^{\frac{\alpha q+2-q}{\alpha q+2}}$, then $2\leq \bigl(\lambda\sigma_\xi^{q-2}\bigr)^{-\frac{2}{\alpha q}} = c_1^{-\frac{2}{\alpha q}}\sigma_\xi^{-\frac{2q}{\alpha q+2}}\leq N$, where we have used that the assumption that $\sigma_\xi\leq c$; the definition of $c$ such that $c\leq 2^{-\frac{\alpha q+2}{2q}}C_1^{-\frac{\alpha q+2}{\alpha q^2}}$; the assumption that $N\geq C\sigma_\xi^{-\frac{2q}{\alpha q+2}}$; and the definition that $C\geq c_1^{-\frac{2}{\alpha q}}$. Notice that $|t|\mapsto |\eta_{\lambda,j}(t)|$ is monotonically increasing; and notice that $\bP(|\beta_j^\star+\xi_j|>t)\geq\bP(|\xi_j|>t)$ for any $t\geq 0$. Then
\begin{align*}
    \sum_{j=1}^N j^{\alpha q/2}\bE\bigl|\eta_{\lambda,j}(y_j)\bigr|^q\geq \sum_{j=1}^Nj^{\alpha q/2}\bE\bigl|\eta_{\lambda,j}(\xi_j)\bigr|^q\geq c_0\sigma_\xi^q\bigl(\lambda\sigma_\xi^{q-2}\bigr)^{-1-\frac{2}{\alpha q}}\geq c_0c_1^{-1-\frac{2}{\alpha q}}\geq 2\geq \bE\min\{1,\sum_{j=1}^N j^{\alpha q/2}\bigl|y_j\bigr|^q\},
\end{align*}where we have used the lower side of \eqref{eq:two_sides_estimate_sum_eta_j_xi}, and the definition that $c_0c_1^{-1-\frac{2}{\alpha q}}\geq 2$.

Taking $\lambda = C_1\bigl(\sigma_\xi^2\bigr)^{\frac{\alpha q + 2-q}{\alpha q +2}}$ then $2\leq \bigl(\lambda\sigma_\xi^{q-2}\bigr)^{-\frac{2}{\alpha q}}\leq N$, where we have used the assumption that $\sigma_\xi\leq c$; the definition of $c$ such that $c\leq 2^{-\frac{\alpha q+2}{2q}}C_1^{-\frac{\alpha q+2}{\alpha q^2}}$; the assumption that $N\geq C\sigma_\xi^{-\frac{2q}{\alpha q+2}}$.  By triangular inequality the upper side of \eqref{eq:two_sides_estimate_sum_eta_j_xi}, and the definition of $C_1$ such that $C_0C_1^{-1-\frac{2}{\alpha q}}\leq (\frac{\delta}{2q})^q$, we obtain
\begin{align*}
    \sum_{j=1}^N j^{\alpha q/2}\bE\bigl|\eta_{\lambda,j}(y_j)\bigr|^q &\leq \bigg(\Bigl(\sum_{j=1}^N j^{\alpha q/2}\bigl|\beta_j^\star\bigr|^q\Bigr)^{1/q} + \Bigl(\sum_{j=1}^N j^{\alpha q/2}\bE\bigl|\eta_{\lambda,j}(\xi_j)\bigr|^q\Bigr)^{1/q}\biggr)^q \leq \biggl((1-\delta)^{1/q} + \frac{\delta}{2q}\biggr)^q \\
    &\leq 1-\frac{1}{2}\delta\leq 1-\frac{\delta}{4}\leq \bE\min\Bigl\{1,\sum_{j=1}^Nj^{\alpha q/2}\bigl|y_j\bigr|^q\Bigr\}.
\end{align*}As a result, $\lambda_\self\asymp (\sigma_\xi^2)^{\frac{\alpha q+2-q}{\alpha q+2}}$. We now compute $r_\star^2$. Similar to the proof of \eqref{eq:two_sides_estimate_sum_eta_j_xi}, that is, separating the sum at $j=\lfloor\bigl(\lambda_\self\sigma_\xi^{q-2}\bigr)^{-\frac{2}{\alpha q}}\rfloor$, one may compute that the variance term 
\begin{align*}
    \sum_{j=1}^N \min\{\sigma_\xi^2, \bigl(\frac{\sigma_\xi}{\lambda_\self j^{\alpha q/2}}\bigr)^{\frac{2}{q-1}}\}\asymp_{\alpha,q} \sigma_\xi^2\bigl(\lambda_\self \sigma_\xi^{q-2}\bigr)^{-\frac{2}{\alpha q}}\asymp_{\alpha,q,\delta} (\sigma_\xi^2)^{\frac{\alpha q+2 -q}{\alpha q+2}}.
\end{align*}
Moreover, by $\min\{a,b\}\leq\sqrt{ab}$ for $a,b\geq0$, we obtain $\sum_{j=1}^N\min\{|\beta_j^\star|^2, \lambda_\self^2j^{\alpha q}|\beta_j^\star|^{2q-2}\}\leq \sum_{j=1}^N\sqrt{\bigl|\beta_j^\star\bigr|^2\lambda_\self^2j^{\alpha q}\bigl|\beta_j^\star\bigr|^{2q-2}} = \lambda_\self\sum_{j=1}^Nj^{\alpha q/2}|\beta_j^\star|^q\leq\lambda_\self$. Therefore $r_\star^2\asymp \lambda_\self$. It is easy to see that $\sigma_\xi^{-2}r_\star^2\geq \frac{1}{2}C_1^{-\frac{2}{\alpha q}}\sigma_\xi^{-\frac{2q}{\alpha q+2}}$. By the assumption that $\sigma_\xi\leq c$, and by the definition of $c$, there holds $\sigma_\xi^{-2}r_\star^2\geq (C')^2$, and hence Theorem~\ref{theo:ell_q} applies.

\paragraph{Computation of Example~\texorpdfstring{\ref{example:k_interpolated}}{interpolated body}.}For any $t\geq0$, $s\in\{1,2\}$, one may compute that
\begin{align}
    \bE\Bigl[\min\Bigl\{1,\bigl(\bigl|\xi_1\bigr|-t\bigr)_+\Bigr\}^s\Bigr] &= \begin{cases}
        \sqrt{\frac{2}{\pi}}\frac{\exp(-\frac{1}{2\sigma_\xi^2}t^2)}{\sigma_\xi}\int_0^\infty\min\{1,u\}^s\exp\bigl( -\frac{tu}{\sigma_\xi^2} - \frac{u^2}{2\sigma_\xi^2} \bigr)\mathrm{d}u, &\mbox{ if }\xi\mbox{ is Gaussian},\\
        \frac{\sqrt{2}}{\sigma_\xi}\exp\bigl(-\frac{\sqrt2}{\sigma_\xi}t\bigr)\int_0^\infty \min\{1,u\}^s\exp\bigl(-\frac{\sqrt{2}u}{\sigma_\xi}\bigr)\mathrm{d}u, & \mbox{ if }\xi\mbox{ is exponential}.
    \end{cases}\notag\\
    &=\begin{cases}
        \frac{\exp\bigl(-\frac{1}{2\sigma_\xi^2}t^2\bigr)}{1+\frac{t}{\sigma_\xi}}\min\bigl\{1,\frac{\sigma_\xi}{1+\frac{t}{\sigma_\xi}}\bigr\}^s,&\mbox{ if Gaussian},\\
        \exp\bigl(-\frac{\sqrt{2}t}{\sigma_\xi}\bigr)\min\bigl\{\sigma_\xi,1\bigr\}^s,&\mbox{ if exponential}.
    \end{cases}\label{eq:first_order_moment_estimate}
\end{align}Take $t=0$, $s=1$, then together with $\bE|\xi_1|\leq\sigma_\xi$, in both cases we have $\min\{\sigma_\xi,1\}\asymp \bE\min\{1,|\xi_1|\}$.

By $\min\{1,t\}\geq 1-\exp(-t)$ for any $t\geq0$; $\exp(-\frac{u}{k})\leq 1-\frac{u}{2k}$, $0\leq u\leq 1$; $(1-t)\leq \exp(-Nt)$; $\bE\{1,|\xi_1|\}\gtrsim \min\{\sigma_\xi,1\}$; and $1-\exp(-t)\geq\frac12\min\{1,t\}$, we have
\begin{align*}
    &\bE\min\bigl\{k,\sum_{j=1}^N\min\{1,|\xi_j|\}\bigr\}  = k\bE\min\Bigl\{ 1, \frac{1}{k}\sum_{j=1}^N\min\{1,\bigl|\xi_j\bigr|\}\Bigr\} \geq k\Bigl( 1-\bE\exp\bigl(-\frac{1}{k}\sum_{j=1}^N\min\{1,|\xi_j|\}\bigr)\Bigr)\\
    &=k\Bigl( 1-\prod_{j=1}^N\bE\exp\bigl( -\frac{\min\{1,|\xi_j|\}}{k} \bigr) \Bigr) \geq k\Bigl( 1-\prod_{j=1}^N \Bigl( 1-\frac{\bE\min\{1,|\xi_j|\}}{2k}\Bigr)\Bigr)\geq k\Bigl( 1 - \exp\Bigl( -\frac{N\bE\min\{1,|\xi_1|\}}{2k} \Bigr) \Bigr) \\
    &\geq \frac{1}{2}\min\Bigl\{k,\frac{N}{2}\bE\min\{1,|\xi_1|\}\Bigr\}\geq c\min\bigl\{k,N\min\{\sigma_\xi,1\}\bigr\},
\end{align*}for some absolute constant $c$. Moreover, by Jensen's inequality applied to the function $t\in\bR_+\mapsto\min\{k,t\}$, this inequality can be reversed, that is,
\begin{align*}
    \bE\min\Bigl\{k,\sum_{j=1}^N\min\{1,|\xi_j|\}\Bigr\}\leq\min\Bigl\{k,\bE\sum_{j=1}^N\min\bigl\{1,|\xi_j|\bigr\}\Bigr\} \leq \min\{k,N\min\{1,\sigma_\xi\}\}.
\end{align*}Therefore
\begin{align}\label{eq:right_hand_side_k_interpolated}
    \bE\min\bigl\{k,\sum_{j=1}^N\min\{1,|\xi_j|\}\bigr\} \in \Bigl[c,1\Bigr]\min\bigl\{k,N\min\{\sigma_\xi,1\}\bigr\}.
\end{align}
Since $\bbeta^\star=\vzero$, \eqref{eq:eta_self_ell_1_ball} and $r_\star^2$ are equivalent to
\begin{align}\label{eq:simplified_equation}
    N\bE\min\bigl\{1,(|\xi_1|-\lambda_\self)_+\bigr\} = \bE\min\Bigl\{k,\sum_{j=1}^N\bigl\{1,|\xi_j|\bigr\}\Bigr\},\mbox{ and }r_\star^2 = N\bE\bigl[\min\{1,\bigl(|\xi_1|-\lambda_\self\bigr)_+\}^2\bigr],
\end{align}where we have used the fact that $(|\xi_1|-\lambda_\self)_+^2\leq |\xi_1|^2$ almost surely.

Next we estimate $\lambda_\self$.
\begin{enumerate}
    \item When $\xi$ is Gaussian. By \eqref{eq:first_order_moment_estimate} applied to $s=1$, together with $\frac{1}{1+\frac{\lambda_\self}{\sigma_\xi}}\leq\frac{\min\{1,\frac{\sigma_\xi}{1+\lambda_\self/\sigma_\xi}\}}{\min\{\sigma_\xi,1\}}\leq 1$, from \eqref{eq:simplified_equation} we obtain
    \begin{align*}
        \frac{\exp\bigl(-\frac{\lambda_\self^2}{2\sigma_\xi^2}\bigr)}{\bigl(1+\frac{\lambda_\self}{\sigma_\xi}\bigr)^2} \lesssim \min\Bigl\{1,\frac{k}{N\min\{\sigma_\xi,1\}}\Bigr\} \lesssim \frac{\exp\bigl(-\frac{\lambda_\self^2}{2\sigma_\xi^2}\bigr)}{\bigl(1+\frac{\lambda_\self}{\sigma_\xi}\bigr)}.
    \end{align*}Taking logarithm on both sides, we obtain
    \begin{align*}
        1+\frac{\lambda_\self^2}{\sigma_\xi^2} \asymp \log\Bigl(e+\frac{N\min\{\sigma_\xi,1\}}{k}\Bigr)\quad \Rightarrow\quad 1+\frac{\lambda_\self}{\sigma_\xi} \asymp \sqrt{\log\Bigl(e+\frac{N\min\{\sigma_\xi,1\}}{k}\Bigr)},
    \end{align*}where we have used the fact that $\sqrt{1+x^2}\asymp 1+x$ for $x\geq 0$.
    Applying \eqref{eq:first_order_moment_estimate} to $s=2$,
    \begin{align*}
        r_\star^2&\asymp \frac{N\exp\bigl(-\frac{\lambda_\self^2}{2\sigma_\xi^2}\bigr)}{1+\frac{\lambda_\self}{\sigma_\xi}}\min\Bigl\{1,\frac{\sigma_\xi}{1+\frac{\lambda_\self}{\sigma_\xi}}\Bigr\}^2 \asymp \min\bigl\{k,N\min\{\sigma_\xi,1\}\bigr\}\min\Bigl\{1,\frac{\sigma_\xi}{1+\frac{\lambda_\self}{\sigma_\xi}}\Bigr\}\\
        &\asymp \min\bigl\{k,N\min\{\sigma_\xi,1\}\bigr\}\min\Bigl\{1,\frac{\sigma_\xi}{\sqrt{\log\Bigl(e+\frac{N\min\{\sigma_\xi,1\}}{k}\Bigr)}}\Bigr\}.
    \end{align*}
    \item When $\xi$ is exponential. Similarly by \eqref{eq:first_order_moment_estimate}, applied to $s=1$, from \eqref{eq:simplified_equation} and \eqref{eq:right_hand_side_k_interpolated} we obtain 
    \begin{align*}
        \exp\Bigl(-\frac{\sqrt{2}\lambda_\self}{\sigma_\xi}\Bigr)\asymp \min\Bigl\{1,\frac{k}{N\min\{\sigma_\xi,1\}}\Bigr\}\quad\Rightarrow\quad 1+\frac{\lambda_{\self}}{\sigma_\xi}\asymp \log\Bigl( e+\frac{N\min\{\sigma_\xi,1\}}{k} \Bigr).
    \end{align*}Similarly applying \eqref{eq:first_order_moment_estimate} to $s=2$, we obtain $r_\star^2\asymp \min\{N\sigma_\xi^2,k\sigma_\xi,k\}$.
\end{enumerate}
\bibliographystyle{alpha}
\bibliography{other}
\end{document}